\documentclass[reqno]{amsart}
\usepackage[margin=1in]{geometry}
\usepackage[x11names,dvipsnames]{xcolor}
\usepackage{amsmath,amssymb,amsthm}
\usepackage{lmodern}
\usepackage{microtype}
\usepackage{enumitem}
\setlist[enumerate]{label=\upshape(\arabic*)}
\usepackage{tikz}
\usepackage{float}
\usetikzlibrary{arrows.meta,calc,decorations.markings}
\usepackage[colorlinks]{hyperref}
\hypersetup{
  colorlinks=true,
  linkcolor=,
  citecolor=,
  filecolor=,
  urlcolor=
}
\definecolor{mtone}{RGB}{0,114,178}
\definecolor{mttwo}{RGB}{213,94,0}
\definecolor{mtthree}{RGB}{0,140,100}
\definecolor{mtstab}{RGB}{150,70,130}
\definecolor{mtslice}{RGB}{160,110,0}
\numberwithin{equation}{section}

\newtheorem{lemma}[equation]{Lemma}
\newtheorem{proposition}[equation]{Proposition}
\newtheorem{corollary}[equation]{Corollary}
\theoremstyle{definition}
\newtheorem{definition}[equation]{Definition}
\newtheorem{example}[equation]{Example}
\theoremstyle{remark}
\newtheorem{remark}[equation]{Remark}
\theoremstyle{plain}
\newtheorem{introtheorem}{Theorem}

\tikzset{
  vx/.style={circle,inner sep=1.5pt},
  vname/.style={font=\footnotesize,text=black!55},
  midarrow/.style={postaction={decorate},decoration={markings,mark=at position 0.55 with {\arrow{Stealth[length=2.1mm]}}}}
}

\newcommand{\kk}{\Bbbk}
\DeclareMathOperator{\supp}{Supp}
\DeclareMathOperator{\degt}{deg_t}
\DeclareMathOperator{\lc}{lc}
\DeclareMathOperator{\Ext}{Ext}
\DeclareMathOperator{\Tor}{Tor}
\DeclareMathOperator{\Hom}{Hom}
\DeclareMathOperator{\pd}{pd}
\DeclareMathOperator{\injdim}{injdim}
\DeclareMathOperator{\gldim}{gldim}
\DeclareMathOperator{\GKdim}{GKdim}
\DeclareMathOperator{\hdet}{hdet}
\DeclareMathOperator{\Tr}{Tr}
\DeclareMathOperator{\tr}{tr}
\DeclareMathOperator{\gci}{gci}

\DeclareMathOperator{\Tot}{Tot}

\title{Dual reflection groups of graded Ore extensions}

\author{Yassine Ait Mohamed}
\address{D\'epartement de math\'ematiques, Universit\'e de Sherbrooke, 2500 Bd de l'Universit\'e, Sherbrooke, QC, J1K 2R1, Canada}
\email{yassine.ait.mohamed@usherbrooke.ca}

\author{Asmae El-Baz}
\address{D\'epartement de math\'ematiques et de statistique, Universit\'e Laval, Qu\'ebec, QC, G1V 0A6, Canada}
\email{aselb@ulaval.ca}

\begin{document}

\begin{abstract}
We prove that the Artin--Schelter Gorenstein and Artin--Schelter regular properties descend from a graded Ore extension to its coefficient algebra, for every twisting automorphism and every skew derivation. We then determine the dual reflection groups of graded Ore extensions for all gradings compatible with the Ore structure. They are exactly the semidirect products of a dual reflection group of the coefficient algebra by the cyclic group generated by the degree of the new variable, meeting it trivially, and a nonzero skew derivation never contributes a new degree. For iterated Ore extensions we characterize the dual reflection groups by an order formula analogous to the theorem of Shephard, Todd and Chevalley, and we prove the binomial conjecture of Buzaglo and Rogalski for these algebras. We show that their covariant rings are complete intersections of growth type, whose Ext-algebras grow like a polynomial ring in as many variables as there are generators of nontrivial degree. This answers a question of Kirkman, Kuzmanovich and Zhang for these algebras. Finally, we describe the smash product as a Cayley quiver with relations, factor its matrix Hilbert series, and show that the extension becomes Hopf--Galois once the Jacobian is inverted.
\end{abstract}

\maketitle



\section{Introduction}\label{YAMs0Ik2}

Let $\kk$ be a field of characteristic zero, let $V$ be a finite-dimensional $\kk$-vector space, and let $W\subseteq GL(V)$ be a finite group acting on the polynomial ring $\kk[V]$. The theorem of Shephard, Todd and Chevalley \cite{ST,Che} states that the invariant ring $\kk[V]^W$ is again a polynomial ring if and only if $W$ is generated by pseudo-reflections. In this case the invariant ring is generated by homogeneous elements $u_1,\dots,u_n$ of degrees $d_1,\dots,d_n$, and
\begin{equation}\label{YAMq1Ab3}
|W|=d_1d_2\cdots d_n,\quad H_{\kk[V]^W}(s)=\prod_{i=1}^n\frac1{1-s^{d_i}},
\end{equation}
where $H$ denotes the Hilbert series. Stanley described the relative invariants of such a group \cite[Theorems 2.3 and 3.1]{Sta}. For every linear character $\chi$ of $W$, the space of $\chi$-relative invariants is a free $\kk[V]^W$-module of rank one,
\begin{equation}\label{YAMq2Bc4}
\kk[V]^W_\chi=f_\chi\,\kk[V]^W,
\end{equation}
generated by a product $f_\chi$ of linear forms defining reflecting hyperplanes. The Jacobian $\det(\partial u_j/\partial x_i)$ is a relative invariant. Its character is the inverse of the determinant of the action of $W$ on the linear forms, and up to a nonzero scalar it is the corresponding generator $f_\chi$. The group $W$ acts freely on the complement of its reflecting hyperplanes \cite[Theorem 1.5]{Ste}. Geometrically, the quotient map $V\to V/W$ is flat, its fibre over the origin is the spectrum of the coinvariant algebra $\kk[V]/(u_1,\dots,u_n)$, which is a complete intersection, and $W$ acts freely away from the reflecting hyperplanes. The results of this paper are noncommutative analogues of these three facts.

Artin--Schelter (AS) regular algebras are the noncommutative analogues of polynomial rings \cite{AS}. Their invariant theory studies actions of finite groups, and more generally of semisimple Hopf algebras, on AS regular algebras. In this setting the determinant is replaced by the homological determinant of J{\o}rgensen and Zhang \cite{JZ}, and the theorem of Shephard, Todd and Chevalley has analogues for skew polynomial rings \cite{KKZ0}. Kirkman, Kuzmanovich and Zhang studied the dual situation, in which a group algebra $\kk G$ coacts on an AS regular algebra $B$ \cite{KKZ}. A coaction of $\kk G$ on $B$ is the same as a $G$-grading
\[
B=\bigoplus_{g\in G}B_g,\quad B_gB_h\subseteq B_{gh},
\]
and its coinvariants form the identity component $B_e$. Let $B$ be a noetherian AS regular domain generated in degree one. A finite group $G$ is a \emph{dual reflection group} for $B$ if $B$ has a faithful $G$-grading, refining its $\mathbb N$-grading, whose identity component $B_e$ is again AS regular \cite[Definition 0.1]{KKZ}.

The structure theory of reflection groups extends to dual reflection groups. If $G$ is a dual reflection group for $B$, every component $B_g$ is a free $B_e$-module of rank one on both sides,
\begin{equation}\label{YAMq4De6}
B_g=f_gB_e=B_ef_g,
\end{equation}
which is the analogue of \eqref{YAMq2Bc4}. There is a unique element $m_G\in G$, the \emph{mass element}, at which $\deg f_g$ is maximal, and $f_{m_G}$ is the Jacobian of the coaction in the sense of Kirkman and Zhang \cite[Definition 2.1]{KZ}. The \emph{covariant ring} $C(B)=B/B(B_e)_{\ge1}$, where $(B_e)_{\ge1}$ is the augmentation ideal of $B_e$, is a Frobenius algebra of dimension $|G|$, with basis the images of the elements $f_g$, and its Hilbert series is a palindromic product of cyclotomic polynomials \cite[Theorems 0.3, 0.4 and 3.5]{KKZ}. In the commutative case the covariant ring is a complete intersection, and Kirkman, Kuzmanovich and Zhang asked whether $C(B)$ is always a complete intersection of growth type, in the sense that its Ext-algebra has finite Gelfand--Kirillov dimension \cite[Question 0.7]{KKZ}. They noted that the question is open even for coactions of group algebras.

In low dimension, dual reflection groups are known through classification. In dimension two, Crawford classified the gradings of AS regular algebras by nonabelian finite groups and determined when the identity component is AS regular \cite{Cra1}. Buzaglo and Rogalski classified the coactions of cocommutative Hopf algebras on one-parameter skew polynomial rings in two variables, and in three variables for most values of the parameter \cite{BR2}, and the nonabelian dual reflection groups of three-dimensional AS regular algebras \cite[Theorem 1.2]{BR}. They conjectured that an AS regular algebra with a dual reflection group whose identity component has no elements of degree one can be presented by binomial relations \cite[Conjecture 1.4]{BR}, and they proved it for nonabelian gradings of three-dimensional AS regular algebras \cite[Corollary 5.15]{BR}. In dimension four, Goetz, Kirkman, Moore and Vashaw constructed families of dual reflection groups, among them groups of wreath product type coacting on skew polynomial extensions \cite[Examples 2.1.6 and 2.1.7]{GKMV}. More recently, Crawford constructed dual reflection groups from solutions to the Yang--Baxter equation \cite{Cra2}.

Ore extensions are noncommutative polynomial rings in one variable \cite{Ore}. If $A$ is AS regular and $\sigma$ is a graded automorphism of $A$, then every graded Ore extension $B=A[t;\sigma,\delta]$ is again AS regular by the graded version of \cite[Theorem 3.3]{LWW}. Many of the algebras in the classifications above are iterated Ore extensions, and one wants to obtain the dual reflection groups of $B$ from those of $A$. When $\delta=0$, Buzaglo and Rogalski proved that a dual reflection group of $B$ for which $A$ and $t$ are homogeneous is a semidirect product built from $A$ \cite[Proposition 6.1]{BR}. In dimension three, the Ore extensions with $\delta\neq0$ are absent from their list of algebras with nonabelian dual reflection groups \cite[Theorem 1.2]{BR}. This absence is established through the classification, and L. Buzaglo asked whether there is a structural reason behind it. This paper gives such a reason for gradings compatible with the Ore structure. We then classify the dual reflection groups of iterated Ore extensions with homogeneous Ore variables. We describe their covariant rings and determine the bigraded Poincar\'e series and the growth of the Ext-algebras of these rings. We also compute the associated McKay quivers and localizations.

Throughout, $\kk$ is a field. In the structure theory of dual reflection groups recalled at the end of Section \ref{YAMs1Pr0}, and from Section \ref{YAMs4Or0} on, $\kk$ is algebraically closed of characteristic zero. Gradings by groups are faithful and refine the $\mathbb N$-grading. If $b\in B_g$ is nonzero for a $G$-grading of $B$, we call $g$ the \emph{$G$-degree} of $b$. The \emph{support} $\supp(B)$ is the set of $G$-degrees of nonzero homogeneous elements of $B$. The symbol $\deg$ always denotes the $\mathbb N$-degree. The identity element of a group is $e$, and the order of an element $\gamma$ is $|\gamma|$. We write $[r]_s=1+s+\cdots+s^{r-1}$.

Our first result is a descent theorem for the Gorenstein property. Let $R$ be a noetherian connected graded algebra, let $\sigma$ be a graded automorphism of $R$, and let $\delta$ be a $\sigma$-derivation of $R$ such that
\begin{equation}\label{YAMq5Ef7}
\delta(R_n)\subseteq R_{n+m}\quad\text{for all }n\ge0,
\end{equation}
for a fixed integer $m\ge1$. The Ore extension $S=R[t;\sigma,\delta]$ is then connected graded with $\deg t=m$. The AS Gorenstein algebras and their Gorenstein parameters are recalled in Definition \ref{YAMd2As2}.

\begin{introtheorem}\label{YAMa7Qd2}
Suppose that $S=R[t;\sigma,\delta]$ is AS Gorenstein of dimension $D$ with Gorenstein parameter $\ell$. Then $R$ is AS Gorenstein of dimension $D-1$ with Gorenstein parameter $\ell-m$. If moreover $S$ is AS regular, then $R$ is AS regular.
\end{introtheorem}

Together with the preservation theorem of Liu, Wang and Wu, Theorem \ref{YAMa7Qd2} shows that $S$ is AS regular if and only if $R$ is AS regular. Zhang and Zhang showed, in a setting that includes graded Ore extensions, that the Gorenstein condition passes from $R$ to $S$ \cite[Proposition 3.2]{ZZ}, and Theorem \ref{YAMa7Qd2} gives the converse. When $\delta=0$, the variable $t$ is a normal regular element and $R\cong S/tS$, so the Gorenstein part of the theorem follows from the Rees lemma \cite[Proposition 3.4(b)]{Lev}, and the regular part also uses the descent of the global dimension \cite[\S7.5]{MR}. The new case is $\delta\neq0$. Then $t$ is in general not normal, as for the Jordan plane, and the finiteness of the global dimension and the Hilbert series still descend by elementary arguments, but the Gorenstein condition requires another argument, which we obtain from local duality for $S$. We are not aware of a previous statement of this descent for $\delta\neq0$, either in the Gorenstein or in the regular case. Theorem \ref{YAMa7Qd2} is the Gorenstein counterpart of Phan's theorem that Koszulity descends along graded Ore extensions with $\deg t=1$ \cite[Corollary 1.3]{Pha}. It holds for every degree $m$ of the new variable and for AS Gorenstein algebras that are not regular.

We fix the setting for dual reflection groups. Let $A$ be a noetherian AS regular domain generated in degree one, let $\sigma$ be a graded automorphism of $A$, let $\delta$ be a $\sigma$-derivation of $A$ with $\delta(A_n)\subseteq A_{n+1}$ for all $n$, and put
\[
B=A[t;\sigma,\delta],\quad\deg t=1.
\]
A grading of $B$ by a group $G$ is \emph{compatible with the Ore structure} if $A$ is a graded subalgebra and $t$ is homogeneous. We studied gradings of Ore extensions by arbitrary groups in \cite{AM1} for differential polynomial rings and in \cite{AM2} for general Ore extensions. By \cite[Theorem A]{AM2}, a $G$-grading of $A$ extends to a grading of $B$ in which $t$ is homogeneous of $G$-degree $\gamma$ if and only if
\[
\sigma(A_g)\subseteq A_{\gamma g\gamma^{-1}}\quad\text{and}\quad\delta(A_g)\subseteq A_{\gamma g}\quad\text{for all }g\in G.
\]
When $\sigma=\mathrm{id}$, the degree $\gamma$ commutes with the support of $A$ \cite[Remark 2.2(2)]{AM2}. For general $\sigma$ this fails, and this is the source of the semidirect products below. Our second result determines the dual reflection groups of $B$.

\begin{introtheorem}\label{YAMb4Lx9}
Let $G$ be a finite group with a grading of $B=A[t;\sigma,\delta]$ compatible with the Ore structure. Let $\gamma\in G$ be the $G$-degree of $t$, let $K=\supp(A)$ and let $r=|\gamma|$. Then $G$ is a dual reflection group for $B$ if and only if $K$ is a dual reflection group for $A$, for the grading of $A$ by $K$, and
\[
K\cap\langle\gamma\rangle=\{e\}.
\]
In this case $G=K\rtimes\langle\gamma\rangle$, the identity component is the graded Ore extension
\[
B_e=A_e[t^r;\sigma^r|_{A_e},\delta|_{A_e}],
\]
where $\delta=0$ whenever $r>1$.

If $G$ is a dual reflection group for $B$ and $\delta\neq0$, then $\gamma=e$. In this case $\sigma$ and $\delta$ preserve every component $A_g$, and $(A_e)_1\neq0$.
\end{introtheorem}

When $\delta=0$, Theorem \ref{YAMb4Lx9} is \cite[Proposition 6.1]{BR}. The case $\delta\neq0$ is new, and the criterion is the same for every skew derivation. In particular a nonzero skew derivation never contributes a new $G$-degree to a dual reflection group, and the identity component of a dual reflection group of an Ore extension with $\delta\neq0$ always has elements of degree one. If moreover $\dim A=2$, then $A_1$ is spanned by two homogeneous elements, one of which lies in $A_e$, so $G=K$ is cyclic. This explains the absence of Ore extensions with $\delta\neq0$ observed by Buzaglo and Rogalski in dimension three. Theorem \ref{YAMb4Lx9} also explains the wreath product examples of Goetz, Kirkman, Moore and Vashaw (Example \ref{YAMex7Wr3}).

Applying Theorem \ref{YAMb4Lx9} to each step of an iterated Ore extension determines the group completely. Let
\begin{equation}\label{YAMq8Hi1}
B=\kk[x_1][x_2;\sigma_2,\delta_2]\cdots[x_n;\sigma_n,\delta_n],\quad\deg x_i=1,
\end{equation}
where $\sigma_i$ is a graded automorphism of the subalgebra $B^{(i-1)}$ generated by $x_1,\dots,x_{i-1}$ and $\delta_i$ is a $\sigma_i$-derivation of $B^{(i-1)}$ raising degrees by one. Every grading of $B$ in which each $x_i$ is homogeneous is compatible with each step of \eqref{YAMq8Hi1}. For such a grading we write $c_i\in G$ for the $G$-degree of $x_i$, and $r_i=|c_i|$.

\begin{introtheorem}\label{YAMtc3Mn8}
Let $G$ be a finite group with a grading of the iterated Ore extension $B$ of \eqref{YAMq8Hi1} in which every $x_i$ is homogeneous. Then $G=\langle c_1\rangle\langle c_2\rangle\cdots\langle c_n\rangle$, and the following are equivalent.
\begin{enumerate}
\item $G$ is a dual reflection group for $B$.
\item $|G|=|c_1|\,|c_2|\cdots|c_n|$.
\item Every $g\in G$ has a unique expression $g=c_1^{a_1}c_2^{a_2}\cdots c_n^{a_n}$ with $0\le a_i<|c_i|$.
\end{enumerate}
In this case $G$ is solvable, $c_i=e$ whenever $\delta_i\neq0$, and
\[
\begin{aligned}
f_{c_1^{a_1}\cdots c_n^{a_n}}&=x_1^{a_1}\cdots x_n^{a_n}\quad(0\le a_i<r_i),\\
f_{m_G}&=x_1^{r_1-1}\cdots x_n^{r_n-1},\\
H_{B_e}(s)&=\prod_{i=1}^n\frac1{1-s^{r_i}}.
\end{aligned}
\]
If moreover $(B_e)_1=0$, then every $\delta_i$ is zero, the degrees $c_1,\dots,c_n$ are pairwise distinct, and $B$ is presented by the binomial relations
\[
x_ix_j=\lambda_{ij}\,x_{\pi_i(j)}x_i\quad(1\le j<i\le n),
\]
where $\lambda_{ij}\in\kk^\times$, $\pi_i$ is a permutation of $\{1,\dots,i-1\}$ and $c_{\pi_i(j)}=c_ic_jc_i^{-1}$.
\end{introtheorem}

Condition (2) is the analogue of $|W|=d_1\cdots d_n$ in \eqref{YAMq1Ab3}, with the orders $r_i$ playing the role of the degrees of the basic invariants, and the formula for $H_{B_e}$ is the analogue of the Hilbert series of $\kk[V]^W$. The last part proves the conjecture of Buzaglo and Rogalski \cite[Conjecture 1.4]{BR} for iterated Ore extensions with homogeneous Ore variables, in every dimension and for every finite group. The identity component is again an iterated Ore extension, on the generators $x_i^{r_i}$ (Proposition \ref{YAMp7Ie3}).

Condition (2) reduces the search for dual reflection groups to the finite quotients of an explicit finitely presented group. For instance, the dual reflection groups of $\kk[x_1,\dots,x_n][t;\sigma]$, with $\sigma$ a cyclic permutation of the variables, for gradings in which $x_1,\dots,x_n$ and $t$ are homogeneous, are exactly the cyclic groups and the groups $\mathbb Z_r^n\rtimes\mathbb Z_{np}$ of Goetz, Kirkman, Moore and Vashaw (Corollary \ref{YAMc7Wr4}).

Our fourth result determines the bigraded Poincar\'e series of the covariant ring and the growth of its Ext-algebra. For a connected graded algebra $C$ we write $E(C)=\bigoplus_{n\ge0}\Ext^n_C(\kk,\kk)$ for its Ext-algebra and
\[
P_C(y,s)=\sum_{n,j\ge0}\dim_\kk\Ext^n_C(\kk,\kk)_{-j}\;y^ns^j
\]
for its bigraded Poincar\'e series. We write $\gci(C)=\GKdim E(C)$, so that $C$ is a complete intersection of growth type exactly when $\gci(C)$ is finite \cite[Definition 1.5]{KKZ2}.

\begin{introtheorem}\label{YAMtd9Ci4}
In the setting of Theorem \ref{YAMtc3Mn8}, suppose that $G$ is a dual reflection group for $B$. Then $C(B)$ is an iterated truncated skew polynomial ring
\[
C(B)\cong\kk[\bar x_{i_1}]/(\bar x_{i_1}^{\,r_{i_1}})\,[\bar x_{i_2};\bar\sigma_{i_2}]/(\bar x_{i_2}^{\,r_{i_2}})\cdots[\bar x_{i_m};\bar\sigma_{i_m}]/(\bar x_{i_m}^{\,r_{i_m}}),
\]
where $i_1<\cdots<i_m$ are the indices with $c_i\neq e$ and $\bar\sigma_i$ is the automorphism induced by $\sigma_i$ on the covariant ring of $B^{(i-1)}$. Its Poincar\'e series is
\[
P_{C(B)}(y,s)=\prod_{c_i\neq e}\frac{1+ys}{1-y^2s^{r_i}}.
\]
In particular $C(B)$ is a complete intersection of growth type, and
\[
\gci\big(C(B)\big)=\#\{i\mid c_i\neq e\}.
\]
\end{introtheorem}

Theorem \ref{YAMtd9Ci4} answers \cite[Question 0.7]{KKZ} positively for iterated Ore extensions with homogeneous Ore variables, for abelian and nonabelian groups alike. The Poincar\'e series coincides with that of the commutative complete intersection $\bigotimes_i\kk[x_i]/(x_i^{r_i})$, the covariant ring of the diagonal reflection group $\prod_i\mu_{r_i}$, where $\mu_r$ is the group of $r$-th roots of unity. The proof uses an explicit minimal resolution of the trivial module over a truncated skew polynomial ring (Proposition \ref{YAMp3Tr7}).

Our fifth result concerns the smash products and McKay quivers of these gradings. A $G$-grading of $B$ is an action of the dual group algebra $\kk^G$, and the smash product $\Lambda=B\#\kk^G$ is the analogue of the skew group algebra of a group action. It has orthogonal idempotents $e_g$ indexed by $G$, with $e_g\Lambda e_h\cong B_{gh^{-1}}$ \cite[Corollary 1.5]{CoMo}. For $c\in G$ let $\Pi_c$ be the permutation matrix of left multiplication by $c$ on $G$. The \emph{Cayley quiver} of $G$ with respect to $c_1,\dots,c_n$ has vertex set $G$ and an arrow $h\to c_ih$ for every $h\in G$ and every $i$. The matrix Hilbert series of $\Lambda$ is the matrix $H_\Lambda(s)$ of the Hilbert series of the spaces $e_g\Lambda e_h$, and $P_G(s)$ is the matrix with entries $s^{\deg f_{gh^{-1}}}$.

\begin{introtheorem}\label{YAMtd5Qv1}
Let $G$ be a finite group with a grading of the iterated Ore extension $B$ of \eqref{YAMq8Hi1} in which every $x_i$ is homogeneous, and let $\Lambda=B\#\kk^G$.
\begin{enumerate}
\item The algebra $\Lambda$ is the path algebra of the Cayley quiver of $G$ with respect to $c_1,\dots,c_n$, modulo the lifts of the relations of $B$, and
\[
H_\Lambda(s)=(I-s\Pi_{c_1})^{-1}(I-s\Pi_{c_2})^{-1}\cdots(I-s\Pi_{c_n})^{-1}.
\]
\item If $G$ is a dual reflection group for $B$, then $H_\Lambda(s)=H_{B_e}(s)\,P_G(s)$ and
\[
P_G(s)=\prod_{i=1}^n\big(I+s\Pi_{c_i}+s^2\Pi_{c_i}^2+\cdots+s^{r_i-1}\Pi_{c_i}^{r_i-1}\big).
\]
The algebra $C(B)\#\kk^G$ is Frobenius of dimension $|G|^2$, and the socle of its left projective module $\bar\Lambda e_h$ at the vertex $h$ is the simple module at $m_Gh$.
\item If $G$ is a dual reflection group for $B$ and $L=B[f_{m_G}^{-1}]$, then $L$ is a crossed product of $L_e$ by $G$, and
\[
L\#\kk^G\cong M_{|G|}(L_e).
\]
\end{enumerate}
\end{introtheorem}

The factorization of $P_G$ is the matrix form of the identity $H_{C(B)}(s)=\prod_i[r_i]_s$, and part (3) is the analogue of the fact that a reflection group acts freely off its reflecting hyperplanes. The description of $\Lambda$ by a quiver with relations is an instance of covering theory \cite{Gre,CM2,Asa} and of Crawford's description of smash products of Koszul AS regular algebras \cite[Theorem 1.2]{Cra22}. For iterated Ore extensions the relations are the lifts of the Ore relations, and the matrix Hilbert series factors through the generators, also for the nonabelian groups of Theorem \ref{YAMtc3Mn8}. For the grading of $\kk[x_1,x_2][t;\sigma]$ by the dihedral group of order eight, the quiver is the twisted cube of Figure \ref{YAMfg2Dh3}. If $G$ is abelian and $B=\kk[V]$, then $\Lambda$ is the skew group algebra of a diagonal reflection group. When this group is generated by reflections of order two, the quotient of $\Lambda$ by the ideal generated by $e_e$ was studied by Buchweitz, Faber and Ingalls \cite[Theorem 4.17 and Corollary 4.19]{BFI}.

These results have a geometric meaning, which we explain in Section \ref{YAMs9Gq1}. In the language of noncommutative algebraic geometry \cite{AZ}, we regard $B$ and $B_e$ as algebraic analogues of the coordinate rings of an affine space $X$ and of its quotient $X/G$. Then Theorem \ref{YAMtc3Mn8} is a Chevalley--Shephard--Todd theorem for iterated Ore extensions. The extension $B_e\subseteq B$ is free of rank $|G|$, and Theorem \ref{YAMtd9Ci4} says that its base change at the augmentation of $B_e$ is a complete intersection of growth type. After inverting the Jacobian the extension becomes Hopf--Galois, and for iterated Ore extensions this localization is obtained by inverting the generators of nontrivial $G$-degree (Proposition \ref{YAMp12Gl3}). For the dihedral group of order eight, the point scheme of $B$ is a triangle of lines, whereas $B_e$ is a commutative polynomial ring (Example \ref{YAMex21Pt5}).

We also obtain formulas for homological determinants. If $R$ is noetherian and AS regular and $\theta$ is a graded automorphism of $S=R[t;\sigma,\delta]$ with $\theta(R)=R$, then $\theta(t)=\lambda t+a$ with $\lambda\in\kk^\times$, and
\[
\hdet_S(\theta)=\lambda\,\hdet_R(\theta|_R).
\]
When $\delta=0$ and $a=0$ this follows from a result of Jing and Zhang on normal elements \cite[Proposition 2.4]{JiZ}. For a dual reflection group and a graded automorphism $\theta$ of $A$ preserving every component, we prove
\[
\hdet_{A_e}(\theta|_{A_e})=\lambda_{m_G}\,\hdet_A(\theta),\quad\text{where }\theta(f_{m_G})=\lambda_{m_G}f_{m_G}.
\]
For iterated Ore extensions both determinants become explicit products.

All results on dual reflection groups assume that the grading is compatible with the Ore structure. Once $A$ is a graded subalgebra, a homogeneous Ore variable can always be chosen after an affine change $t'=\lambda t+a$, at the cost of changing $\delta$ by an inner $\sigma$-derivation (Proposition \ref{YAMu9Cj4}). For the three-dimensional AS regular algebras of Ore type graded by nonabelian groups, the whole assumption holds up to a graded change of variables by the classification of Buzaglo and Rogalski \cite[Theorem 3.4]{BR}. In dimension four, the dual reflection group of order sixteen of \cite[\S3.1]{GKMV} acts on an algebra whose presentation as a double Ore extension is not homogeneous for the grading, so this example lies outside the scope of the present methods (Remark \ref{YAMv3Ei8}).

\subsection*{Organization of the paper} Section \ref{YAMs1Pr0} recalls AS regular algebras, graded Ore extensions, group gradings and the structure theory of dual reflection groups. Theorems \ref{YAMa7Qd2} to \ref{YAMtd5Qv1} are proved in Sections \ref{YAMs1Pz4} to \ref{YAMs8Mq0}, respectively. Section \ref{YAMs9Gq1} interprets these results in terms of noncommutative quotients, and Section \ref{YAMs7Hd0} contains the formulas for homological determinants.

\subsection*{Acknowledgments}

We thank Lucas Buzaglo for the question that started this work. We also thank the Institut des sciences math\'ematiques (ISM) for financial support.

\section{Preliminaries}\label{YAMs1Pr0}

In this section we fix the homological conventions and recall the structure theory of graded Ore extensions and dual reflection groups.

A $\kk$-algebra $R$ is \emph{connected graded} if $R=\bigoplus_{n\ge0}R_n$ with $R_iR_j\subseteq R_{i+j}$, $R_0=\kk$ and $\dim_\kk R_n<\infty$ for all $n$. We write $R_+=\bigoplus_{n>0}R_n$ for the augmentation ideal and $\kk=R/R_+$ for the trivial module, on the left or on the right. Let $M=\bigoplus_{n\in\mathbb Z}M_n$ be a graded right $R$-module. For $j\in\mathbb Z$, the \emph{shift} $M(j)$ is the graded module with
\[
M(j)_n=M_{j+n}.
\]
If $M$ is bounded below and $\dim M_n<\infty$ for all $n$, its \emph{Hilbert series} is the Laurent series
\[
H_M(s)=\sum_{n\in\mathbb Z}(\dim_\kk M_n)\,s^n,
\]
so that $H_{M(j)}(s)=s^{-j}H_M(s)$. For graded modules $M$ and $N$ we put
\[
\Hom_R(M,N)=\bigoplus_{j\in\mathbb Z}\Hom_{\mathrm{gr}}\big(M,N(j)\big),
\]
where $\Hom_{\mathrm{gr}}$ denotes the degree-preserving homomorphisms, and $\Ext^i_R(M,N)$ are its right derived functors, computed in the category of graded modules. They are graded vector spaces, and $\Ext^i_R(M,N(j))=\Ext^i_R(M,N)(j)$. All homological dimensions are computed in the category of graded modules. We write $\pd$, $\injdim$ and $\gldim$ for the projective, injective and global dimensions, and $\GKdim$ for the Gelfand--Kirillov dimension \cite[Chapter 8]{MR}.

\begin{definition}\label{YAMd2As2}
A connected graded algebra $R$ is \emph{AS Gorenstein of dimension $d$} if $\injdim R_R=\injdim{}_RR=d<\infty$ and there is $\ell\in\mathbb Z$ such that
\[
\Ext^i_R(\kk,R)\cong\begin{cases}\kk(\ell)&i=d,\\0&i\neq d\end{cases},
\]
both for right and for left modules. The integers $d=\dim R$ and $\ell=\ell_R$ are the \emph{dimension} and the \emph{Gorenstein parameter} of $R$. If moreover $\gldim R=d$ and $\GKdim R<\infty$, then $R$ is \emph{AS regular of dimension $d$}.
\end{definition}

The polynomial ring $\kk[x_1,\dots,x_n]$ with $\deg x_i=1$ is AS regular of dimension $n$ with Gorenstein parameter $n$. For $q\in\kk^\times$ the quantum plane $\kk_q[x,y]=\kk\langle x,y\rangle/(yx-qxy)$ is AS regular of dimension two with Gorenstein parameter two.

We use the following standard facts about a right noetherian connected graded algebra $R$. Every finitely generated graded right $R$-module $M$ has a minimal graded free resolution $P_\bullet\to M$ by finitely generated free modules, obtained from the graded Nakayama lemma, and minimality means that the differentials of $P_\bullet\otimes_R\kk$ vanish. Hence $\Tor^R_i(M,\kk)\cong P_i\otimes_R\kk$, and $\pd M=\max\{i\mid\Tor^R_i(M,\kk)\neq0\}\le\pd{}_R\kk$.

The graded Baer criterion shows that $\gldim R$ is the supremum of $\pd R/I$ over the graded right ideals $I$, so $\gldim R=\pd\kk_R=\pd{}_R\kk$. The same criterion gives
\begin{equation}\label{YAMe1Ba1}
\injdim R_R\le n\iff\Ext^{n+1}_R(R/I,R)=0\text{ for every graded right ideal }I.
\end{equation}
If $p=\pd\kk_R<\infty$ and $P_\bullet$ is a minimal resolution of $\kk_R$, then the image of $\Hom_R(P_{p-1},R)\to\Hom_R(P_p,R)$ lies in $R_+\Hom_R(P_p,R)$, so the graded Nakayama lemma gives
\begin{equation}\label{YAMe1Ba2}
\Ext^p_R(\kk,R)\neq0\quad\text{for }p=\pd\kk_R<\infty.
\end{equation}

Let $R$ be a $\kk$-algebra and $\sigma$ a $\kk$-algebra endomorphism of $R$. A \emph{$\sigma$-derivation} of $R$ is a $\kk$-linear map $\delta\colon R\to R$ such that
\begin{equation}\label{YAMe7Sd7}
\delta(ab)=\sigma(a)\delta(b)+\delta(a)b\quad\text{for all }a,b\in R.
\end{equation}
For $u\in R$ the map $\delta_u(a)=ua-\sigma(a)u$ is a $\sigma$-derivation, called \emph{inner}, and a $\sigma$-derivation is \emph{outer} if it is not inner. The \emph{Ore extension} $R[t;\sigma,\delta]$ is the $\kk$-algebra which is a free left $R$-module with basis $\{t^n\}_{n\ge0}$, with multiplication determined by
\begin{equation}\label{YAMe9Or9}
ta=\sigma(a)t+\delta(a)\quad\text{for all }a\in R,
\end{equation}
and we write $R[t;\sigma]$ when $\delta=0$ \cite[\S1.2]{MR}. For $0\neq f=\sum_ia_it^i$ with $a_i\in R$, the \emph{$t$-degree} $\degt f$ is the largest $n$ with $a_n\neq0$, and $\lc(f)=a_n$ is the \emph{leading coefficient}. An induction on $n$ using \eqref{YAMe9Or9} gives, for $a\in R$ and $n\ge1$,
\begin{equation}\label{YAMe10Tn1}
t^na\in\sigma^n(a)t^n+\sum_{0<j<n}Rt^j+\delta^n(a).
\end{equation}
Multiplying \eqref{YAMe10Tn1} on the left by $t$ gives the terms $\sigma^{n+1}(a)t^{n+1}+\delta\sigma^n(a)t^n$, terms in $\sum_{0<j<n+1}Rt^j$, and $t\delta^n(a)=\sigma\delta^n(a)t+\delta^{n+1}(a)$. By \eqref{YAMe10Tn1}, the coefficient of $t^{n+k}$ in $at^n\cdot bt^k$ is $a\sigma^n(b)$, and all other terms have lower $t$-degree. So if $R$ is a domain and $\sigma$ is injective, then
\begin{equation}\label{YAMe29Dt3}
\degt(fg)=\degt f+\degt g,\quad\lc(fg)=\lc(f)\,\sigma^{\degt f}\big(\lc(g)\big),
\end{equation}
for nonzero $f,g$, and $R[t;\sigma,\delta]$ is a domain.

Let now $R$ be connected graded, $\sigma$ a graded automorphism and $\delta$ a $\sigma$-derivation satisfying \eqref{YAMq5Ef7} for some $m\ge1$. The \emph{graded Ore extension} $S=R[t;\sigma,\delta]$ is graded by $\deg t=m$, that is,
\begin{equation}\label{YAMe11Go2}
S_n=\bigoplus_{k\ge0}R_{n-mk}\,t^k.
\end{equation}
For $a\in R_n$, both sides of \eqref{YAMe9Or9} have degree $n+m$, so \eqref{YAMe11Go2} makes $S$ a connected graded algebra, with Hilbert series
\begin{equation}\label{YAMe12Hs3}
H_S(s)=\frac{H_R(s)}{1-s^m}.
\end{equation}
If $R$ is right noetherian, then so is $S$ by the skew Hilbert basis theorem \cite[Theorem 1.2.9]{MR}. Conversely, if $S$ is right noetherian then so is $R$, since $JS=\bigoplus_nJt^n$ and hence $JS\cap R=J$ for every right ideal $J$ of $R$.

Graded Ore extensions preserve AS regularity. A connected graded algebra is AS regular if and only if it is twisted Calabi--Yau in the graded sense \cite[Lemma 1.2]{RRZ}, and twisted Calabi--Yau algebras are preserved by Ore extensions, with the dimension raised by one, by the graded version of \cite[Theorem 3.3]{LWW}, whose proof applies verbatim by \cite[Remark 2.2]{LWW}. Since $S$ is generated by $t$ and generators of $R$, \eqref{YAMe12Hs3} gives $\GKdim S\le\GKdim R+1$. Hence, if $R$ is noetherian and AS regular, then
\begin{equation}\label{YAMt1Lw9}
S=R[t;\sigma,\delta]\text{ is AS regular and }\dim S=\dim R+1.
\end{equation}

An \emph{iterated Ore extension} is an algebra of the form \eqref{YAMq8Hi1}. We put $B^{(0)}=\kk$, $B^{(1)}=\kk[x_1]$, and for $2\le i\le n$ we require that $\sigma_i$ be a graded automorphism of
\[
B^{(i-1)}=\kk[x_1][x_2;\sigma_2,\delta_2]\cdots[x_{i-1};\sigma_{i-1},\delta_{i-1}],
\]
and that $\delta_i$ be a $\sigma_i$-derivation of $B^{(i-1)}$ raising degrees by one. It is convenient to regard $B^{(1)}$ as the Ore extension $\kk[x_1;\sigma_1,\delta_1]$ with $\sigma_1=\mathrm{id}_\kk$ and $\delta_1=0$. By \eqref{YAMt1Lw9} and \eqref{YAMe29Dt3}, every $B^{(i)}$ is a noetherian AS regular domain of dimension $i$ generated in degree one, with Gorenstein parameter $i$ by Corollary \ref{YAMc1Rg5}, and with basis $x_1,\dots,x_i$ of $(B^{(i)})_1$. The monomials $x_1^{a_1}\cdots x_i^{a_i}$ form a basis of $B^{(i)}$, and since $\sigma_i(x_j)$ and $\delta_i(x_j)$ are expressions in $x_1,\dots,x_{i-1}$, the algebra $B$ is presented by the generators $x_1,\dots,x_n$ and the relations
\begin{equation}\label{YAMe13Pr4}
x_ix_j=\sigma_i(x_j)x_i+\delta_i(x_j)\quad(1\le j<i\le n).
\end{equation}

\begin{example}\label{YAMx2Oe1}
The quantum plane is $\kk[x][y;\sigma]$ with $\sigma(x)=qx$. The \emph{Jordan plane} $J=\kk\langle x,y\rangle/(yx-xy-x^2)$ is $\kk[x][y;\delta]$ with $\delta(x)=x^2$. Since $\delta(x^n)=nx^{n+1}$, the relation \eqref{YAMe9Or9} for $a=x^n$ reads
\[
yx^n=x^ny+nx^{n+1}\quad(n\ge0).
\]
For $n=1$ this gives $yx=x(x+y)$ and $xy=(y-x)x$, so $xJ=Jx$ and $x$ is normal. On the other hand $yx=xy+x^2$ does not lie in $Jy\cap J_2=\kk xy+\kk y^2$, so $y$ is not normal. Both planes have Hilbert series $(1-s)^{-2}$, dimension two and Gorenstein parameter two.
\end{example}

Let $G$ be a group. A \emph{$G$-grading} of a $\kk$-algebra $R$ is a decomposition $R=\bigoplus_{g\in G}R_g$ into subspaces with $R_gR_h\subseteq R_{gh}$, and $\supp(R)=\{g\in G\mid R_g\neq0\}$ is its support. If $R$ is $\mathbb N$-graded, the $G$-grading \emph{refines} the $\mathbb N$-grading if $R_n=\bigoplus_g(R_n\cap R_g)$ for all $n$, and it is \emph{faithful} if $\supp(R)$ generates $G$. If $R$ is a domain, then $\supp(R)$ is closed under products, and a submonoid of a finite group is a subgroup. So a faithful grading of a domain by a finite group has support equal to the group. If the grading refines the $\mathbb N$-grading, then $R_0=\kk\subseteq R_e$.

The identity component of a $G$-graded right noetherian ring is right noetherian. Indeed, for a right ideal $J$ of $R_e$, the right ideal $JR=\bigoplus_gJR_g$ of $R$ satisfies $(JR)_e=J$, so $J\mapsto JR$ preserves strict inclusions and ascending chains of right ideals of $R_e$ stabilize.

Let $R$ be $G$-graded and $\gamma\in G$. A pair $(\sigma,\delta)$ formed by an endomorphism $\sigma$ of $R$ and a $\sigma$-derivation $\delta$ is \emph{$\gamma$-compatible} if
\begin{equation}\label{YAMg6Tw1}
\sigma(R_g)\subseteq R_{\gamma g\gamma^{-1}}\quad\text{and}\quad\delta(R_g)\subseteq R_{\gamma g}\quad\text{for all }g\in G.
\end{equation}
By \cite[Theorem A]{AM2}, the grading of $R$ extends to a grading of $R[t;\sigma,\delta]$ in which $t$ is homogeneous of degree $\gamma$ if and only if $(\sigma,\delta)$ is $\gamma$-compatible, and then the extension is unique and given by
\[
R[t;\sigma,\delta]_\beta=\bigoplus_{n\ge0}R_{\beta\gamma^{-n}}\,t^n\quad\text{for all }\beta\in G.
\]

Recall that a grading of $A[t;\sigma,\delta]$ is compatible with the Ore structure if $A$ is a graded subalgebra and $t$ is homogeneous. For an iterated Ore extension \eqref{YAMq8Hi1}, every grading in which $x_1,\dots,x_n$ are homogeneous makes every $B^{(i)}$ a graded subalgebra, since $B^{(i)}$ is generated by $x_1,\dots,x_i$. Such a grading is therefore compatible with each Ore extension $B^{(i)}=B^{(i-1)}[x_i;\sigma_i,\delta_i]$, and it refines the $\mathbb N$-grading since the $x_i$ have degree one.

Let $G$ be finite. A $G$-grading of $R$ is the same as a right $\kk G$-comodule algebra structure $\rho\colon R\to R\otimes\kk G$, given by $\rho(a)=\sum_ga_g\otimes g$ for $a=\sum_ga_g$ with $a_g\in R_g$, and the coinvariants of $\rho$ form $R_e$. Dually, it is the same as a left module algebra structure over $\kk^G=\Hom_\kk(\kk G,\kk)$, in which the idempotent $p_g$ dual to $g$ acts as the projection onto $R_g$. If $G$ is abelian and $\kk$ is algebraically closed of characteristic zero, then $\kk^G\cong\kk G^\vee$ for the dual group $G^\vee=\Hom(G,\kk^\times)$, and a $G$-grading is an action of $G^\vee$ by automorphisms, in which $\chi\in G^\vee$ acts on $R_g$ by $\chi(g)$. Then $R_e=R^{G^\vee}$. This explains the terminology of the introduction.

If $B=\kk[V]$ and $G$ is abelian, then $G$ is a dual reflection group for $B$ if and only if the image of $G^\vee$ in $GL(V)$ is generated by pseudo-reflections \cite{ST,Che}. We recall the structure theory of Kirkman, Kuzmanovich and Zhang. Let $G$ be a dual reflection group for $B$. For every $g\in G$ there is an $\mathbb N$-homogeneous element $f_g\in B_g$, with $f_e=1$, such that
\begin{equation}\label{YAMq9Fb2}
B_g=f_gB_e=B_ef_g,
\end{equation}
and the lowest nonzero homogeneous component of $B_g$ is the line $\kk f_g$ \cite[Theorems 0.3(1) and 3.5(1)]{KKZ}. The mass element $m_G\in G$ is the unique element with
\begin{equation}\label{YAMt3Kk6}
\deg f_{m_G}>\deg f_g\quad\text{for all }g\neq m_G
\end{equation}
\cite[Theorem 0.3(4)]{KKZ}. We put
\[
p_G(s)=\sum_{g\in G}s^{\deg f_g},\quad N_G=\deg f_{m_G}.
\]
Recall that the covariant ring of $G$ is
\[
C(B)=B/B(B_e)_{\ge1},\quad(B_e)_{\ge1}=\bigoplus_{n\ge1}(B_e)_n,
\]
and we write $\bar b$ for the image of $b\in B$ in $C(B)$. Kirkman, Kuzmanovich and Zhang proved that $B(B_e)_{\ge1}$ is a two-sided ideal, that the images of the elements $f_g$ form a basis of $C(B)$, that $C(B)$ is a Frobenius algebra with Hilbert series $p_G(s)$, which is palindromic, and that $f_{m_G}$ is a normal element of $B$ \cite[Theorems 0.4, 0.5(1) and 3.5]{KKZ}. That is,
\begin{equation}\label{YAMt4Kv0}
\begin{gathered}
B(B_e)_{\ge1}=(B_e)_{\ge1}B,\quad C(B)=\bigoplus_{g\in G}\kk\bar f_g,\\
H_{C(B)}(s)=p_G(s)=s^{N_G}p_G(s^{-1}),\quad f_{m_G}B=Bf_{m_G}.
\end{gathered}
\end{equation}
In particular $\dim_\kk C(B)=|G|$, and $C(B)$ is $G$-graded with one-dimensional components $C(B)_g=\kk\bar f_g$. Kirkman and Zhang call $f_{m_G}$ the \emph{Jacobian} of the coaction \cite[Definition 2.1]{KZ}, and they proved that for every $g\in G$ there is $\mu_g\in\kk^\times$ with
\begin{equation}\label{YAMe27Fm1}
f_gf_{g^{-1}m_G}=\mu_g\,f_{m_G}
\end{equation}
\cite[Lemma 2.11(2)]{KZ}, applied to $g^{-1}m_G$. Comparing degrees, $\deg f_g+\deg f_{g^{-1}m_G}=N_G$, which is the palindromicity of $p_G$.

We record the elementary consequences of \eqref{YAMq9Fb2} that we use. Let $g\in G$. If $g\neq e$, then $B_g\cap B_0=0$ since $B_0\subseteq B_e$, so $\deg f_g\ge1$. Since $B$ is a domain, left multiplication by $f_g$ is injective, and
\begin{equation}\label{YAMe3Lf4}
B_n\cap B_g=f_g\,(B_e)_{n-\deg f_g}\quad(n\ge0).
\end{equation}
Summing the Hilbert series of the components gives
\begin{equation}\label{YAMh3Nc7}
H_B(s)=p_G(s)\,H_{B_e}(s).
\end{equation}
If $g\neq e$, then \eqref{YAMe3Lf4} shows that $B_1\cap B_g$ is $\kk f_g$ if $\deg f_g=1$ and zero otherwise. If moreover $(B_e)_1=0$ and $B_1\cap B_g\neq0$, then $B_2\cap B_g=f_g(B_e)_1=0$. We record this as
\begin{equation}\label{YAMe3Lf5}
\dim(B_1\cap B_g)\le1\quad(g\neq e),\quad (B_e)_1=0\ \text{and}\ B_1\cap B_g\neq0\ \Longrightarrow\ B_2\cap B_g=0.
\end{equation}
The first statement is also \cite[Theorems 0.3(2) and 3.5(3)]{KKZ}. If $\theta$ is a graded automorphism of $B$ with $\theta(B_g)=B_g$, then $\theta(f_g)=\lambda_gf_g$ for some $\lambda_g\in\kk^\times$, since $\kk f_g$ is the lowest component of $B_g$.

\begin{example}\label{YAMex1Cm4}
Let $B=\kk[x_1,\dots,x_n]$ and $G=\prod_{i=1}^n\langle c_i\rangle$ with $|c_i|=r_i$, and let $x_i$ have $G$-degree $c_i$. Then
\[
B_e=\kk[x_1^{r_1},\dots,x_n^{r_n}],\quad f_g=x_1^{a_1}\cdots x_n^{a_n}\quad\text{for }g=c_1^{a_1}\cdots c_n^{a_n},\ 0\le a_i<r_i,
\]
so $G$ is a dual reflection group, $m_G=c_1^{r_1-1}\cdots c_n^{r_n-1}$ and $f_{m_G}=\prod_ix_i^{r_i-1}$. Since $(1-s)[r]_s=1-s^r$, the identity \eqref{YAMh3Nc7} reads
\[
\frac{1}{(1-s)^n}=\prod_{i=1}^n[r_i]_s\cdot\prod_{i=1}^n\frac{1}{1-s^{r_i}},\quad p_G(s)=\prod_{i=1}^n[r_i]_s.
\]
The Jacobian of the basic invariants is $\det\big(\partial x_j^{r_j}/\partial x_i\big)=\prod_ir_ix_i^{r_i-1}$, a nonzero multiple of $f_{m_G}$. Here $G$ is the character group of the reflection group $\prod_i\mu_{r_i}$ acting diagonally on $\kk^n$, and $C(B)=\bigotimes_i\kk[x_i]/(x_i^{r_i})$.
\end{example}

\begin{example}\label{YAMex8Qv2}
Let $A=\kk_q[x,y]$ be the quantum plane, graded by $\mathbb Z_2=\{e,\tau\}$ with $x\in A_\tau$ and $y\in A_e$. Then $A_e$ is spanned by the monomials $x^{2a}y^b$, and $yx^2=q^2x^2y$, so
\[
A_e=\kk\langle x^2,y\rangle\cong\kk_{q^2}[x^2,y],\quad f_\tau=x,\quad m_G=\tau,\quad C(A)=\kk[\bar x]/(\bar x^2).
\]
So $\mathbb Z_2$ is a dual reflection group for $A$, and
\[
H_A(s)=\frac{1}{(1-s)^2}=(1+s)\cdot\frac{1}{(1-s^2)(1-s)}=p_G(s)\,H_{A_e}(s).
\]
The ideal $A(A_e)_{\ge1}$ is generated by $x^2$ and $y$. The nontrivial character of $\mathbb Z_2$ acts by $x\mapsto-x$ and $y\mapsto y$, and for $q=1$ this is the reflection of $\kk^2$ in the line $x=0$.
\end{example}

\section{Descent of the Gorenstein property}\label{YAMs1Pz4}

In this section we prove Theorem \ref{YAMa7Qd2}. Throughout, $R$ is a noetherian connected graded algebra, $\sigma$ is a graded automorphism of $R$, $\delta$ is a $\sigma$-derivation satisfying \eqref{YAMq5Ef7} for a fixed $m\ge1$, and $S=R[t;\sigma,\delta]$ is graded by $\deg t=m$. By Section \ref{YAMs1Pr0}, $S$ is a right noetherian connected graded algebra, and it is left noetherian by the left version of the same argument, since $\sigma$ is an automorphism \cite[Theorem 1.2.9]{MR}. The proof compares $\Ext_R(\kk,R)$ with $\Ext_S$ of the cyclic module $S/R_+S$.

\begin{lemma}\label{YAMl5Gd0}
The following hold.
\begin{enumerate}
\item The algebra $S$ is free as a left and as a right $R$-module on the powers of $t$,
\begin{equation}\label{YAMe15Fr6}
S=\bigoplus_{n\ge0}Rt^n=\bigoplus_{n\ge0}t^nR,\quad S_R\cong\bigoplus_{n\ge0}R(-mn)\cong{}_RS.
\end{equation}
\item The ideal $I=R_+S$ is two-sided and equal to $SR_+$, and $Q=S/I$ is the polynomial ring $\kk[\bar t]$ with $\deg\bar t=m$. As a right and as a left $R$-module,
\begin{equation}\label{YAMe16Qr7}
Q\cong\bigoplus_{n\ge0}\kk(-mn).
\end{equation}
\item There are exact sequences of graded right and left $S$-modules, respectively,
\begin{equation}\label{YAMe17Ex8}
0\to Q(-m)\xrightarrow{\ \bar t\,\cdot\ }Q\to\kk_S\to0\quad\text{and}\quad0\to Q(-m)\xrightarrow{\ \cdot\,\bar t\ }Q\to{}_S\kk\to0.
\end{equation}
\end{enumerate}
\end{lemma}

\begin{proof}
(1) The first decomposition is the definition of $S$. By \eqref{YAMe10Tn1}, $t^na-\sigma^n(a)t^n\in\sum_{j<n}Rt^j$ for $a\in R$. Since $\sigma^n$ is surjective, induction on $n$ gives $\sum_{j\le n}t^jR=\sum_{j\le n}Rt^j$. If $\sum_{j\le n}t^ja_j=0$ with $a_n\neq0$, then the coefficient of $t^n$ in the basis $\{t^j\}$ is $\sigma^n(a_n)\neq0$, since $\sigma^n$ is injective, a contradiction. So the sum $\sum_nt^nR$ is direct, and $a\mapsto t^na$ is an isomorphism of graded right $R$-modules $R(-mn)\to t^nR$.

(2) Since $\sigma$ is a graded bijection and $\delta$ raises degrees, $\sigma(R_+)=R_+$ and $\delta(R_+)\subseteq R_+$. For $a\in R_+$ put $a'=\sigma^{-1}(a)\in R_+$. Then
\[
ta=\sigma(a)t+\delta(a)\in R_+S,\quad at=ta'-\delta(a')\in SR_+.
\]
So $tR_+\subseteq R_+S$ and $R_+t\subseteq SR_+$. Since $S$ is generated by $R$ and $t$, the right ideal $R_+S$ and the left ideal $SR_+$ are two-sided, and both equal the ideal generated by $R_+$. By (1), $R_+S=\bigoplus_nR_+t^n$, so $Q$ has basis $\{\bar t^n\}$, with $\deg\bar t^n=mn$, and $Q=\kk[\bar t]$. Since $QR_+=0=R_+Q$, the algebra $R$ acts on $Q$ through $R/R_+=\kk$ on both sides, which gives \eqref{YAMe16Qr7}.

(3) Left multiplication by $\bar t$ is injective on $\kk[\bar t]$, it raises degrees by $m$, and it is right $S$-linear because $I$ is two-sided. Its cokernel is $S/(I+tS)=S/S_+=\kk$, since $S_+=R_+S+tS$. The left sequence is obtained in the same way from $S_+=SR_++St$.
\end{proof}

\begin{lemma}\label{YAMl8Mb3}
Let $M$ be a finitely generated graded right $R$-module. For every $i\ge0$ there is an isomorphism of graded vector spaces
\[
\Ext^i_S(M\otimes_RS,S)\cong\bigoplus_{n\ge0}\Ext^i_R(M,R)(-mn).
\]
In particular, if $X^i=\Ext^i_R(\kk_R,R)$, then $\Ext^i_S(Q,S)\cong\bigoplus_{n\ge0}X^i(-mn)$. The same holds for left modules, with $S\otimes_RM$ in place of $M\otimes_RS$.
\end{lemma}

\begin{proof}
Let $P_\bullet\to M$ be a resolution by finitely generated graded free right $R$-modules. Since ${}_RS$ is free, $P_\bullet\otimes_RS$ is a free resolution of $M\otimes_RS$. For a free right $R$-module $P$ we have $\Hom_S(P\otimes_RS,S)\cong\Hom_R(P,S_R)$, hence $\Ext^i_S(M\otimes_RS,S)\cong\Ext^i_R(M,S_R)$. By \eqref{YAMe15Fr6}, $S_R\cong\bigoplus_nR(-mn)$, and $\Ext^i_R(M,-)$ commutes with this direct sum because the $P_i$ are finitely generated. Altogether,
\[
\Ext^i_S(M\otimes_RS,S)\cong\Ext^i_R(M,S_R)\cong\Ext^i_R\Big(M,\bigoplus_{n\ge0}R(-mn)\Big)\cong\bigoplus_{n\ge0}\Ext^i_R(M,R)(-mn).
\]
The second statement is the case $M=\kk$, since $\kk\otimes_RS=S/R_+S=Q$. For left modules, the same argument applies to $S\otimes_R-$, using that $S_R$ is free, so that $S\otimes_RP_\bullet$ is a free resolution of $S\otimes_RM$, and that ${}_RS\cong\bigoplus_nR(-mn)$. Both uses rely on the two free decompositions \eqref{YAMe15Fr6}, hence on the bijectivity of $\sigma$.
\end{proof}

We bound the injective dimension of $R$ with local duality for $S$. Let $T$ be a noetherian connected graded algebra and $N$ a graded right $T$-module. The \emph{torsion submodule} of $N$ is
\[
\Gamma_{\mathfrak m}(N)=\{b\in N\mid bT_{\ge k}=0\text{ for some }k\ge0\},
\]
and the \emph{local cohomology} of $N$ is $H^i_{\mathfrak m}(N)=\mathrm R^i\Gamma_{\mathfrak m}(N)$, where $\mathfrak m=T_+$, computed in the category of graded modules. The \emph{graded dual} of a graded vector space $N$ is $N'$ with $(N')_n=\Hom_\kk(N_{-n},\kk)$. If $N$ is finitely generated, then $\Gamma_{\mathfrak m}(N)$ is a finitely generated submodule annihilated by some $T_{\ge k}$, hence finite-dimensional. If a noetherian connected graded algebra $T$ is AS Gorenstein of dimension $D$ with Gorenstein parameter $\ell$, then it has a balanced dualizing complex of the form ${}^{\nu}T(-\ell)[D]$ for a graded automorphism $\nu$ of $T$ \cite[Section 4]{YZ}, see also \cite[Theorem 4.12]{LW} and \cite{Ye}. Then $T$ satisfies the $\chi$-condition and has finite local cohomological dimension \cite[Corollary 4.3]{YZ}, see also \cite[Corollary 4.13]{LW} and \cite[Theorem 6.3]{VdB}, and local duality \cite[Theorem 5.1]{VdB} gives isomorphisms of graded vector spaces
\begin{equation}\label{YAMe63Ld1}
\Ext^i_T(N,T)\cong H^{D-i}_{\mathfrak m}(N)'(\ell),
\end{equation}
for every finitely generated graded right or left $T$-module $N$. For $i=D$ this gives $\Ext^D_T(N,T)\cong\Gamma_{\mathfrak m}(N)'(\ell)$, which is finite-dimensional. Only this finite-dimensionality is used below.

\begin{proof}[\textnormal{\bfseries Proof of Theorem \ref{YAMa7Qd2}}]
If $b\in S$ satisfies $bS_+=0$, then $bt=0$, hence $b=0$ by \eqref{YAMe15Fr6}. Thus $\Hom_S(\kk,S)=0$, and $D\ge1$.

We compute the spaces $X^i=\Ext^i_R(\kk,R)$ of Lemma \ref{YAMl8Mb3}. Since $R$ is noetherian, each $X^i$ is bounded below and finite-dimensional in each degree. Put $Y^i=\Ext^i_S(Q,S)$. Then $\Ext^i_S(Q(-m),S)=Y^i(m)$, and the first sequence of \eqref{YAMe17Ex8} gives a long exact sequence of graded vector spaces
\begin{equation}\label{YAMe19Le2}
\cdots\to\Ext^i_S(\kk,S)\to Y^i\xrightarrow{\ \varphi_i\ }Y^i(m)\to\Ext^{i+1}_S(\kk,S)\to Y^{i+1}\to\cdots,
\end{equation}
where each $\varphi_i$ preserves degrees. Since $S$ is AS Gorenstein, $\Ext^j_S(\kk,S)=0$ for $j\neq D$ and $\Ext^D_S(\kk,S)=\kk(\ell)$.

Let $i\neq D-1$. Then $\Ext^{i+1}_S(\kk,S)=0$, so $\varphi_i$ is surjective. Suppose that $X^i\neq0$, and let $a$ be its lowest degree. By Lemma \ref{YAMl8Mb3}, $Y^i\cong\bigoplus_{n\ge0}X^i(-mn)$ has lowest degree $a$, so it vanishes in degree $a-m$, while $Y^i(m)$ does not. This contradicts the surjectivity of $\varphi_i$, so $X^i=0$. For $i=D-1$ we have $Y^D=0$ by the previous case, and \eqref{YAMe19Le2} reduces to the exact sequence
\begin{equation}\label{YAMe20Sh3}
0\to Y^{D-1}\xrightarrow{\ \varphi_{D-1}\ }Y^{D-1}(m)\to\kk(\ell)\to0.
\end{equation}
Let $h(s)$ be the Hilbert series of $X^{D-1}$. By Lemma \ref{YAMl8Mb3}, the Hilbert series of $Y^{D-1}$ and of $Y^{D-1}(m)$ are $h(s)/(1-s^m)$ and $s^{-m}h(s)/(1-s^m)$. Exactness of \eqref{YAMe20Sh3} gives
\begin{equation}\label{YAMe21Hs4}
\begin{aligned}
\frac{s^{-m}h(s)}{1-s^m}-\frac{h(s)}{1-s^m}&=s^{-\ell},\\
s^{-m}h(s)&=s^{-\ell},\quad h(s)=s^{m-\ell}.
\end{aligned}
\end{equation}
So $X^{D-1}$ is one-dimensional, concentrated in degree $m-\ell$. A one-dimensional graded module is annihilated by $R_+$, hence $X^{D-1}\cong\kk(\ell-m)$. The second sequence of \eqref{YAMe17Ex8} and the left version of Lemma \ref{YAMl8Mb3} give the same statement for left modules. Altogether,
\begin{equation}\label{YAMe21Hs5}
\Ext^i_R(\kk,R)\cong\begin{cases}\kk(\ell-m)&i=D-1,\\0&i\neq D-1\end{cases},
\end{equation}
on both sides.

It remains to show that $\injdim R_R=\injdim{}_RR=D-1$. Let $M$ be a finitely generated graded right $R$-module, and put $N=M\otimes_RS$, a finitely generated graded right $S$-module. By \eqref{YAMe63Ld1} with $T=S$, the space $\Ext^D_S(N,S)$ is finite-dimensional. By Lemma \ref{YAMl8Mb3},
\[
\Ext^D_S(N,S)\cong\bigoplus_{n\ge0}\Ext^D_R(M,R)(-mn).
\]
If $\Ext^D_R(M,R)$ were nonzero in some degree $a$, this space would be nonzero in every degree $a+mn$ with $n\ge0$. Hence $\Ext^D_R(M,R)=0$ for every such $M$, so $\injdim R_R\le D-1$ by \eqref{YAMe1Ba1}, and equality holds since $X^{D-1}\neq0$. The left statement is the same, using the left version of Lemma \ref{YAMl8Mb3} and local duality for left modules. This proves that $R$ is AS Gorenstein of dimension $D-1$ with Gorenstein parameter $\ell-m$.

Suppose now that $S$ is AS regular. It remains to show that $\gldim R=D-1$ and $\GKdim R<\infty$. Since $S_R$ is free, every projective right $S$-module is projective over $R$. A projective resolution of $Q$ over $S$ of length at most $D$ therefore restricts to a projective resolution of $Q_R$ over $R$. By \eqref{YAMe16Qr7}, $\kk_R$ is a direct summand of $Q_R$, so $p=\pd\kk_R\le D$. By \eqref{YAMe1Ba2}, $\Ext^p_R(\kk,R)\neq0$, so $p=D-1$ by \eqref{YAMe21Hs5}, and $\gldim R=\pd\kk_R=D-1$. Since $R$ is a subalgebra of $S$, we have $\GKdim R\le\GKdim S<\infty$.
\end{proof}

\begin{corollary}\label{YAMc1Rg5}
Let $S=R[t;\sigma,\delta]$ be a graded Ore extension as above. Then $S$ is AS regular if and only if $R$ is AS regular. In this case
\begin{equation}\label{YAMe22Dm5}
\dim S=\dim R+1\quad\text{and}\quad\ell_S=\ell_R+m.
\end{equation}
\end{corollary}

\begin{proof}
The implication from $R$ to $S$ is \eqref{YAMt1Lw9}, and the converse together with \eqref{YAMe22Dm5} is Theorem \ref{YAMa7Qd2}.
\end{proof}

Iterating Theorem \ref{YAMa7Qd2} gives the descent along iterated Ore extensions, with arbitrary degrees of the variables.

\begin{corollary}\label{YAMc9It2}
Let $R$ be a noetherian connected graded algebra, put $S^{(0)}=R$, and for $1\le j\le k$ let
\[
S^{(j)}=S^{(j-1)}[t_j;\sigma_j,\delta_j],\quad\deg t_j=m_j,
\]
where $\sigma_j$ is a graded automorphism of $S^{(j-1)}$ and $\delta_j$ is a $\sigma_j$-derivation of $S^{(j-1)}$ raising degrees by $m_j\ge1$. If $S^{(k)}$ is AS Gorenstein of dimension $D$ with Gorenstein parameter $\ell$, then $R$ is AS Gorenstein of dimension $D-k$ with Gorenstein parameter
\[
\ell_R=\ell-(m_1+m_2+\cdots+m_k).
\]
If moreover $S^{(k)}$ is AS regular, then $R$ is AS regular.
\end{corollary}

\begin{proof}
Each $S^{(j)}$ is noetherian and connected graded. Applying Theorem \ref{YAMa7Qd2} to $S^{(j)}=S^{(j-1)}[t_j;\sigma_j,\delta_j]$ for $j=k,k-1,\dots,1$ shows that $S^{(j-1)}$ is AS Gorenstein of dimension $D-(k-j+1)$ with Gorenstein parameter $\ell-(m_j+\cdots+m_k)$, and AS regular if $S^{(k)}$ is.
\end{proof}

\begin{remark}\label{YAMrm3Xa6}
When $\delta=0$, the element $t$ is normal and regular of degree $m$ in $S$, and $R\cong S/tS$. The Rees lemma \cite[Proposition 3.4(b)]{Lev} then gives $\Ext^{i+1}_S(\kk,S)\cong\Ext^i_R(\kk,R)(m)$, which is \eqref{YAMe21Hs5} directly. The proof above uses no normal element. Its point is the case $\delta\neq0$, where $t$ is in general not normal. The hypotheses enter as follows. The decomposition $S=\bigoplus_nRt^n$ holds by definition, while $S=\bigoplus_nt^nR$ uses the bijectivity of $\sigma$. Together they give the right and left versions of Lemma \ref{YAMl8Mb3}. The noetherian hypothesis gives finitely generated free resolutions, and it makes $S$ noetherian, which is needed for local duality. In the regular case local duality can be avoided, since $\gldim R=D-1$ follows from \eqref{YAMe21Hs5} and \eqref{YAMe1Ba2} alone and bounds the injective dimension. The condition $m\ge1$ makes $S$ connected graded, and the shift by $m$ in \eqref{YAMe17Ex8} produces the Gorenstein parameter $\ell-m$.
\end{remark}

\begin{example}\label{YAMex10Jd2}
Let $S=\kk[x][y;\delta]$ be the Jordan plane of Example \ref{YAMx2Oe1}. Here $R=\kk[x]$, $m=1$, $D=2$ and $\ell=2$, and $Q=S/xS=\kk[\bar y]$. Since $X^1=\Ext^1_{\kk[x]}(\kk,\kk[x])=\kk(1)$ and $X^i=0$ for $i\neq1$, Lemma \ref{YAMl8Mb3} gives $Y^1=\Ext^1_S(Q,S)\cong\bigoplus_{n\ge0}\kk(1-n)$, with Hilbert series $s^{-1}/(1-s)$. The exact sequence \eqref{YAMe20Sh3} reads
\[
0\to Y^1\xrightarrow{\ \varphi_1\ }Y^1(1)\to\kk(2)\to0,\quad\frac{s^{-2}}{1-s}-\frac{s^{-1}}{1-s}=s^{-2},
\]
which is \eqref{YAMe21Hs4} with $h(s)=s^{-1}=s^{m-\ell}$. Here $R$ is not the quotient of $S$ by the Ore variable $y$, which is not normal, so this computation does not come from the Rees lemma.
\end{example}

\begin{example}\label{YAMex11Dm3}
The degree of $t$ may exceed one. Let $R=\kk[x]$ and let $\delta$ be the derivation with $\delta(x)=x^3$, so that \eqref{YAMq5Ef7} holds with $m=2$. The algebra
\[
S=\kk[x][t;\delta]=\kk\langle x,t\rangle/(tx-xt-x^3),\quad\deg x=1,\quad\deg t=2,
\]
is AS regular of dimension two by \eqref{YAMt1Lw9}, with Hilbert series $H_S(s)=(1-s)^{-1}(1-s^2)^{-1}$, and \eqref{YAMe22Dm5} gives $\ell_S=\ell_{\kk[x]}+2=3$. Here $t$ is not normal, since $tx=xt+x^3$ does not lie in $St\cap S_3=\kk xt$.
\end{example}

\begin{example}\label{YAMex17Gn1}
Theorem \ref{YAMa7Qd2} also applies to algebras that are not AS regular. Let $R=\kk[x]/(x^2)$ with $\deg x=1$, and $S=R[t]=\kk[x,t]/(x^2)$ with $\deg t=1$. The algebra $S$ is a hypersurface in $\kk[x,t]$, hence a commutative graded Gorenstein ring of dimension one, and its $a$-invariant is $\deg x^2-2=0$ \cite[\S3.6]{BH}. In our convention the Gorenstein parameter is $\ell=-a$. So $S$ is AS Gorenstein of dimension one with Gorenstein parameter $0$. By Theorem \ref{YAMa7Qd2}, $R$ is AS Gorenstein of dimension $0$ with Gorenstein parameter $-1$. This agrees with the direct computation. The algebra $R$ is self-injective, and $\Hom_R(\kk,R)$ is spanned by the inclusion of the socle $\kk x$, which lies in degree one, so $\Hom_R(\kk,R)\cong\kk(-1)$.
\end{example}

\section{Dual reflection groups of an Ore extension}\label{YAMs4Or0}

In this section we prove Theorem \ref{YAMb4Lx9}, describe the elements $f_g$ of $B$, and discuss the compatibility assumption. From now on $\kk$ is algebraically closed of characteristic zero.

Throughout this section, $A$ is a noetherian AS regular domain generated in degree one, $\sigma$ is a graded automorphism of $A$, $\delta$ is a $\sigma$-derivation with $\delta(A_n)\subseteq A_{n+1}$, and $B=A[t;\sigma,\delta]$ with $\deg t=1$. By \eqref{YAMt1Lw9} and \eqref{YAMe29Dt3}, $B$ is a noetherian AS regular domain generated in degree one, and the $t$-degree is additive on $B$. We fix a finite group $G$ and a $G$-grading of $B$ compatible with the Ore structure, with $t\in B_\gamma$, and we put
\[
K=\supp(A),\quad r=|\gamma|.
\]
Since $A$ is a domain, $K$ is a subgroup of $G$. The grading of $A$ refines its $\mathbb N$-grading, since $A$ is a graded subalgebra of $B$. By Section \ref{YAMs1Pr0}, the pair $(\sigma,\delta)$ is $\gamma$-compatible and
\begin{equation}\label{YAMe31Bc5}
B_\beta=\bigoplus_{n\ge0}A_{\beta\gamma^{-n}}\,t^n\quad\text{for all }\beta\in G.
\end{equation}

\begin{lemma}\label{YAMl4Ae9}
The following hold for every grading of $B$ compatible with the Ore structure.
\begin{enumerate}
\item $\sigma(A_g)=A_{\gamma g\gamma^{-1}}$ for every $g\in G$. In particular $\sigma(A_e)=A_e$, and $\gamma K\gamma^{-1}=K$.
\item $\supp(B)=K\langle\gamma\rangle$. Hence $G=K\langle\gamma\rangle$, the subgroup $K$ is normal in $G$, and $|G|=|K|\,r/|K\cap\langle\gamma\rangle|$.
\item If $\delta\neq0$, then $\gamma\in K$.
\end{enumerate}
\end{lemma}

\begin{proof}
(1) The equality $\sigma(A_g)=A_{\gamma g\gamma^{-1}}$ is \cite[Lemma 2.3(3)]{AM2}. Since $\sigma$ is injective, $A_g\neq0$ if and only if $A_{\gamma g\gamma^{-1}}\neq0$.

(2) By \eqref{YAMe31Bc5}, $B_\beta\neq0$ if and only if $\beta\gamma^{-n}\in K$ for some $n\ge0$, that is, $\beta\in K\gamma^n$. Since $\gamma$ has finite order, $\supp(B)=\bigcup_nK\gamma^n=K\langle\gamma\rangle$. By (1), $\gamma$ normalizes $K$, so $K\langle\gamma\rangle$ is a subgroup in which $K$ is normal, and faithfulness gives $G=K\langle\gamma\rangle$. The order formula is the product formula for the subgroups $K$ and $\langle\gamma\rangle$.

(3) Since $\delta\neq0$ and $A=\bigoplus_gA_g$, there is a homogeneous $a\in A_g$ with $\delta(a)\neq0$. By \eqref{YAMg6Tw1}, $0\neq\delta(a)\in A_{\gamma g}$. So $g\in K$ and $\gamma g\in K$, and $\gamma=(\gamma g)g^{-1}\in K$.
\end{proof}

\begin{lemma}\label{YAMi6Ub2}
If $G$ is a dual reflection group for $B$, then $K\cap\langle\gamma\rangle=\{e\}$.
\end{lemma}

\begin{proof}
Suppose that $h=\gamma^a\in K$ with $h\neq e$ and $a\ge1$. Let $f_h$ be as in \eqref{YAMq9Fb2} for $B$. Choose $0\neq b\in A_h$ and write $b=f_hc$ with $c\in B_e$. By \eqref{YAMe29Dt3}, $\degt f_h+\degt c=\degt b=0$, so $f_h\in A$, and $\deg f_h\ge1$ since $h\neq e$. On the other hand $t^a\in B_h=f_hB_e$, so
\[
t^a=f_hc'=\sum_if_hc'_it^i\quad\text{with }c'=\sum_ic'_it^i\in B_e,\ c'_i\in A.
\]
Comparing the coefficients of $t^a$ in the free left $A$-module $\bigoplus_iAt^i$ gives $f_hc'_a=1$ in $A$. This is impossible, since $A$ is connected graded and $f_h\in A_{\ge1}$.
\end{proof}

The same argument appears in the proof of \cite[Proposition 6.1]{BR} when $\delta=0$. Here we only use that the $t$-degree is additive, which holds for every $\delta$. Suppose now that $K\cap\langle\gamma\rangle=\{e\}$. For $k\in K$ and $0\le j<r$, the element $k\gamma^{j-n}$ lies in $K$ if and only if $r$ divides $n-j$. So \eqref{YAMe31Bc5} gives
\begin{equation}\label{YAMe76Cp4}
B_{k\gamma^j}=\bigoplus_{n\ge0}A_{k\gamma^{j-n}}\,t^n=\bigoplus_{l\ge0}A_k\,t^{j+rl}\quad(k\in K,\ 0\le j<r).
\end{equation}
In particular $B_e=\bigoplus_{l\ge0}A_et^{rl}$. If $\delta\neq0$, then $\gamma\in K$ by Lemma \ref{YAMl4Ae9}(3), so $\gamma=e$ and $r=1$. In this case $\sigma(A_e)=A_e$ and $\delta(A_e)\subseteq A_{\gamma}=A_e$. If $\delta=0$, then $t^ra=\sigma^r(a)t^r$ for $a\in A$. In both cases
\begin{equation}\label{YAMe34Be8}
B_e=A_e[t^r;\sigma^r|_{A_e},\delta|_{A_e}],\quad\deg t^r=r,
\end{equation}
is a graded Ore extension of $A_e$, where $\delta=0$ whenever $r>1$. The algebra $A_e$ is noetherian by Section \ref{YAMs1Pr0}, since $A$ is.

\begin{proof}[\textnormal{\bfseries Proof of Theorem \ref{YAMb4Lx9}}]
Assume that $G$ is a dual reflection group for $B$. By Lemma \ref{YAMi6Ub2}, $K\cap\langle\gamma\rangle=\{e\}$, so \eqref{YAMe34Be8} holds. The algebra $B_e$ is AS regular, $A_e$ is noetherian, $\sigma^r|_{A_e}$ is a graded automorphism of $A_e$, and $\delta|_{A_e}$ raises degrees by $r$ since it vanishes when $r>1$. So Theorem \ref{YAMa7Qd2} applies with $m=r$ and shows that $A_e$ is AS regular. The grading of $A$ by $K$ is faithful and refines the $\mathbb N$-grading. Hence $K$ is a dual reflection group for $A$.

Conversely, assume that $K$ is a dual reflection group for $A$ and that $K\cap\langle\gamma\rangle=\{e\}$. Then \eqref{YAMe34Be8} exhibits $B_e$ as a graded Ore extension of the AS regular algebra $A_e$, so $B_e$ is AS regular by \eqref{YAMt1Lw9}. The grading of $B$ is faithful and refines the $\mathbb N$-grading by our conventions, so $G$ is a dual reflection group for $B$.

In this case $K$ is normal in $G=K\langle\gamma\rangle$ by Lemma \ref{YAMl4Ae9}(2), and $K\cap\langle\gamma\rangle=\{e\}$, so $G=K\rtimes\langle\gamma\rangle$. If $\delta\neq0$, we saw that $\gamma=e$, and then \eqref{YAMg6Tw1} shows that $\sigma$ and $\delta$ preserve every component $A_g$. It remains to prove that $(A_e)_1\neq0$ when $\delta\neq0$. Suppose that $(A_e)_1=0$. Then $A_1=\bigoplus_{g\neq e}(A_1\cap A_g)$. Let $g\neq e$ with $A_1\cap A_g\neq0$. By \eqref{YAMe3Lf5} applied to $A$, we have $A_2\cap A_g=0$. Since $\delta$ preserves $A_g$ and raises degrees by one,
\[
\delta(A_1\cap A_g)\subseteq A_2\cap A_g=0.
\]
Hence $\delta$ vanishes on $A_1$. Since $A$ is generated by $A_1$, \eqref{YAMe7Sd7} gives $\delta=0$, a contradiction.
\end{proof}

We describe the elements $f_g$ of $B$ in terms of those of $A$.

\begin{proposition}\label{YAMn5Fh3}
In the setting of Theorem \ref{YAMb4Lx9}, suppose that $G$ is a dual reflection group for $B$, and let $f_k$, for $k\in K$, be the elements of \eqref{YAMq9Fb2} for $A$. Then for $k\in K$ and $0\le j<r$,
\begin{equation}\label{YAMe35Fg9}
B_{k\gamma^j}=f_kt^j\,B_e=B_e\,f_kt^j.
\end{equation}
So the elements $f_{k\gamma^j}=f_kt^j$ play the role of \eqref{YAMq9Fb2} for $B$, and
\begin{equation}\label{YAMe45Pg1}
p_G(s)=p_K(s)\,[r]_s,\quad m_G=m_K\gamma^{r-1},\quad f_{m_G}=f_{m_K}t^{r-1},\quad N_G=N_K+r-1.
\end{equation}
\end{proposition}

\begin{proof}
By \eqref{YAMe76Cp4} and \eqref{YAMq9Fb2} for $A$,
\[
B_{k\gamma^j}=\bigoplus_{l\ge0}f_kA_e\,t^{j+rl}.
\]
If $r=1$, then $j=0$. If $r>1$, then $\delta=0$ and $t^jA_e=\sigma^j(A_e)t^j=A_et^j$ by Lemma \ref{YAMl4Ae9}(1). In both cases $A_et^j=t^jA_e$, hence $f_kA_et^{j+rl}=f_kt^jA_et^{rl}$, and $B_{k\gamma^j}=f_kt^jB_e$. For the other side, $B_e\,f_kt^j=\bigoplus_lA_et^{rl}f_kt^j$. If $r>1$, then $\delta=0$ and $t^{rl}f_k=\sigma^{rl}(f_k)t^{rl}$, and $\sigma^{rl}(f_k)$ is a nonzero multiple of $f_k$, since $\sigma^r$ is a graded automorphism preserving $A_k$ by Lemma \ref{YAMl4Ae9}(1) and $\kk f_k$ is the lowest component of $A_k$. If $r=1$, then $j=0$ and $B_ef_k\subseteq B_k=f_kB_e$, and both sides have Hilbert series $s^{\deg f_k}H_{B_e}(s)$ because $B$ is a domain. In both cases $B_ef_kt^j=B_{k\gamma^j}$. The lowest component of $B_{k\gamma^j}$ is therefore $\kk f_kt^j$, of degree $\deg f_k+j$. Since every element of $G=K\rtimes\langle\gamma\rangle$ is uniquely of the form $k\gamma^j$ with $0\le j<r$, summing over $G$ gives $p_G(s)=p_K(s)[r]_s$, and the maximal degree is attained exactly at $k=m_K$ and $j=r-1$.
\end{proof}

When $\delta\neq0$, the elements $f_g$ of $B$ are those of $A$, and the $\sigma$-derivation acts trivially on the covariant ring of $A$, as follows.

\begin{proposition}\label{YAMet9Lw3}
In the setting of Theorem \ref{YAMb4Lx9}, suppose that $G$ is a dual reflection group for $B$ and that $\delta\neq0$. For every $g\in G$ there are $\lambda_g\in\kk^\times$ and $a_g\in(A_e)_1$ with
\begin{equation}\label{YAMe36Et1}
\sigma(f_g)=\lambda_gf_g,\quad\delta(f_g)=f_ga_g,\quad tf_g=f_g(\lambda_gt+a_g).
\end{equation}
Moreover $\delta(A)\subseteq A(A_e)_{\ge1}$, so $\delta$ induces the zero map on $C(A)$.
\end{proposition}

\begin{proof}
By Theorem \ref{YAMb4Lx9}, $\sigma$ and $\delta$ preserve $A_g$. Since $\sigma$ preserves degrees, $\sigma(f_g)=\lambda_gf_g$ with $\lambda_g\neq0$. The element $\delta(f_g)\in A_g$ has degree $\deg f_g+1$, so $\delta(f_g)=f_ga_g$ with $a_g\in(A_e)_1$ by \eqref{YAMe3Lf4}. Then $tf_g=\sigma(f_g)t+\delta(f_g)=f_g(\lambda_gt+a_g)$. For $b\in A_e$, the element $\delta(b)$ lies in $(A_e)_{\ge1}$, since $\delta$ preserves $A_e$ and raises degrees. Hence
\[
\delta(f_gb)=\sigma(f_g)\delta(b)+\delta(f_g)b=f_g\big(\lambda_g\delta(b)+a_gb\big)\in f_g(A_e)_{\ge1}.
\]
Since $A=\bigoplus_gf_gA_e$, this proves $\delta(A)\subseteq A(A_e)_{\ge1}$.
\end{proof}

The relation $B_ef_g=f_gB_e$ of \eqref{YAMe35Fg9} defines a graded automorphism $\eta_g$ of $B_e$ by $bf_g=f_g\eta_g(b)$, and \eqref{YAMe36Et1} says that $\eta_g(t)=\lambda_gt+a_g$. This is the automorphism $\eta_g$ of \cite[(E3.5.1)]{KKZ}.

The first example is one of the two Ore extensions with $\delta\neq0$ in the list of Buzaglo and Rogalski.

\begin{example}\label{YAMx7Lp1}
Let $A_{-1}=\kk\langle u,v\rangle/(u^2-v^2)$, let $\sigma$ be the automorphism exchanging $u$ and $v$, and let $\delta$ be the $\sigma$-derivation with $\delta(u)=uv$ and $\delta(v)=vu$. Then
\[
\delta(u^2-v^2)=\sigma(u)\delta(u)+\delta(u)u-\sigma(v)\delta(v)-\delta(v)v=vuv+uvu-uvu-vuv=0,
\]
so $\delta$ is well defined, and $B=A_{-1}[z;\sigma,\delta]$ is the algebra with relations
\[
u^2=v^2,\quad zu=vz+uv,\quad zv=uz+vu
\]
of \cite[Table 1]{BR}. Consider a grading of $B$ by a finite group $G$ for which $u,v,z$ are homogeneous, of degrees $a,b,c$. By \cite[Theorem 3.4]{BR}, every grading of $B$ by a nonabelian group is of this form up to a graded change of variables. The subalgebra $A_{-1}$ is graded, so the grading is compatible with the Ore structure, with $\gamma=c$. By \eqref{YAMg6Tw1} applied to $\sigma(u)=v$ and $\delta(u)=uv$, we get $b=cac^{-1}$ and $ab=ca$. If $G$ is a dual reflection group for $B$, then $c=e$ by Theorem \ref{YAMb4Lx9}, hence $b=a$ and $a^2=a$. So $a=b=c=e$ and $G$ is trivial. Hence $B$ has no nonabelian dual reflection group, in agreement with \cite[Theorem 1.2]{BR}, and the argument uses the classification of gradings \cite[Theorem 3.4]{BR} but not the classification of dual reflection groups.

On the other hand, $B$ has the $\mathbb Z_2$-grading of \cite[Remark 5.16]{BR}, defined in the basis $x=(u+v)/2$, $y=(u-v)/2$, $z$. In this basis $u^2=v^2$ reads $yx=-xy$, and adding and subtracting the two other relations gives
\begin{equation}\label{YAMe12Br6}
yx=-xy,\quad zx=xz+x^2-y^2,\quad zy=-yz-2xy.
\end{equation}
So $B$ is the iterated Ore extension $\kk[x][y;\sigma_2][z;\sigma_3,\delta_3]$ with $\sigma_2(x)=-x$, $\sigma_3(x)=x$, $\sigma_3(y)=-y$, $\delta_3(x)=x^2-y^2$ and $\delta_3(y)=-2xy$. Write $\mathbb Z_2=\{e,\tau\}$ and grade $B$ by $x,z\in B_e$ and $y\in B_\tau$. The subalgebra $A=\kk\langle x,y\rangle$ is graded, $z$ is homogeneous of degree $\gamma=e$, and $\sigma_3,\delta_3$ preserve both components. Put $w=y^2$. Since $yx=-xy$, the element $w$ is central in $A$, and $A_e=\kk[x,w]$ is a polynomial ring with $\deg w=2$. So $\mathbb Z_2$ is a dual reflection group for $A$ with $f_\tau=y$, and by Theorem \ref{YAMb4Lx9} also for $B$, with $B_e=A_e[z;\sigma_3|_{A_e},\delta_3|_{A_e}]$. We have $\sigma_3(w)=w$, and using $zy=-yz-2xy$ twice,
\[
zw=(zy)y=-y(zy)-2xyy=-y(-yz-2xy)-2xw=wz+2yxy-2xw=wz-4xw.
\]
Hence
\begin{equation}\label{YAMe12Br7}
B_e=\kk[x,w][z;\delta_e],\quad\delta_e(x)=x^2-w,\quad\delta_e(w)=-4xw,
\end{equation}
with $\deg x=\deg z=1$ and $\deg w=2$. Its Hilbert series is $(1-s)^{-2}(1-s^2)^{-1}$, and $\ell_{B_e}=4=\ell_B+\deg f_\tau$. Here $(A_e)_1=\kk x\neq0$, as predicted by Theorem \ref{YAMb4Lx9}, and in \eqref{YAMe36Et1} we have $\lambda_\tau=-1$ and $\delta_3(y)=y\cdot(2x)$, so $a_\tau=2x$ and $\eta_\tau(z)=-z+2x$.
\end{example}

\begin{example}\label{YAMex3Od5}
Let $\alpha\in\kk^\times$ and $A=\kk_{-1}[x,y]=\kk\langle x,y\rangle/(yx+xy)$. Let $\sigma$ be the graded automorphism with $\sigma(x)=\alpha x$ and $\sigma(y)=-y$, and let $\delta$ be the $\sigma$-derivation of the free algebra with
\[
\delta(x)=y^2,\quad\delta(y)=0.
\]
By \eqref{YAMe7Sd7}, $\delta(yx+xy)=\sigma(y)\delta(x)+\delta(y)x+\sigma(x)\delta(y)+\delta(x)y=-y^3+y^3=0$, so $\delta$ induces a $\sigma$-derivation of $A$ raising degrees by one. This derivation is outer. If $\delta=\delta_u$, then also $\delta=\delta_{u_1}$ for the component $u_1$ of $u$ of degree one, since $\delta$ raises degrees by one. For $u=px+qy\in A_1$,
\[
\delta_u(x)=ux-\sigma(x)u=p(1-\alpha)x^2-q(1+\alpha)xy,
\]
which is never $y^2$. Put $B=A[z;\sigma,\delta]$, that is,
\[
B=\kk\langle x,y,z\rangle/\big(yx+xy,\ zx-\alpha xz-y^2,\ zy+yz\big).
\]
Grade $B$ by $\mathbb Z_2=\{e,\tau\}$ with $x,z\in B_e$ and $y\in B_\tau$. This grading is compatible with the Ore structure, with $\gamma=e$. Since $A_e=\kk[x,w]$ with $w=y^2$ central in $A$, the group $\mathbb Z_2$ is a dual reflection group for $A$, and by Theorem \ref{YAMb4Lx9} also for $B$. We have $\sigma(w)=\sigma(y)^2=w$ and $\delta(w)=\sigma(y)\delta(y)+\delta(y)y=0$, so
\[
B_e=\kk[x,w][z;\sigma|_{A_e},\delta|_{A_e}],\quad zx=\alpha xz+w,\quad zw=wz.
\]
The algebra $B_e$ is AS regular of dimension three with generators of degrees $1,2,1$, so $H_{B_e}(s)=(1-s)^{-2}(1-s^2)^{-1}$ and $\ell_{B_e}=4=\ell_B+\deg f_\tau$, with $f_\tau=y$, $m_G=\tau$ and $(A_e)_1=\kk x\neq0$. In \eqref{YAMe36Et1} we have $\lambda_\tau=-1$ and $a_\tau=0$, so $\eta_\tau(z)=-z$.
\end{example}

\begin{example}\label{YAMex9Jp4}
Let $J=\kk[x][y;\delta]$ be the Jordan plane, with $\delta(x)=x^2$. Every grading of $J$ compatible with the Ore structure has $x\in J_c$ and $y\in J_\gamma$, and $\delta(x)=x^2$ forces $\gamma c=c^2$, that is, $\gamma=c$. Such a grading exists for every $c$ of finite order, since $\delta(x^k)=kx^{k+1}$ has degree $c^{k+1}=\gamma c^k$. By Theorem \ref{YAMb4Lx9}, a dual reflection group requires $\gamma=e$, hence $c=e$ and $G=\{e\}$. So the Jordan plane has no nontrivial compatible dual reflection group, although it has nontrivial compatible gradings. The hypothesis that $G$ be a dual reflection group is therefore needed for the conclusion $\gamma=e$.
\end{example}

Once $A$ is a graded subalgebra, a homogeneous Ore variable can be chosen after an affine change of variable.

\begin{proposition}\label{YAMu9Cj4}
Let $B=A[t;\sigma,\delta]$ be a graded Ore extension with $\deg t=1$, and suppose that $B$ has a $G$-grading for which $A$ is a graded subalgebra. Then there is a homogeneous element $t'=\lambda t+a$ with $\lambda\in\kk^\times$ and $a\in A_1$, and
\begin{equation}\label{YAMe37Tp2}
B=A[t';\sigma,\delta'],\quad\delta'=\lambda\delta+\delta_a.
\end{equation}
In particular the grading is compatible with the Ore structure $A[t';\sigma,\delta']$, and $\delta'$ is outer if and only if $\delta$ is outer.
\end{proposition}

\begin{proof}
Since $B_1=A_1\oplus\kk t$ and $A_1$ is a graded subspace of $B_1$, the quotient $B_1/A_1$ is a one-dimensional graded space, concentrated in some degree $\gamma$. An element of $B_1\cap B_\gamma$ with nonzero image in $B_1/A_1$ has the form $t'=\lambda t+a$ with $\lambda\neq0$ and $a\in A_1$. For $b\in A$,
\[
t'b=\lambda\sigma(b)t+\lambda\delta(b)+ab=\sigma(b)t'+\lambda\delta(b)+\big(ab-\sigma(b)a\big)=\sigma(b)t'+(\lambda\delta+\delta_a)(b).
\]
Expanding $(\lambda t+a)^n$ with \eqref{YAMe9Or9} gives $t'^n\in\lambda^nt^n+\sum_{j<n}At^j$, so $B=\bigoplus_nAt'^n$. This proves \eqref{YAMe37Tp2}. Since $\delta'-\lambda\delta=\delta_a$ is inner, $\delta'$ is outer if and only if $\delta$ is.
\end{proof}

\begin{example}\label{YAMex18Af2}
The Ore variable need not be homogeneous even when $A$ is graded. Let $B=\kk[x][t]=\kk[x,t]$, and grade it by $\mathbb Z_2=\{e,\tau\}$ with $x\in B_\tau$ and $u=t+x\in B_e$. Since $B=\kk[x,u]$, this is a grading of $B$, and $A=\kk[x]$ is a graded subalgebra, but $t=u-x$ is not homogeneous. Proposition \ref{YAMu9Cj4} replaces $t$ by $t'=u$, with $\delta'=\delta_x=0$, and the grading is compatible with $B=\kk[x][u]$. Here $\gamma=e$, $B_e=\kk[x^2,u]$, and $\mathbb Z_2$ is a dual reflection group for $B$.
\end{example}

The choice of $t'$ in Proposition \ref{YAMu9Cj4} is unique up to adding an element of $A_1\cap A_\gamma$ and rescaling, which changes $\delta'$ by an inner $\sigma$-derivation. So the hypothesis $\delta\neq0$ in Theorem \ref{YAMb4Lx9} refers to a homogeneous choice of Ore variable, and it holds for every such choice when $\delta$ is outer. The remaining part of the compatibility assumption is that $A$ be a graded subalgebra.

\begin{remark}\label{YAMv3Ei8}
For the three-dimensional AS regular algebras of Ore type in \cite[Table 1]{BR}, graded by nonabelian groups, the grading is compatible with an Ore structure after a graded change of variables. By \cite[Theorem 3.4]{BR} the grading corresponds to a quotient of the maximal grading group of the corresponding row of \cite[Table 1]{BR}, for which the basis $u,v,z$ is homogeneous, and for the Ore-type rows this basis gives an Ore decomposition with $A$ graded and $z$ homogeneous. In dimension four, the algebra with a dual reflection group of order sixteen constructed in \cite[\S3.1]{GKMV} is a double Ore extension \cite{ZZ} of an AS regular algebra of dimension two \cite[Proposition 3.1.1]{GKMV}, but this presentation is not homogeneous for the grading, as noted after \cite[Proposition 3.1.1]{GKMV}. So this example lies outside the scope of the present methods.
\end{remark}

\section{Iterated Ore extensions}\label{YAMs5It0}

In this section we prove Theorem \ref{YAMtc3Mn8} and describe the identity component as an iterated Ore extension. Let $B$ be the iterated Ore extension \eqref{YAMq8Hi1}, and fix a $G$-grading of $B$ in which every $x_i$ is homogeneous. As explained in Section \ref{YAMs1Pr0}, this grading is compatible with each Ore extension $B^{(i)}=B^{(i-1)}[x_i;\sigma_i,\delta_i]$. We write $c_i\in G$ for the $G$-degree of $x_i$, we put $r_i=|c_i|$ and $K_i=\supp(B^{(i)})$, and we grade $B^{(i)}$ by $K_i$. Its identity component is $B^{(i)}_e=B^{(i)}\cap B_e$. Since $B$ is a domain and the grading is faithful, $K_n=G$.

\begin{lemma}\label{YAMl6Kn3}
The following hold.
\begin{enumerate}
\item $K_1=\langle c_1\rangle$, and for $2\le i\le n$ the subgroup $K_{i-1}$ is normal in $K_i=K_{i-1}\langle c_i\rangle$.
\item $G=\langle c_1\rangle\langle c_2\rangle\cdots\langle c_n\rangle$ as a product of subsets, and
\[
|G|=\frac{r_1r_2\cdots r_n}{\prod_{i=2}^n|K_{i-1}\cap\langle c_i\rangle|}\le r_1r_2\cdots r_n.
\]
\item If $\delta_i\neq0$, then $c_i\in K_{i-1}$.
\item The equality $|G|=r_1\cdots r_n$ holds if and only if $K_{i-1}\cap\langle c_i\rangle=\{e\}$ for $2\le i\le n$, if and only if every $g\in G$ has a unique expression $g=c_1^{a_1}\cdots c_n^{a_n}$ with $0\le a_i<r_i$.
\end{enumerate}
\end{lemma}

\begin{proof}
The component $\kk[x_1]_\beta$ is spanned by the powers $x_1^k$ with $c_1^k=\beta$, so $K_1=\langle c_1\rangle$. For $i\ge2$, Lemma \ref{YAMl4Ae9}(1),(2) applied to $B^{(i)}=B^{(i-1)}[x_i;\sigma_i,\delta_i]$ with its $K_i$-grading give (1), and Lemma \ref{YAMl4Ae9}(3) gives (3). By induction on $i$, (1) gives $K_i=\langle c_1\rangle\cdots\langle c_i\rangle$ as a product of subsets, and $|K_i|=|K_{i-1}|\,r_i/|K_{i-1}\cap\langle c_i\rangle|$. Since $K_n=G$, this proves (2). The map $\prod_i\{0,\dots,r_i-1\}\to G$, $(a_i)\mapsto c_1^{a_1}\cdots c_n^{a_n}$, is surjective by (2), so it is bijective if and only if $|G|=r_1\cdots r_n$, which by (2) happens if and only if every intersection $K_{i-1}\cap\langle c_i\rangle$ is trivial.
\end{proof}

\begin{proof}[\textnormal{\bfseries Proof of Theorem \ref{YAMtc3Mn8}}]
By Lemma \ref{YAMl6Kn3}(2), the group $G$ is the product of the cyclic subgroups $\langle c_i\rangle$, and the equivalence of (2) and (3) is Lemma \ref{YAMl6Kn3}(4).

Assume (1). We prove by descending induction on $i$ that $K_i$ is a dual reflection group for $B^{(i)}$ and that $K_{i-1}\cap\langle c_i\rangle=\{e\}$. For $i=n$ the first statement is the hypothesis. If $K_i$ is a dual reflection group for $B^{(i)}=B^{(i-1)}[x_i;\sigma_i,\delta_i]$, then Theorem \ref{YAMb4Lx9} shows that $K_{i-1}$ is a dual reflection group for $B^{(i-1)}$ and that $K_{i-1}\cap\langle c_i\rangle=\{e\}$. Hence (2) holds by Lemma \ref{YAMl6Kn3}(4).

Assume (2). By Lemma \ref{YAMl6Kn3}(4), every intersection $K_{i-1}\cap\langle c_i\rangle$ is trivial. If $\delta_i\neq0$, then $c_i\in K_{i-1}$ by Lemma \ref{YAMl6Kn3}(3), so $\langle c_i\rangle\subseteq K_{i-1}\cap\langle c_i\rangle$ and $c_i=e$. We prove by induction on $i$ that $K_i$ is a dual reflection group for $B^{(i)}$. For $i=1$, the identity component of $\kk[x_1]$ for its $\langle c_1\rangle$-grading is $\kk[x_1^{r_1}]$, which is AS regular. For $i\ge2$, Theorem \ref{YAMb4Lx9} applies to $B^{(i)}=B^{(i-1)}[x_i;\sigma_i,\delta_i]$, since $K_{i-1}$ is a dual reflection group for $B^{(i-1)}$ and $K_{i-1}\cap\langle c_i\rangle=\{e\}$. For $i=n$ this is (1).

Assume now that (1), (2) and (3) hold. We saw that $c_i=e$ whenever $\delta_i\neq0$. By Theorem \ref{YAMb4Lx9}, $K_i=K_{i-1}\rtimes\langle c_i\rangle$, so the series
\[
\{e\}\trianglelefteq K_1\trianglelefteq K_2\trianglelefteq\cdots\trianglelefteq K_n=G,
\]
has cyclic factors $K_i/K_{i-1}\cong\langle c_i\rangle$, and $G$ is solvable. Applying Proposition \ref{YAMn5Fh3} to each step, with $f_{c_1^{a}}=x_1^{a}$ for $0\le a<r_1$ in $\kk[x_1]$, we get by induction on $i$
\begin{equation}\label{YAMe78Nf6}
f_{c_1^{a_1}\cdots c_i^{a_i}}=x_1^{a_1}\cdots x_i^{a_i}\quad(0\le a_j<r_j),
\end{equation}
for the grading of $B^{(i)}$ by $K_i$. The step with $\delta_i\neq0$ has $r_i=1$ and does not change the elements $f_g$. By \eqref{YAMe45Pg1} and induction,
\begin{equation}\label{YAMe79Nf7}
\begin{aligned}
p_G(s)&=\prod_{i=1}^n[r_i]_s,&m_G&=c_1^{r_1-1}\cdots c_n^{r_n-1},\\
N_G&=\sum_{i=1}^n(r_i-1),&f_{m_G}&=x_1^{r_1-1}\cdots x_n^{r_n-1}.
\end{aligned}
\end{equation}
Since $H_B(s)=(1-s)^{-n}$ by \eqref{YAMe12Hs3}, the identity \eqref{YAMh3Nc7} gives
\[
H_{B_e}(s)=\frac{1}{(1-s)^n\prod_i[r_i]_s}=\prod_{i=1}^n\frac{1}{1-s^{r_i}}.
\]

Suppose finally that $(B_e)_1=0$. No $x_i$ lies in $B_e$, so $c_i\neq e$ for all $i$, and therefore $\delta_i=0$ for $i\ge2$. By \eqref{YAMe3Lf5}, each component $B_1\cap B_g$ with $g\neq e$ is at most one-dimensional, and $B_1\cap B_e=0$. Since $x_1,\dots,x_n$ are homogeneous and linearly independent, their $G$-degrees are pairwise distinct, $B_1\cap B_{c_j}=\kk x_j$, and the other components of $B_1$ vanish. Let $j<i$. Since $(\sigma_i,0)$ is $c_i$-compatible, $\sigma_i(x_j)$ is a nonzero element of $(B^{(i-1)})_1\cap B_{c_ic_jc_i^{-1}}$, so $c_ic_jc_i^{-1}=c_{\pi_i(j)}$ for some $\pi_i(j)<i$, and $\sigma_i(x_j)=\lambda_{ij}x_{\pi_i(j)}$ with $\lambda_{ij}\neq0$. Since $\sigma_i$ is injective on $(B^{(i-1)})_1$, the map $\pi_i$ is a permutation of $\{1,\dots,i-1\}$. The binomial presentation follows from \eqref{YAMe13Pr4}.
\end{proof}

Theorem \ref{YAMtc3Mn8} turns the search for compatible dual reflection groups into a problem about finite quotients of one finitely presented group. A grading of $B$ with $x_i\in B_{c_i}$ exists if and only if the relations \eqref{YAMe13Pr4} are homogeneous. Since the monomials $x_1^{a_1}\cdots x_n^{a_n}$ are linearly independent in $B$, this means that $\sigma_i(x_j)$ is a combination of the $x_k$ with $c_k=c_ic_jc_i^{-1}$, and that $\delta_i(x_j)$, written in the monomial basis, is a combination of monomials $x_kx_l$ with $k\le l$ and $c_kc_l=c_ic_j$. Let $\Gamma_B$ be the group generated by symbols $\hat c_1,\dots,\hat c_n$ subject to the relations
\[
\hat c_k=\hat c_i\hat c_j\hat c_i^{-1},\quad\hat c_k\hat c_l=\hat c_i\hat c_j,
\]
for all $1\le j<i\le n$, all variables $x_k$ occurring in $\sigma_i(x_j)$ and all monomials $x_kx_l$ occurring in $\delta_i(x_j)$ with nonzero coefficient. Then the gradings of $B$ in which every $x_i$ is homogeneous are exactly the gradings induced by the surjections $\Gamma_B\to G$ onto finite groups. By Theorem \ref{YAMtc3Mn8}, the dual reflection groups among them are exactly the quotients in which the order of $G$ is the product of the orders of the images of the $\hat c_i$.

\begin{example}\label{YAMex7Eu1}
Let $n\ge2$, let $A=\kk[x_1,\dots,x_n]$, and let $\delta=x_1E$, where $E=\sum_ix_i\frac{\partial}{\partial x_i}$ is the Euler derivation. So $\delta(x_1)=x_1^2$ and $\delta(x_i)=x_1x_i$, and $\delta$ raises degrees by one. The Ore extension
\[
B=A[t;\delta]=\kk\langle x_1,\dots,x_n,t\rangle/\big(x_ix_j-x_jx_i,\ tx_1-x_1t-x_1^2,\ tx_i-x_it-x_1x_i\ (i\ge2)\big)
\]
is AS regular of dimension $n+1$, and $t$ is not normal, since $Bt=\bigoplus_{k\ge1}At^k$ does not contain $tx_1=x_1t+x_1^2$. The Jordan plane is the analogous algebra with $n=1$. Let $r_2,\dots,r_n\ge1$ and grade $B$ by $G=\prod_{i\ge2}\mathbb Z_{r_i}$, with $x_1$ and $t$ of degree $e$ and $x_i$ of degree a generator $c_i$ of the factor $\mathbb Z_{r_i}$. The relations are homogeneous since $\delta(x_i)=x_1x_i$ has degree $c_i$. Since $|G|=\prod_{i\ge2}r_i$, Theorem \ref{YAMtc3Mn8} shows that $G$ is a dual reflection group for $B$, with $x_1,t\in B_e$, and by Theorem \ref{YAMb4Lx9}
\[
B_e=\kk[x_1,y_2,\dots,y_n][t;\delta_e],\quad y_i=x_i^{r_i},\quad\delta_e(x_1)=x_1^2,\quad\delta_e(y_i)=r_ix_1y_i.
\]
Here $x_1\in(A_e)_1$, so $(A_e)_1\neq0$, as predicted by Theorem \ref{YAMb4Lx9}, and $\delta$ induces the zero map on $C(A)=\bigotimes_{i\ge2}\kk[\bar x_i]/(\bar x_i^{\,r_i})$, as in Proposition \ref{YAMet9Lw3}, since $\delta(x_i)=x_1x_i\in A(A_e)_{\ge1}$. By Theorem \ref{YAMtd9Ci4}, $C(B)\cong C(A)$ and $\gci(C(B))=\#\{i\ge2\mid r_i\ge2\}$. Moreover $f_{m_G}=x_2^{r_2-1}\cdots x_n^{r_n-1}$ and $\ell_{B_e}=2+\sum_{i\ge2}r_i$.
\end{example}

\begin{remark}\label{YAMrm6Br1}
By \cite[Theorem 0.3(3)]{KKZ}, $\deg f_g$ is the length of $g$ with respect to the set of degrees of the homogeneous elements of $B_1\setminus B_e$. Here this set is $\{c_i\mid c_i\neq e\}$, and \eqref{YAMe78Nf6} shows that $\deg f_g=a_1+\cdots+a_n$. So the normal form of Theorem \ref{YAMtc3Mn8}(3) is a word of minimal length representing $g$.
\end{remark}

The proof also describes the identity component as an iterated Ore extension.

\begin{proposition}\label{YAMp7Ie3}
In the setting of Theorem \ref{YAMtc3Mn8}, suppose that $G$ is a dual reflection group for $B$, and put $y_i=x_i^{r_i}$. Then $B^{(i)}_e$ is generated by $y_1,\dots,y_i$, and
\begin{equation}\label{YAMe77Ie5}
B_e=\kk[y_1][y_2;\sigma'_2,\delta'_2]\cdots[y_n;\sigma'_n,\delta'_n],\quad\deg y_i=r_i,
\end{equation}
where $\sigma'_i=\sigma_i^{r_i}|_{B^{(i-1)}_e}$ and $\delta'_i=\delta_i|_{B^{(i-1)}_e}$, with $\delta_i=0$ whenever $r_i>1$. Moreover
\begin{equation}\label{YAMe77Ie6}
\ell_{B_e}=r_1+r_2+\cdots+r_n=\ell_B+N_G.
\end{equation}
\end{proposition}

\begin{proof}
By the proof of Theorem \ref{YAMtc3Mn8}, $K_i$ is a dual reflection group for $B^{(i)}$ and $K_{i-1}\cap\langle c_i\rangle=\{e\}$. So \eqref{YAMe34Be8} gives $B^{(i)}_e=B^{(i-1)}_e[y_i;\sigma'_i,\delta'_i]$. Since $B^{(1)}_e=\kk[y_1]$, induction on $i$ gives \eqref{YAMe77Ie5}. Since $\kk$ has Gorenstein parameter $0$, Corollary \ref{YAMc1Rg5} gives $\ell_{B^{(i)}_e}=\ell_{B^{(i-1)}_e}+r_i$, hence $\ell_{B_e}=\sum_ir_i$, and $\ell_B=n$ and $N_G=\sum_i(r_i-1)$ by \eqref{YAMe79Nf7}.
\end{proof}

The identity \eqref{YAMe77Ie6} is an instance of \cite[Lemma 3.4(3)]{KKZ}, since the Frobenius algebra $C(B)$ has Gorenstein parameter $-N_G$.

\begin{corollary}\label{YAMc7Qa1}
Let $(q_{ij})$ be a multiplicatively antisymmetric matrix in $\kk^\times$, that is, $q_{ii}=1$ and $q_{ji}=q_{ij}^{-1}$, and let
\begin{equation}\label{YAMe68Qa2}
\mathcal O_q(\kk^n)=\kk\langle x_1,\dots,x_n\rangle/\big(x_ix_j-q_{ij}x_jx_i\big)
\end{equation}
be the quantum affine space, viewed as the iterated Ore extension with $\sigma_i(x_j)=q_{ij}x_j$ for $j<i$ and $\delta_i=0$. Let $G$ be a finite group with a grading of $\mathcal O_q(\kk^n)$ in which every $x_i$ is homogeneous, of degree $c_i$. Then the $c_i$ commute, and $G$ is a dual reflection group if and only if the multiplication map $\prod_i\langle c_i\rangle\to G$ is an isomorphism. In this case
\[
B_e\cong\mathcal O_{q'}(\kk^n),\quad q'_{ij}=q_{ij}^{\,r_ir_j},
\]
on the generators $x_1^{r_1},\dots,x_n^{r_n}$.
\end{corollary}

\begin{proof}
Since $\sigma_i(x_j)=q_{ij}x_j$ is a nonzero element of $B_{c_ic_jc_i^{-1}}\cap B_{c_j}$, the elements $c_i$ and $c_j$ commute. So $G=\langle c_1\rangle\cdots\langle c_n\rangle$ is abelian, the multiplication map $\prod_i\langle c_i\rangle\to G$ is a surjective homomorphism, and it is an isomorphism if and only if $|G|=\prod_ir_i$. The first statement is now Theorem \ref{YAMtc3Mn8}. By Proposition \ref{YAMp7Ie3}, $B_e$ is generated by the elements $x_i^{r_i}$, and the relations \eqref{YAMe68Qa2} give, by induction on the exponents,
\[
x_i^ax_j^b=q_{ij}^{ab}\,x_j^bx_i^a,\quad\text{so}\quad x_i^{r_i}x_j^{r_j}=q_{ij}^{r_ir_j}x_j^{r_j}x_i^{r_i}.
\]
The monomials in the $x_i^{r_i}$ are linearly independent, so these are the defining relations of $B_e$.
\end{proof}

In the dual language, the $G$-grading of $\mathcal O_q(\kk^n)$ is a diagonal action of the character group $G^\vee$, and Corollary \ref{YAMc7Qa1} is the diagonal case of the theorem of Kirkman, Kuzmanovich and Zhang for skew polynomial rings \cite[Theorem 4.5]{KKZ0}.

In dimension two, Theorem \ref{YAMtc3Mn8} applies to the quantum plane and the Jordan plane. For $\kk_q[x,y]=\kk[x][y;\sigma]$, the $G$-degrees $c_1$ of $x$ and $c_2$ of $y$ commute for every compatible grading, and $G$ is a dual reflection group if and only if $\langle c_1\rangle\cap\langle c_2\rangle=\{e\}$. Then $G=\mathbb Z_{r_1}\times\mathbb Z_{r_2}$, $B_e=\kk_{q^{r_1r_2}}[x^{r_1},y^{r_2}]$ and $f_{m_G}=x^{r_1-1}y^{r_2-1}$. For the Jordan plane, the only compatible dual reflection group is trivial (Example \ref{YAMex9Jp4}). These gradings are abelian. For the quantum plane they are among the coactions classified by Buzaglo and Rogalski \cite[Theorem 1.2]{BR2}.

\begin{example}\label{YAMc12Dt0}
We describe the three-dimensional case with $(B_e)_1=0$. Let $B$ be an iterated Ore extension \eqref{YAMq8Hi1} with $n=3$, and let $G$ be a dual reflection group for $B$ with every $x_i$ homogeneous and $(B_e)_1=0$. By Theorem \ref{YAMtc3Mn8}, $\delta_2=\delta_3=0$, the degrees $c_1,c_2,c_3$ are distinct and different from $e$, and the relations are binomial. The permutation $\pi_2$ of $\{1\}$ is trivial, so $x_2x_1=qx_1x_2$ for some $q\in\kk^\times$, and $c_1,c_2$ commute. If $\pi_3$ is the identity, then $G$ is abelian and $B$ is a quantum affine space as in Corollary \ref{YAMc7Qa1}. Otherwise $\pi_3$ is the transposition, and
\begin{equation}\label{YAMe72Dt1}
B=\kk\langle x_1,x_2,x_3\rangle/\big(x_2x_1-qx_1x_2,\ x_3x_1-\lambda x_2x_3,\ x_3x_2-\mu x_1x_3\big)\quad(\lambda,\mu\in\kk^\times).
\end{equation}
The assignment $\sigma_3(x_1)=\lambda x_2$, $\sigma_3(x_2)=\mu x_1$ defines an endomorphism of $\kk_q[x_1,x_2]$ if and only if
\[
\sigma_3(x_2)\sigma_3(x_1)-q\sigma_3(x_1)\sigma_3(x_2)=\lambda\mu\,(x_1x_2-qx_2x_1)=\lambda\mu(1-q^2)\,x_1x_2
\]
vanishes, that is, if and only if $q=\pm1$. Conjugation by $c_3$ exchanges $c_1$ and $c_2$, so $r_1=r_2=:r\ge2$, conjugation by $c_3$ is an automorphism of order two of $\langle c_1\rangle\times\langle c_2\rangle$, and $r_3=2m$ is even. By Theorem \ref{YAMtc3Mn8},
\[
G=\big(\langle c_1\rangle\times\langle c_2\rangle\big)\rtimes\langle c_3\rangle\cong(\mathbb Z_r\times\mathbb Z_r)\rtimes\mathbb Z_{2m},
\]
where the generator of $\mathbb Z_{2m}$ exchanges the two factors, and $c_3^2$ is central. Conversely, for $q=\pm1$, $\lambda,\mu\in\kk^\times$, $r\ge2$ and $m\ge1$, the relations \eqref{YAMe72Dt1} are homogeneous for the grading of $B$ by this group with $x_i\in B_{c_i}$ for $i=1,2,3$, because $c_3c_1=c_2c_3$ and $c_3c_2=c_1c_3$. Since $|G|=r\cdot r\cdot2m$, Theorem \ref{YAMtc3Mn8} shows that $G$ is a dual reflection group, and Proposition \ref{YAMp7Ie3} gives
\[
B_e=\kk\langle y_1,y_2,y_3\rangle/\big(y_2y_1-q^{r^2}y_1y_2,\ y_3y_1-(\lambda\mu)^{mr}y_1y_3,\ y_3y_2-(\lambda\mu)^{mr}y_2y_3\big),
\]
with $y_1=x_1^r$, $y_2=x_2^r$ and $y_3=x_3^{2m}$, since $\sigma_3^{2m}(x_j)=(\lambda\mu)^mx_j$ for $j=1,2$. This recovers, for the algebras of Ore type, the family of groups $(\mathbb Z_r\times\mathbb Z_r)\rtimes\mathbb Z_{2m}$ in the classification of Buzaglo and Rogalski \cite[Theorem 1.2]{BR}, together with an explicit identity component.
\end{example}

The following family was constructed by Goetz, Kirkman, Moore and Vashaw \cite[Example 2.1.7]{GKMV}, who also allow $\kk_{-1}[x_1,\dots,x_n]$ in place of $\kk[x_1,\dots,x_n]$. We treat the commutative case, and the same argument applies to the skew case. Its dual reflection groups follow from Theorems \ref{YAMb4Lx9} and \ref{YAMtc3Mn8}.

\begin{example}\label{YAMex7Wr3}
Let $n\ge2$, $r\ge1$ and $p\ge1$. Let $A=\kk[x_1,\dots,x_n]$, let $\sigma$ be the automorphism with $\sigma(x_i)=x_{i+1}$, indices modulo $n$, and put
\begin{equation}\label{YAMe66Wr8}
B=A[t;\sigma]=\kk\langle x_1,\dots,x_n,t\rangle/\big(x_ix_j-x_jx_i,\ tx_i-x_{i+1}t\big).
\end{equation}
It is the iterated Ore extension $\kk[x_1][x_2]\cdots[x_n][t;\sigma]$, AS regular of dimension $n+1$. Let $\varepsilon_1,\dots,\varepsilon_n$ be generators of $\mathbb Z_r^n=\langle\varepsilon_1\rangle\times\cdots\times\langle\varepsilon_n\rangle$, written multiplicatively, and let
\begin{equation}\label{YAMe67Wg9}
G=\mathbb Z_r^n\rtimes\langle\gamma\rangle,\quad|\gamma|=np,\quad\gamma\varepsilon_i\gamma^{-1}=\varepsilon_{i+1}.
\end{equation}
For $p=1$ this is the wreath product $\mathbb Z_r\wr\mathbb Z_n$. Grade $B$ by $x_i\in B_{\varepsilon_i}$ and $t\in B_\gamma$. The automorphism $\sigma$ maps a monomial of $G$-degree $g=\varepsilon_1^{a_1}\cdots\varepsilon_n^{a_n}$ to one of $G$-degree $\varepsilon_2^{a_1}\cdots\varepsilon_1^{a_n}=\gamma g\gamma^{-1}$, so the relations are homogeneous. Here $K=\mathbb Z_r^n$ is a dual reflection group for $A$ with $A_e=\kk[x_1^r,\dots,x_n^r]$, and $K\cap\langle\gamma\rangle=\{e\}$. Equivalently, $|G|=r^n\cdot np=|\varepsilon_1|\cdots|\varepsilon_n|\,|\gamma|$. By Theorem \ref{YAMb4Lx9}, $G$ is a dual reflection group for $B$, with
\[
B_e=\kk[x_1^r,\dots,x_n^r][t^{np};\sigma^{np}]=\kk[x_1^r,\dots,x_n^r,t^{np}],
\]
since $\sigma^n=\mathrm{id}$. By Theorem \ref{YAMtc3Mn8} and Proposition \ref{YAMp7Ie3},
\[
p_G(s)=[r]_s^n\,[np]_s,\quad|G|=npr^n,\quad f_{m_G}=(x_1\cdots x_n)^{r-1}t^{np-1},\quad\ell_{B_e}=nr+np.
\]
The group $G$ is nonabelian if and only if $r\ge2$. When $r\ge2$ we have $(B_e)_1=0$, and the defining relations \eqref{YAMe66Wr8} are binomial, as predicted by Theorem \ref{YAMtc3Mn8}.
\end{example}

\begin{corollary}\label{YAMc7Wr4}
Let $B=\kk[x_1,\dots,x_n][t;\sigma]$ be as in Example \ref{YAMex7Wr3}, and let $G$ be a finite group with a grading of $B$ in which $x_1,\dots,x_n$ and $t$ are homogeneous. Then $G$ is a dual reflection group for $B$ if and only if either $x_1,\dots,x_n\in B_e$ and $G$ is cyclic of arbitrary order, generated by the $G$-degree of $t$, or $G\cong\mathbb Z_r^n\rtimes\mathbb Z_{np}$ as in \eqref{YAMe67Wg9} for some $r\ge2$ and $p\ge1$, with the grading of Example \ref{YAMex7Wr3}.
\end{corollary}

\begin{proof}
Let $\varepsilon_i$ and $\gamma$ be the $G$-degrees of $x_i$ and $t$. The relations $x_ix_j=x_jx_i$ and $tx_i=x_{i+1}t$ are homogeneous exactly when the $\varepsilon_i$ commute and $\gamma\varepsilon_i\gamma^{-1}=\varepsilon_{i+1}$. So the $\varepsilon_i$ are conjugate and have a common order $r$, and $E=\langle\varepsilon_1,\dots,\varepsilon_n\rangle$ is abelian of exponent dividing $r$, with $|E|\le r^n$. By Lemma \ref{YAMl4Ae9}(2), $G=E\langle\gamma\rangle$ with $E$ normal, so
\[
|G|\le|E|\,|\gamma|\le r^n|\gamma|,
\]
with equality if and only if $E\cong\mathbb Z_r^n$ with basis $\varepsilon_1,\dots,\varepsilon_n$ and $E\cap\langle\gamma\rangle=\{e\}$. By Theorem \ref{YAMtc3Mn8}, $G$ is a dual reflection group exactly in this case. If $r=1$, then $G=\langle\gamma\rangle$. If $r\ge2$, the elements $\varepsilon_1,\dots,\varepsilon_n$ are pairwise distinct, and $\gamma^{|\gamma|}=e$ gives $\varepsilon_{i+|\gamma|}=\varepsilon_i$, so $n$ divides $|\gamma|$ and $G$ is the group \eqref{YAMe67Wg9} with $np=|\gamma|$. Conversely, the groups with $r\ge2$ are dual reflection groups by Example \ref{YAMex7Wr3}. If $r=1$ and $|\gamma|=N$, then $K=\{e\}$, and Theorem \ref{YAMb4Lx9} shows that $\langle\gamma\rangle$ is a dual reflection group, with $B_e=A[t^N;\sigma^N]$.
\end{proof}

\begin{example}\label{YAMex13Dt6}
We make Example \ref{YAMex7Wr3} explicit for $n=r=2$ and $p=1$. Here $B=\kk[x_1,x_2][t;\sigma]$ with $\sigma(x_1)=x_2$ and $\sigma(x_2)=x_1$, that is,
\begin{equation}\label{YAMe66Wr9}
x_2x_1=x_1x_2,\quad tx_1=x_2t,\quad tx_2=x_1t.
\end{equation}
These are the binomial relations of Theorem \ref{YAMtc3Mn8}, with $\pi_3$ the transposition of $\{1,2\}$. The group is $G=(\langle\varepsilon_1\rangle\times\langle\varepsilon_2\rangle)\rtimes\langle\gamma\rangle$ with $\varepsilon_1^2=\varepsilon_2^2=\gamma^2=e$ and $\gamma\varepsilon_1\gamma^{-1}=\varepsilon_2$. The element $\rho=\varepsilon_1\gamma$ satisfies $\rho^2=\varepsilon_1\gamma\varepsilon_1\gamma^{-1}=\varepsilon_1\varepsilon_2$, so it has order four, and $\gamma\rho\gamma^{-1}=\varepsilon_2\gamma=\rho^{-1}$. Since $\gamma\notin\langle\rho\rangle$, the group $G=\langle\rho\rangle\rtimes\langle\gamma\rangle$ is the dihedral group of order eight. By \eqref{YAMe78Nf6} the elements $f_g$ are
\[
\begin{array}{c|cccccccc}
g&e&\varepsilon_1&\varepsilon_2&\varepsilon_1\varepsilon_2&\gamma&\varepsilon_1\gamma&\varepsilon_2\gamma&\varepsilon_1\varepsilon_2\gamma\\\hline
f_g&1&x_1&x_2&x_1x_2&t&x_1t&x_2t&x_1x_2t\\
\deg f_g&0&1&1&2&1&2&2&3
\end{array}
\]
So $p_G(s)=(1+s)^3$, $m_G=\varepsilon_1\varepsilon_2\gamma$ and $f_{m_G}=x_1x_2t$. The identity component $B_e=\kk[x_1^2,x_2^2,t^2]$ is a polynomial ring with Hilbert series $(1-s^2)^{-3}$ and Gorenstein parameter $6=\ell_B+\deg f_{m_G}$. The identity \eqref{YAMe27Fm1} can be checked directly. For instance $f_{\varepsilon_1}f_{\varepsilon_1^{-1}m_G}=x_1\cdot x_2t=f_{m_G}$ and $f_\gamma f_{\gamma^{-1}m_G}=t\cdot x_1x_2=x_2x_1t=f_{m_G}$, since $\gamma^{-1}\varepsilon_1\varepsilon_2\gamma=\varepsilon_1\varepsilon_2$.
\end{example}

\section{Covariant rings and complete intersections}\label{YAMs6Cv0}

In this section we determine the structure of the covariant rings of iterated Ore extensions and the growth of their Ext-algebras, and we prove Theorem \ref{YAMtd9Ci4}.

\begin{proposition}\label{YAMcv1Ka7}
In the setting of Theorem \ref{YAMb4Lx9}, suppose that $G$ is a dual reflection group for $B$. Then $\sigma$ induces a graded automorphism $\bar\sigma$ of $C(A)$, and the inclusion $A\subseteq B$ induces an isomorphism of graded algebras
\[
C(B)\cong C(A)[\bar t;\bar\sigma]/(\bar t^{\,r}).
\]
In particular $C(B)\cong C(A)$ if $\delta\neq0$.
\end{proposition}

\begin{proof}
Put $J=A(A_e)_{\ge1}$, a two-sided ideal of $A$ by \eqref{YAMt4Kv0}. Since $\sigma(A_e)=A_e$, we have $\sigma(J)=J$, so $\sigma$ induces $\bar\sigma$ on $C(A)=A/J$. By \eqref{YAMe34Be8}, $(B_e)_{\ge1}=(A_e)_{\ge1}\oplus\bigoplus_{l\ge1}A_et^{rl}$, so
\[
B(B_e)_{\ge1}=B(A_e)_{\ge1}+Bt^r,\quad Bt^r=\bigoplus_{n\ge r}At^n.
\]
We claim that
\begin{equation}\label{YAMe44Cb8}
B(B_e)_{\ge1}=\bigoplus_{n<r}Jt^n\oplus\bigoplus_{n\ge r}At^n.
\end{equation}
The right side is contained in the left side, since $Jt^n=A(A_e)_{\ge1}t^n$ and $(A_e)_{\ge1}t^n=t^n\sigma^{-n}\big((A_e)_{\ge1}\big)$ when $\delta=0$, while $r=1$ when $\delta\neq0$. For the reverse inclusion it suffices to show that $At^nb$ lies in the right side for $b\in(A_e)_{\ge1}$ and $n\ge0$. If $\delta=0$, then $At^nb=A\sigma^n(b)t^n\subseteq Jt^n$. If $\delta\neq0$, then $r=1$, and \eqref{YAMe10Tn1} gives $t^nb\in\delta^n(b)+Bt$ for $n\ge1$, with $\delta^n(b)\in(A_e)_{\ge1}$ since $\delta$ preserves $A_e$ and raises degrees. So $At^nb\subseteq Jt^n$ if $\delta=0$, and $At^nb\subseteq J+Bt$ if $\delta\neq0$, where $r=1$. In both cases $At^nb$ lies in the right side of \eqref{YAMe44Cb8}, which proves it. The map $B=\bigoplus_nAt^n\to C(A)[\bar t;\bar\sigma]/(\bar t^{\,r})$, $at^n\mapsto\bar a\bar t^n$, is a surjective homomorphism of graded algebras by the universal property of Ore extensions, because $\bar\delta=0$ on $C(A)$ when $\delta\neq0$ by Proposition \ref{YAMet9Lw3}. Its kernel is the right side of \eqref{YAMe44Cb8}.
\end{proof}

For example, in Example \ref{YAMex3Od5} we have $\delta\neq0$, and
\[
B(B_e)_{\ge1}=A(A_e)_{\ge1}\oplus Bz,\quad A(A_e)_{\ge1}=Ax+Ay^2,\quad C(B)\cong C(A)=\kk[\bar y]/(\bar y^{\,2}),
\]
with Hilbert series $1+s=p_G(s)$. The Ore variable $z$ and the generator $x$ of $A_e$ both disappear in the covariant ring.

Let $C'$ be a finite-dimensional connected graded algebra, let $\bar\sigma$ be a graded automorphism of $C'$, let $r\ge2$, and put
\[
C=C'[x;\bar\sigma]/(x^r),\quad\deg x=1.
\]
The element $x$ is normal in $C'[x;\bar\sigma]$, so $(x^r)=x^rC'[x;\bar\sigma]$, and $C=\bigoplus_{j<r}C'x^j=\bigoplus_{j<r}x^jC'$ is free as a left and as a right $C'$-module. The ideal $xC=Cx$ is two-sided with $C/xC\cong C'$, and we write $\pi\colon C\to C'$ for the quotient map. Let $F_\bullet\to\kk$ be a minimal graded free resolution of $\kk_{C'}$. Since $C'$ is finite-dimensional, each $F_p$ has finite rank, and $\Ext^p_{C'}(\kk,\kk)\cong\Hom_{C'}(F_p,\kk)$.

Let $F^{\bar\sigma}$ be the complex $F_\bullet$ with the right action $f\ast c=f\bar\sigma(c)$. It is a free resolution of $\kk^{\bar\sigma}=\kk$, so the identity of $\kk$ lifts to a $C'$-linear chain map $F\to F^{\bar\sigma}$, that is, to a degree-preserving chain map $\theta\colon F_\bullet\to F_\bullet$ with
\begin{equation}\label{YAMe5Tr2}
\theta(fc)=\theta(f)\,\bar\sigma(c)\quad(f\in F_p,\ c\in C'),
\end{equation}
and $\theta_0(1)=1$ on $F_0=C'$. For $k\ge1$, the map
\[
\varphi_k\colon F_p\otimes_{C'}C\to F_p\otimes_{C'}C,\quad\varphi_k(f\otimes c)=\theta^k(f)\otimes x^kc,
\]
is well defined, since $\theta^k(fc')\otimes x^kc=\theta^k(f)\otimes\bar\sigma^k(c')x^kc=\theta^k(f)\otimes x^kc'c$ for $c'\in C'$. It is right $C$-linear, it raises the internal degree by $k$, and it commutes with $d_F\otimes1$. Moreover $\varphi_1\varphi_{r-1}=\varphi_{r-1}\varphi_1=\varphi_r=0$, since $x^r=0$ in $C$.

\begin{proposition}\label{YAMp3Tr7}
For $p,q\ge0$, put $D_{p,q}=F_p\otimes_{C'}C(-s_q)$, where $s_{2j}=jr$ and $s_{2j+1}=jr+1$. Let the horizontal differential $D_{p,q}\to D_{p-1,q}$ be $d_F\otimes1$, and let the vertical differential $D_{p,q}\to D_{p,q-1}$ be $\varphi_1$ for $q$ odd and $\varphi_{r-1}$ for $q\ge2$ even.
\begin{enumerate}
\item The total complex $\Tot D$, with differential $d_F\otimes1+(-1)^p\varphi$ on $D_{p,q}$, is a minimal graded free resolution of $\kk_C$.
\item The Poincar\'e series of $C$ and $C'$ are related by
\begin{equation}\label{YAMe5Tr4}
\begin{aligned}
P_C(y,s)&=P_{C'}(y,s)\,\frac{1+ys}{1-y^2s^r},\\
\dim\Ext^n_C(\kk,\kk)&=\sum_{p=0}^n\dim\Ext^p_{C'}(\kk,\kk).
\end{aligned}
\end{equation}
\item The algebra map $\pi^*\colon E(C')\to E(C)$ induced by $\pi$ is injective. There is $\eta\in\Ext^2_C(\kk,\kk)_{-r}$ such that the linear map
\[
E(C')\otimes\kk[\eta]\to E(C),\quad a\otimes\eta^j\mapsto\pi^*(a)\,\eta^j,
\]
is injective, with image $\bigoplus_{j\ge0}\Hom_C(\bigoplus_pD_{p,2j},\kk)$.
\end{enumerate}
\end{proposition}

\begin{proof}
(1) The vertical maps commute with the horizontal ones and compose to zero, so $D$ is a double complex concentrated in the first quadrant, and the sign makes $\Tot D$ a complex. Since $C$ is free as a left $C'$-module, each row $D_{\bullet,q}$ is a free resolution of $N(-s_q)$, where
\[
N=\kk\otimes_{C'}C=C/C'_+C\cong\kk[\bar x]/(\bar x^r).
\]
Here $C'_+C=CC'_+$ is a two-sided ideal, since $xC'_+=\bar\sigma(C'_+)x=C'_+x$. On $F_0\otimes_{C'}C=C$ we have $\varphi_k(c)=x^kc$, since $\theta_0=\bar\sigma$. So the induced maps on $H_0$ of the rows are left multiplication by $\bar x$ and by $\bar x^{r-1}$, and they form the complex
\[
\cdots\to N(-r-1)\xrightarrow{\ \bar x\ }N(-r)\xrightarrow{\ \bar x^{r-1}\ }N(-1)\xrightarrow{\ \bar x\ }N\to\kk\to0.
\]
It is exact, since in $\kk[\bar x]/(\bar x^r)$ the kernel of multiplication by $\bar x$ is $\bar x^{r-1}N$ and the kernel of multiplication by $\bar x^{r-1}$ is $\bar xN$. A first-quadrant double complex whose rows are resolutions has total complex quasi-isomorphic to the complex of the $H_0$ of the rows, so $\Tot D$ is a free resolution of $\kk_C$. It is minimal, because $d_F\otimes1$ has entries in $C'_+$ by minimality of $F$, while $\varphi_1$ and $\varphi_{r-1}$ have image in $(\Tot D)\,xC$ and $(\Tot D)\,x^{r-1}C$, which lie in $(\Tot D)\,C_+$ since $r\ge2$.

(2) By (1), $\Ext^n_C(\kk,\kk)=\Hom_C(\Tot_nD,\kk)=\bigoplus_{p+q=n}\Hom_C(D_{p,q},\kk)$, and
\[
\Hom_C(D_{p,q},\kk)\cong\Hom_{C'}(F_p,\kk)(s_q)\cong\Ext^p_{C'}(\kk,\kk)(s_q).
\]
Hence $P_C(y,s)=P_{C'}(y,s)\sum_{q\ge0}y^qs^{s_q}$, and $\sum_qy^qs^{s_q}=(1+ys)\sum_jy^{2j}s^{jr}=(1+ys)/(1-y^2s^r)$. Setting $s=1$ gives the second formula of \eqref{YAMe5Tr4}, since $(1+y)/(1-y^2)=(1-y)^{-1}$.

(3) The maps $f\otimes c\mapsto f\pi(c)$ on $D_{p,0}$, and zero on $D_{p,q}$ for $q>0$, form a chain map $\Phi\colon\Tot D\to F$ over $\pi$ lifting the identity of $\kk$. Indeed, $\Phi\varphi_1(f\otimes c)=\theta(f)\pi(xc)=0$. So $\pi^*(a)$ is represented by the cocycle $a\circ\Phi$, which is $a$ on $D_{p,0}$ and zero elsewhere, and $\pi^*$ is injective. The identity maps $D_{p,q+2}\to D_{p,q}(-r)$ form a chain map $\Theta\colon\Tot D\to\Tot D$ of homological degree $-2$, because the vertical differentials are $2$-periodic. Its component $D_{0,2}=C(-r)\to D_{0,0}(-r)=C(-r)$ is the identity, so $\Theta$ lifts the cocycle dual to the generator of $D_{0,2}$, that is, the class $\eta\in\Ext^2_C(\kk,\kk)$. By the chain-map description of the Yoneda product, $\Theta^j$ lifts $\eta^j$, and the Yoneda product $\pi^*(a)\eta^j$ is represented by $a\circ\Phi\circ\Theta^j$. This cocycle is $a$ on $D_{p,2j}$ and zero elsewhere. As $a$ runs through a basis of $\Ext^p_{C'}(\kk,\kk)$ and $j\ge0$, these cocycles are linearly independent, which proves (3).
\end{proof}

\begin{example}\label{YAMex5Tr6}
Let $C'=\kk[a]/(a^2)$ with $\bar\sigma(a)=-a$, and $r=3$, so that
\[
C=\kk\langle a,t\rangle/(a^2,\ ta+at,\ t^3),\quad H_C(s)=(1+s)(1+s+s^2).
\]
The minimal resolution of $\kk_{C'}$ has $F_p=e_pC'$ with $e_p$ of degree $p$ and $d(e_p)=e_{p-1}a$. The chain map $\theta$ of \eqref{YAMe5Tr2} is $\theta(e_pc)=(-1)^pe_p\bar\sigma(c)$, since $d\theta(e_p)=\theta(e_{p-1}a)=\theta(e_{p-1})\bar\sigma(a)$ forces the signs to alternate. So $\Tot_nD$ is free on generators $e_{p,q}$ with $p+q=n$, of internal degree $p+s_q$, and
\[
d(e_{p,q})=e_{p-1,q}\,a+(-1)^{p+pk_q}\,e_{p,q-1}\,t^{k_q},\quad k_q=\begin{cases}1&q\text{ odd},\\2&q\ge2\text{ even},\end{cases}
\]
with the convention $e_{-1,q}=e_{p,-1}=0$. For instance $d(e_{1,1})=e_{0,1}a+e_{1,0}t$, and $d^2(e_{1,1})=e_{0,0}(ta+at)=0$. Here $P_{C'}(y,s)=(1-ys)^{-1}$, and \eqref{YAMe5Tr4} gives
\[
P_C(y,s)=\frac{1+ys}{(1-ys)(1-y^2s^3)},\quad\dim\Ext^n_C(\kk,\kk)=n+1,
\]
so $\gci(C)=2$. This is the covariant ring of the dual reflection group $\mathbb Z_2\times\mathbb Z_3\cong\mathbb Z_6$ of the quantum plane $\kk_{-1}[x_1,x_2]=\kk[x_1][x_2;\sigma]$, $\sigma(x_1)=-x_1$, with $G$-degrees of $x_1$ and $x_2$ of orders two and three, where $a=\bar x_1$ and $t=\bar x_2$.
\end{example}

The resolution of Proposition \ref{YAMp3Tr7}(1) is the graded case $\delta=0$ of the resolution of McPhate for truncated Ore extensions \cite[Theorem 3.10]{McP}, which is built from the twisted product resolutions of Shepler and Witherspoon \cite[Theorem 5.12]{SW}. We need its minimality and the explicit form of Proposition \ref{YAMp3Tr7}(3) to control the growth of $E(C)$.

\begin{proof}[\textnormal{\bfseries Proof of Theorem \ref{YAMtd9Ci4}}]
For $0\le i\le n$, let $C^{(i)}$ be the covariant ring of $B^{(i)}$ for its $K_i$-grading, with $C^{(0)}=\kk$. By the proof of Theorem \ref{YAMtc3Mn8}, $K_i$ is a dual reflection group for $B^{(i)}=B^{(i-1)}[x_i;\sigma_i,\delta_i]$ and $K_{i-1}\cap\langle c_i\rangle=\{e\}$, and $c_i=e$ whenever $\delta_i\neq0$. So Proposition \ref{YAMcv1Ka7} gives
\[
C^{(i)}\cong C^{(i-1)}[\bar x_i;\bar\sigma_i]/(\bar x_i^{\,r_i}),
\]
which is $C^{(i-1)}$ when $c_i=e$. Since $C^{(n)}=C(B)$, this proves the description of $C(B)$. Let $m=\#\{i\mid c_i\neq e\}$. Applying Proposition \ref{YAMp3Tr7}(2) at each step with $r_i\ge2$, starting from $P_\kk=1$, gives
\[
P_{C(B)}(y,s)=\prod_{c_i\neq e}\frac{1+ys}{1-y^2s^{r_i}},\quad\dim\Ext^n_{C(B)}(\kk,\kk)=\binom{n+m-1}{m-1}\quad(m\ge1).
\]
If $m=0$, then every $c_i$ is trivial, so $G=\{e\}$, $C(B)=\kk$ and $\gci(C(B))=0$. Assume now that $m\ge1$. It remains to show that $\GKdim E(C(B))=m$. Write $E=E(C(B))$. Let $W\subseteq E$ be a finite-dimensional subspace containing $1$, contained in $\bigoplus_{n\le d}E^n$. Then $W^N\subseteq\bigoplus_{n\le dN}E^n$ has dimension at most $\sum_{n\le dN}(n+1)^{m-1}\le(dN+1)^m$, so $\GKdim E\le m$. For the reverse inequality, let $i_1<\cdots<i_m$ be the indices with $c_i\neq e$, and let $\pi_k^*\colon E(C^{(i_{k-1})})\to E(C^{(i_k)})$ and $\eta_k\in E^2(C^{(i_k)})$ be as in Proposition \ref{YAMp3Tr7}(3), where $C^{(i_k-1)}=C^{(i_{k-1})}$. Put $W_0=\kk$ and $W_k=\pi_k^*(W_{k-1})+\kk\eta_k$. By Proposition \ref{YAMp3Tr7}(3), $W_k^{2N}$ contains the linearly independent families $\pi_k^*(W_{k-1}^N)\eta_k^j$ for $0\le j\le N$, so
\[
\dim W_k^{2N}\ge(N+1)\dim W_{k-1}^N.
\]
By induction, $\dim W_m^{2^mN}\ge(N+1)^m$, and $\GKdim E\ge m$. So $C(B)$ is a complete intersection of growth type with $\gci(C(B))=m$.
\end{proof}

The bound $\gci(T)\le\gci(T/(\Omega))\le\gci(T)+1$ for a regular normal element $\Omega$ \cite[Theorem 1.11(a)]{KKZ2}, applied to $T=C^{(i-1)}[\bar x_i;\bar\sigma_i]$ and $\Omega=\bar x_i^{\,r_i}$, gives only inequalities at each step. The equality uses that $\bar x_i^{\,r_i}$ has degree at least two, through the periodicity class $\eta$. In the same way, Proposition \ref{YAMp3Tr7} shows that the Ext-algebra of $C(B)$ has the same bigraded Hilbert series as the Ext-algebra of the commutative complete intersection $\bigotimes_{c_i\neq e}\kk[x_i]/(x_i^{r_i})$, which is the covariant ring of the diagonal reflection group $\prod_i\mu_{r_i}$.

\begin{example}\label{YAMsh5Ge1}
In Example \ref{YAMex1Cm4}, the grading of $\kk[x_1,\dots,x_n]$ has homogeneous variables, and $C(B)=\bigotimes_i\kk[x_i]/(x_i^{r_i})$. For $\Gamma=\kk[x]/(x^r)$ with $r\ge2$, the resolution of Proposition \ref{YAMp3Tr7} with $C'=\kk$ is the periodic resolution
\[
\cdots\to\Gamma(-2r)\xrightarrow{\ x^{r-1}\ }\Gamma(-r-1)\xrightarrow{\ x\ }\Gamma(-r)\xrightarrow{\ x^{r-1}\ }\Gamma(-1)\xrightarrow{\ x\ }\Gamma\to\kk\to0,
\]
and $E(\Gamma)$ is generated by the class $\xi\in\Ext^1$ and the periodicity class $\eta\in\Ext^2$. Theorem \ref{YAMtd9Ci4} gives $\gci(C(B))=\#\{i\mid r_i\ge2\}$, which agrees with the K\"unneth formula.
\end{example}

\begin{example}\label{YAMex16Qc9}
Let $B=\mathcal O_q(\kk^n)$ with a grading as in Corollary \ref{YAMc7Qa1}, so that $G=\prod_i\langle c_i\rangle$ and $B_e$ is generated by $x_1^{r_1},\dots,x_n^{r_n}$. Then $B(B_e)_{\ge1}=\sum_iBx_i^{r_i}$, and
\[
C(B)=\mathcal O_q(\kk^n)/\big(x_1^{r_1},\dots,x_n^{r_n}\big),\quad H_{C(B)}(s)=\prod_{i=1}^n[r_i]_s,
\]
a quantum complete intersection. Each $x_i^{r_i}$ is normal, since $x_i^{r_i}x_j=q_{ij}^{r_i}x_jx_i^{r_i}$. For $0\le k<n$, the quotient $\mathcal O_q(\kk^n)/(x_1^{r_1},\dots,x_k^{r_k})$ has the basis of monomials $x_1^{a_1}\cdots x_n^{a_n}$ with $a_i<r_i$ for $i\le k$, and multiplication by $x_{k+1}^{r_{k+1}}$ sends such a monomial to a nonzero multiple of another one. So $x_1^{r_1},\dots,x_n^{r_n}$ is a regular sequence of normal elements, and $C(B)$ is a classical complete intersection in the sense of \cite[Definition 1.3]{KKZ2}. The Ext-algebras of these quantum complete intersections were determined by Bergh and Oppermann \cite[Theorem 5.3]{BO}. By Theorem \ref{YAMtd9Ci4},
\[
P_{C(B)}(y,s)=\prod_{r_i\ge2}\frac{1+ys}{1-y^2s^{r_i}},\quad\gci\big(C(B)\big)=\#\{i\mid r_i\ge2\}.
\]
\end{example}

\begin{example}\label{YAMex5Cb3}
For the algebra $B=\kk[x][y;\sigma_2][z;\sigma_3,\delta_3]$ of Example \ref{YAMx7Lp1} with its $\mathbb Z_2$-grading, we have $c_1=c_3=e$ and $c_2=\tau$. Theorem \ref{YAMtd9Ci4} gives $C(B)\cong\kk[\bar y]/(\bar y^{\,2})$, $P_{C(B)}(y,s)=(1-ys)^{-1}$ and $\gci(C(B))=1$. The step with $\delta_3\neq0$ contributes nothing.
\end{example}

\begin{example}\label{YAMex6Cd7}
Let $B$ be the algebra \eqref{YAMe72Dt1} with $q=\pm1$, graded by $G=(\mathbb Z_r\times\mathbb Z_r)\rtimes\mathbb Z_{2m}$ as in Example \ref{YAMc12Dt0}, and write $a=\bar x_1$, $b=\bar x_2$ and $c=\bar x_3$. By Theorem \ref{YAMtd9Ci4}, $C(B)$ is presented by
\[
a^r=b^r=c^{2m}=0,\quad ba=q\,ab,\quad ca=\lambda\,bc,\quad cb=\mu\,ac,
\]
it has the basis $a^ib^jc^k$ with $0\le i,j<r$ and $0\le k<2m$, and
\[
\begin{gathered}
H_{C(B)}(s)=[r]_s^2\,[2m]_s,\quad\gci\big(C(B)\big)=3,\\
P_{C(B)}(y,s)=\Big(\frac{1+ys}{1-y^2s^{r}}\Big)^2\,\frac{1+ys}{1-y^2s^{2m}}.
\end{gathered}
\]
The algebra $C(B)$ is quadratic exactly when $r=2$ and $m=1$. In that case, for the degree-lexicographic order with $c>b>a$, the leading words of the six relations are $aa,bb,cc,ba,ca,cb$, the words avoiding them are the eight words $a^ib^jc^k$ with $i,j,k\in\{0,1\}$, and $\dim C(B)=|G|=8$. So the relations form a quadratic Gr\"obner basis, $C(B)$ is Koszul \cite[Theorem 5.3]{Pri}, and for the pairing $\langle\alpha\beta,ab\rangle=1$ and the convention $E(C)\cong C^!$ for the Yoneda product its Ext-algebra is
\[
E\big(C(B)\big)=\kk\langle\alpha,\beta,\kappa\rangle/\big(\alpha\beta+q\beta\alpha,\ \kappa\alpha+\lambda^{-1}\beta\kappa,\ \kappa\beta+\mu^{-1}\alpha\kappa\big),
\]
with Hilbert series $(1-s)^{-3}$. For $r>2$ or $m>1$, the algebra $C(B)$ is not quadratic, and Theorem \ref{YAMtd9Ci4} still gives its Poincar\'e series.
\end{example}

\begin{example}\label{YAMex6Cd8}
For the family of Example \ref{YAMex7Wr3}, Theorem \ref{YAMtd9Ci4} gives
\[
C(B)\cong\Big(\bigotimes_{i=1}^n\kk[\bar x_i]/(\bar x_i^{\,r})\Big)[\bar t;\bar\sigma]/(\bar t^{\,np}),\quad\bar\sigma(\bar x_i)=\bar x_{i+1},
\]
with Hilbert series $[r]_s^n[np]_s$, and for $r\ge2$,
\[
P_{C(B)}(y,s)=\Big(\frac{1+ys}{1-y^2s^r}\Big)^n\frac{1+ys}{1-y^2s^{np}},\quad\gci\big(C(B)\big)=n+1.
\]
We make the case $n=r=2$, $p=1$ of Example \ref{YAMex13Dt6} explicit. Write $a=\bar x_1$ and $b=\bar x_2$. Then $C=C(B)$ is presented by
\[
a^2=b^2=t^2=0,\quad ba=ab,\quad ta=bt,\quad tb=at.
\]
For the degree-lexicographic order with $t>b>a$, the leading words of these six relations are $aa$, $bb$, $tt$, $ba$, $ta$, $tb$. The words avoiding them are those with strictly increasing letters, that is, $a^ib^jt^k$ with $i,j,k\in\{0,1\}$. Reducing with the six relations shows that these eight words span $C$. Since $\dim C=|G|=8$, they are linearly independent, so no further leading word occurs in the ideal of relations. Hence the six relations form a quadratic Gr\"obner basis, $C$ has a PBW basis, and $C$ is Koszul \cite[Theorem 5.3]{Pri}. Its Ext-algebra is the quadratic dual
\[
E(C)=C^!=\kk\langle\alpha,\beta,\tau\rangle/(\alpha\beta+\beta\alpha,\ \tau\alpha+\beta\tau,\ \tau\beta+\alpha\tau),
\]
since the orthogonal complement of the span of $aa,bb,tt,ba-ab,ta-bt,tb-at$ is spanned by $\alpha\beta+\beta\alpha$, $\tau\alpha+\beta\tau$ and $\tau\beta+\alpha\tau$. The algebra $C^!$ is the Ore extension $\kk_{-1}[\alpha,\beta][\tau;\vartheta]$ with $\vartheta(\alpha)=-\beta$ and $\vartheta(\beta)=-\alpha$, which is AS regular of dimension three. Hence
\[
H_C(s)=(1+s)^3,\quad H_{E(C)}(s)=\frac{1}{H_C(-s)}=\frac{1}{(1-s)^3},\quad\dim\Ext^q_C(\kk,\kk)=\binom{q+2}{2},
\]
in agreement with $P_{C(B)}(y,s)=(1+ys)^3/(1-y^2s^2)^3=(1-ys)^{-3}$. In low degrees,
\[
\begin{array}{c|cccccc}
q&0&1&2&3&4&5\\\hline
\dim\Ext^q_C(\kk,\kk)&1&3&6&10&15&21
\end{array}
\]
and $\gci(C)=3$ for this nonabelian group. In this case $C$ is also a classical complete intersection. The ideal $B(B_e)_{\ge1}$ is generated by $x_1^2$, $x_2^2$ and $t^2$, hence by
\[
\Omega_1=t^2,\quad\Omega_2=x_1^2+x_2^2,\quad\Omega_3=x_1^2-x_2^2.
\]
Here $\Omega_1$ and $\Omega_2$ are central in $B$, since $tx_1^2=x_2^2t$ and $tx_2^2=x_1^2t$, and $t\Omega_3=-\Omega_3t$, so $\Omega_3$ is normal. As left $A$-modules, $B/(\Omega_1)=A\oplus At$ with $A=\kk[x_1,x_2]$, and $\Omega_2,\Omega_3\in A$ act on both summands by left multiplication. Since $\Omega_2,\Omega_3$ is a regular sequence in $A$, the sequence $\Omega_1,\Omega_2,\Omega_3$ of normal elements is regular in $B$, and $C=B/(\Omega_1,\Omega_2,\Omega_3)$.
\end{example}

The following table records the identity components, the covariant rings and the values of $\gci(C(B))$ in the examples above. In each case $\gci(C(B))$ is the number of generators of nontrivial $G$-degree.
\begin{center}
\small
\renewcommand{\arraystretch}{1.3}
\begin{tabular}{lllll}
$B$ & $G$ & $B_e$ & $C(B)$ & $\gci(C(B))$\\\hline
$\kk[x_1,\dots,x_n]$ & $\prod_i\mathbb Z_{r_i}$ & $\kk[x_1^{r_1},\dots,x_n^{r_n}]$ & $\bigotimes_i\kk[\bar x_i]/(\bar x_i^{\,r_i})$ & $\#\{i\mid r_i\ge2\}$\\
$\mathcal O_q(\kk^n)$ & $\prod_i\mathbb Z_{r_i}$ & $\mathcal O_{q'}(\kk^n)$ & $\mathcal O_q(\kk^n)/(x_i^{r_i})$ & $\#\{i\mid r_i\ge2\}$\\
Jordan plane & $\{e\}$ & $B$ & $\kk$ & $0$\\
Example \ref{YAMx7Lp1} & $\mathbb Z_2$ & $\kk[x,w][z;\delta_e]$ & $\kk[\bar y]/(\bar y^{\,2})$ & $1$\\
Example \ref{YAMex3Od5} & $\mathbb Z_2$ & $\kk[x,w][z;\sigma,\delta]$ & $\kk[\bar y]/(\bar y^{\,2})$ & $1$\\
Example \ref{YAMc12Dt0} & $(\mathbb Z_r\times\mathbb Z_r)\rtimes\mathbb Z_{2m}$ & skew polynomial & Example \ref{YAMex6Cd7} & $3$\\
$\kk[x_1,\dots,x_n][t;\sigma]$ & $\mathbb Z_r^n\rtimes\mathbb Z_{np}$, $r\ge2$ & $\kk[x_1^r,\dots,x_n^r,t^{np}]$ & Example \ref{YAMex6Cd8} & $n+1$
\end{tabular}
\end{center}

\section{McKay quivers}\label{YAMs8Mq0}

In this section we describe the smash products $B\#\kk^G$ and $C(B)\#\kk^G$ by quivers with relations, we compute their matrix Hilbert series, and we invert the Jacobian. This proves Theorem \ref{YAMtd5Qv1}.

Let $G$ be a finite group and $B$ a $G$-graded algebra. Let $\{p_g\}_{g\in G}$ be the basis of $\kk^G$ dual to $G$, so that $p_gp_h=\delta_{g,h}p_g$. The smash product $\Lambda=B\#\kk^G$ is the vector space $B\otimes\kk^G$ with the multiplication
\begin{equation}\label{YAMe82Sm1}
(a\#p_g)(b\#p_h)=a\,b_{gh^{-1}}\#p_h\quad(a,b\in B,\ g,h\in G),
\end{equation}
where $b_{gh^{-1}}$ is the component of $b$ of degree $gh^{-1}$ \cite[Proposition 1.4]{CoMo}. Both ways of computing the product of $a\#p_g$, $b\#p_h$ and $c\#p_k$ give the element $a\,b_{gh^{-1}}c_{hk^{-1}}\#p_k$, and the unit is $\sum_g1\#p_g$. The category of graded right $B$-modules is isomorphic to the category of right $\Lambda$-modules \cite[Theorem 2.2]{CoMo}. Put
\[
e_g=1\#p_g\quad(g\in G).
\]
These are orthogonal idempotents with $\sum_ge_g=1$, and \eqref{YAMe82Sm1} gives
\begin{equation}\label{YAMe84Co3}
e_g\Lambda e_h=B_{gh^{-1}}\#p_h\cong B_{gh^{-1}}\quad(g,h\in G).
\end{equation}
For $a,b\in B_e$ we have $(a\#p_h)(b\#p_h)=ab\#p_h$, so every corner $e_h\Lambda e_h$ is isomorphic to $B_e$ as an algebra. In particular, if $B$ is a noetherian AS regular domain generated in degree one with a faithful $G$-grading refining the $\mathbb N$-grading, then $G$ is a dual reflection group for $B$ if and only if one, or equivalently every, corner $e_h\Lambda e_h$ is AS regular. If $B$ is $\mathbb N$-graded and the $G$-grading refines the $\mathbb N$-grading, then $\Lambda$ is $\mathbb N$-graded by $\deg(b\#p_h)=\deg b$, and its \emph{matrix Hilbert series} is
\begin{equation}\label{YAMe85Mh4}
H_\Lambda(s)=\big(H_{e_g\Lambda e_h}(s)\big)_{g,h\in G}=\big(H_{B_{gh^{-1}}}(s)\big)_{g,h\in G}.
\end{equation}
If $G$ is abelian, then $\kk^G\cong\kk G^\vee$ as Hopf algebras, and $\Lambda$ is the skew group algebra $B\#\kk G^\vee$ for the action of the dual group described in Section \ref{YAMs1Pr0}.

Suppose that $B$ is generated by $B_1$, and let $x_1,\dots,x_n$ be a basis of $B_1$ consisting of $G$-homogeneous elements, with $x_i$ of degree $c_i$. Such a basis exists when the $G$-grading refines the $\mathbb N$-grading. Let $Q_G$ be the quiver with vertex set $G$ and with arrows
\[
x_i^{(h)}\colon h\longrightarrow c_ih\quad(h\in G,\ 1\le i\le n).
\]
This is the McKay quiver of the $\kk^G$-module $B_1$ in the sense of \cite{Cra22}, and it is the Cayley quiver of $G$ with respect to $c_1,\dots,c_n$. It is connected when the $c_i$ generate $G$. We compose paths from right to left, and we send the trivial path at $g$ to $e_g$ and the arrow $x_i^{(h)}$ to $x_i\#p_h\in e_{c_ih}\Lambda e_h$. By \eqref{YAMe82Sm1},
\begin{equation}\label{YAMe87Pa6}
(x_{i_k}\#p_{h_k})\cdots(x_{i_1}\#p_{h_1})=x_{i_k}\cdots x_{i_1}\#p_{h_1}\quad\text{if }h_{j+1}=c_{i_j}h_j\text{ for all }j,
\end{equation}
and the product is zero otherwise. For a $G$-homogeneous element $\rho$ of the tensor algebra $T(B_1)$ and $h\in G$, we write $\rho^{(h)}$ for the linear combination of paths starting at $h$ which corresponds to $\rho\#p_h$ under \eqref{YAMe87Pa6}.

Let $B=T(V)/I$ with $V=B_1$, where $I$ is generated by elements $\rho_1,\dots,\rho_k$ that are homogeneous for the $\mathbb N$-grading and for the $G$-grading. Then the assignment above induces an isomorphism of graded algebras
\begin{equation}\label{YAMe89Qr2}
\Lambda=B\#\kk^G\cong\kk Q_G/I_G,
\end{equation}
where $I_G$ is the ideal generated by the elements $\rho_j^{(h)}$ for $1\le j\le k$ and $h\in G$. This is the case of a quiver with one vertex of the presentation of smash products of group weighted bound quivers by Asashiba \cite[Theorem 1.18]{Asa}, the smash product of a $G$-graded algebra with $G$ being identified with $B\#\kk^G$ by \cite[Proposition 3.7]{CM2}. It is also an instance of the theory of coverings of quivers with relations \cite{Gre}, and of Crawford's description of smash products of Koszul AS regular algebras with semisimple Hopf algebras by quivers with relations \cite[Theorem 1.2]{Cra22}. For iterated Ore extensions it takes the following explicit form.

\begin{corollary}\label{YAMc13Cq6}
Let $B$ be the iterated Ore extension \eqref{YAMq8Hi1} with a $G$-grading in which every $x_i$ is homogeneous. Then $\Lambda=B\#\kk^G$ is the path algebra of the Cayley quiver $Q_G$ of $G$ with respect to $c_1,\dots,c_n$, modulo the ideal generated by
\begin{equation}\label{YAMe88Lr7}
x_i^{(c_jh)}x_j^{(h)}-\sigma_i(x_j)^{(c_ih)}x_i^{(h)}-\delta_i(x_j)^{(h)}\quad(1\le j<i\le n,\ h\in G),
\end{equation}
where $\sigma_i(x_j)\in B_1$ and $\delta_i(x_j)$ is replaced by a $G$-homogeneous lift in $V\otimes V$ of degree $c_ic_j$. All three terms of \eqref{YAMe88Lr7} are paths from $h$ to $c_ic_jh$. If moreover $G$ is a dual reflection group for $B$ and $(B_e)_1=0$, then the relations are the binomial squares
\[
x_i^{(c_jh)}x_j^{(h)}=\lambda_{ij}\,x_{\pi_i(j)}^{(c_ih)}x_i^{(h)}\quad(1\le j<i\le n,\ h\in G),
\]
with $\lambda_{ij}$ and $\pi_i$ as in Theorem \ref{YAMtc3Mn8}.
\end{corollary}

\begin{proof}
The relations \eqref{YAMe13Pr4} are homogeneous for both gradings, since $\sigma_i(x_j)$ has degree $c_ic_jc_i^{-1}$ and $\delta_i(x_j)$ has degree $c_ic_j$ by \eqref{YAMg6Tw1}. Since $T(V)\to B$ is a surjection of $G$-graded spaces, homogeneous lifts exist, and \eqref{YAMe89Qr2} applies. The term $\sigma_i(x_j)^{(c_ih)}x_i^{(h)}$ ends at $c_ic_jc_i^{-1}\cdot c_ih=c_ic_jh$. The second statement follows from the binomial presentation of Theorem \ref{YAMtc3Mn8}, since $c_{\pi_i(j)}c_i=c_ic_j$.
\end{proof}

We compute the matrix Hilbert series. For $c\in G$, let $\Pi_c$ be the permutation matrix of left multiplication by $c$, with rows and columns indexed by $G$,
\begin{equation}\label{YAMe90Lc9}
(\Pi_c)_{g,h}=\begin{cases}1&g=ch,\\0&\text{otherwise},\end{cases}\quad\Pi_c\Pi_{c'}=\Pi_{cc'},\quad\Pi_c^{|c|}=I.
\end{equation}
Left multiplication by $c$ has $|G|/|c|$ orbits on $G$, the right cosets of $\langle c\rangle$, each of length $|c|$. Hence
\begin{equation}\label{YAMe91Dt0}
\det(I-s\Pi_c)=(1-s^{|c|})^{|G|/|c|},\quad(I-s\Pi_c)^{-1}=\frac{1}{1-s^{|c|}}\sum_{k=0}^{|c|-1}s^k\Pi_c^k.
\end{equation}

\begin{proposition}\label{YAMp11Mh2}
Let $B$ be the iterated Ore extension \eqref{YAMq8Hi1} with a $G$-grading in which every $x_i$ is homogeneous, and let $\Lambda=B\#\kk^G$. Then
\begin{equation}\label{YAMe92Hm1}
H_\Lambda(s)=(I-s\Pi_{c_1})^{-1}(I-s\Pi_{c_2})^{-1}\cdots(I-s\Pi_{c_n})^{-1},
\end{equation}
so that $H_\Lambda(s)^{-1}=(I-s\Pi_{c_n})\cdots(I-s\Pi_{c_1})$ and $\det H_\Lambda(s)^{-1}=\prod_{i=1}^n(1-s^{r_i})^{|G|/r_i}$.
\end{proposition}

\begin{proof}
The monomials $x_1^{a_1}\cdots x_n^{a_n}$ with $a_i\ge0$ form a basis of $B$, and such a monomial has $\mathbb N$-degree $a_1+\cdots+a_n$ and $G$-degree $c_1^{a_1}\cdots c_n^{a_n}$. By \eqref{YAMe85Mh4} and \eqref{YAMe90Lc9}, the $(g,h)$ entry of $H_\Lambda(s)$ is
\[
\sum_{c_1^{a_1}\cdots c_n^{a_n}=gh^{-1}}s^{a_1+\cdots+a_n}=\sum_{a_1,\dots,a_n\ge0}s^{a_1+\cdots+a_n}\big(\Pi_{c_1}^{a_1}\cdots\Pi_{c_n}^{a_n}\big)_{g,h}.
\]
Summing the geometric series in each $a_i$ gives \eqref{YAMe92Hm1}, and the determinant follows from \eqref{YAMe91Dt0}.
\end{proof}

Thus $H_\Lambda(s)^{-1}$ is a matrix polynomial of degree $n$, whose factors $I-s\Pi_{c_i}$ involve the adjacency matrices $\Pi_{c_i}$ of the arrows of each type, one for each Ore step. For a dual reflection group, the matrix Hilbert series also factors through the identity component.

\begin{proposition}\label{YAMc14Pm3}
In the setting of Proposition \ref{YAMp11Mh2}, suppose that $G$ is a dual reflection group for $B$, and put $P_G(s)=\big(s^{\deg f_{gh^{-1}}}\big)_{g,h\in G}$. Then
\begin{equation}\label{YAMe93Pg2}
\begin{aligned}
H_\Lambda(s)&=H_{B_e}(s)\,P_G(s),\\
P_G(s)&=\prod_{i=1}^n\Big(\sum_{k=0}^{r_i-1}s^k\Pi_{c_i}^k\Big),\\
\det P_G(s)&=\prod_{i=1}^n(1-s^{r_i})^{|G|(r_i-1)/r_i}.
\end{aligned}
\end{equation}
\end{proposition}

\begin{proof}
By \eqref{YAMe3Lf4}, $H_{B_{gh^{-1}}}(s)=s^{\deg f_{gh^{-1}}}H_{B_e}(s)$, which gives the first identity by \eqref{YAMe85Mh4}. By \eqref{YAMe91Dt0} and Proposition \ref{YAMp11Mh2},
\[
H_\Lambda(s)=\prod_{i=1}^n\frac{1}{1-s^{r_i}}\cdot\prod_{i=1}^n\Big(\sum_{k=0}^{r_i-1}s^k\Pi_{c_i}^k\Big),
\]
where the scalar factors have been moved to the front. The first factor is $H_{B_e}(s)$ by Theorem \ref{YAMtc3Mn8}, and it is an invertible power series, which gives the factorization of $P_G$. The determinant follows from \eqref{YAMe91Dt0}, since $\det\big(\sum_{k<r}s^k\Pi_c^k\big)=(1-s^r)^{|G|}\det(I-s\Pi_c)^{-1}=(1-s^r)^{|G|-|G|/r}$ for $r=|c|$.
\end{proof}

The factorization of $P_G$ also follows directly from the normal form of Theorem \ref{YAMtc3Mn8}(3). The $(g,h)$ entry of the product is $\sum s^{a_1+\cdots+a_n}$ over the tuples with $0\le a_i<r_i$ and $c_1^{a_1}\cdots c_n^{a_n}=gh^{-1}$, and by uniqueness of the normal form there is exactly one such tuple, with $a_1+\cdots+a_n=\deg f_{gh^{-1}}$ by \eqref{YAMe78Nf6}. Each column of $P_G(s)$ sums to $p_G(s)=\prod_i[r_i]_s$, so \eqref{YAMe93Pg2} is the matrix form of the factorization of $p_G$.

Let $G$ be a dual reflection group for $B$. The ideal $J=B(B_e)_{\ge1}$ is two-sided by \eqref{YAMt4Kv0}, and it is $G$-graded. So $J\#\kk^G$ is an ideal of $\Lambda$, and
\[
\Lambda/(J\#\kk^G)\cong C(B)\#\kk^G.
\]

\begin{proposition}\label{YAMp9Cq3}
Let $B$ be a noetherian AS regular domain generated in degree one, $G$ a dual reflection group for $B$, and $\bar\Lambda=C(B)\#\kk^G$.
\begin{enumerate}
\item The algebra $\bar\Lambda$ has dimension $|G|^2$, and for all $g,h\in G$, the space $e_g\bar\Lambda e_h=\kk\,\bar f_{gh^{-1}}\#p_h$ is one-dimensional, concentrated in degree $\deg f_{gh^{-1}}$. The matrix Hilbert series of $\bar\Lambda$ is $P_G(s)$.
\item The algebra $\bar\Lambda$ is Frobenius. The socle of the left projective module $\bar\Lambda e_h$ is $e_{m_Gh}\bar\Lambda e_h$, which is isomorphic to the simple module at the vertex $m_Gh$.
\end{enumerate}
\end{proposition}

\begin{proof}
(1) By \eqref{YAMt4Kv0}, $C(B)_{gh^{-1}}=\kk\bar f_{gh^{-1}}$, and \eqref{YAMe84Co3} for $C(B)$ gives the claim.

(2) Since $C(B)$ is Frobenius by \eqref{YAMt4Kv0}, so is $\bar\Lambda$ by \cite[Corollary 2.2]{DNN}. The projective module $\bar\Lambda e_h=\bigoplus_ge_g\bar\Lambda e_h$ is indecomposable since $e_h\bar\Lambda e_h=\kk e_h$. Its component $e_{m_Gh}\bar\Lambda e_h=\kk\,\bar f_{m_G}\#p_h$ has the top degree $N_G$. Since $\bar\Lambda$ is finite-dimensional and $\bar\Lambda_0=\bigoplus_g\kk e_g$ is semisimple, the radical of $\bar\Lambda$ is $\bar\Lambda_{\ge1}$. So $e_{m_Gh}\bar\Lambda e_h$ is annihilated by the radical and lies in the socle. The socle of an indecomposable projective module over a Frobenius algebra is simple, so it equals $e_{m_Gh}\bar\Lambda e_h$.
\end{proof}

For the iterated Ore extensions of Theorem \ref{YAMtc3Mn8}, the algebra $B_e$ is generated by the elements $x_i^{r_i}$ by Proposition \ref{YAMp7Ie3}, so $J$ is the two-sided ideal generated by $x_1^{r_1},\dots,x_n^{r_n}$. Then Corollary \ref{YAMc13Cq6} shows that $\bar\Lambda$ is the quotient of $\kk Q_G/I_G$ by the paths
\[
x_i^{(c_i^{r_i-1}h)}\cdots x_i^{(c_ih)}x_i^{(h)}\quad(1\le i\le n,\ h\in G)
\]
of length $r_i$, where $r_i=1$ when $\delta_i\neq0$.

We invert the Jacobian. Let $G$ be a dual reflection group for $B$. By \eqref{YAMt4Kv0}, $f_{m_G}$ is normal, so there is a graded automorphism $\nu$ of $B$ with $f_{m_G}b=\nu(b)f_{m_G}$, and comparing $G$-degrees in the domain $B$ gives $\nu(B_g)=B_{m_Ggm_G^{-1}}$. The set $\{f_{m_G}^k\}_{k\ge0}$ is a two-sided Ore set, and the localization
\[
L=B[f_{m_G}^{-1}]
\]
is a $G$-graded domain containing $B$, with $f_{m_G}^{-1}$ homogeneous of degree $m_G^{-1}$.

\begin{proposition}\label{YAMlo4Pe9}
Let $B$ be a noetherian AS regular domain generated in degree one, $G$ a dual reflection group for $B$, and $L=B[f_{m_G}^{-1}]$. Let $o$ be the order of $m_G$ and $F=f_{m_G}^{\,o}$.
\begin{enumerate}
\item Every $f_g$ is a unit of $L$, and $L_g=f_gL_e=L_ef_g$ for all $g\in G$. So $L$ is a crossed product $L_e*G$, and $L_gL_h=L_{gh}$ for all $g,h\in G$.
\item The element $F$ is a normal element of $B_e$, and $L_e=B_e[F^{-1}]$.
\item The elements
\[
u_{g,h}=f_gf_h^{-1}\#p_h\in e_g(L\#\kk^G)e_h\quad(g,h\in G)
\]
form a complete set of matrix units of $L\#\kk^G$, and $L\#\kk^G\cong M_{|G|}(L_e)$.
\end{enumerate}
\end{proposition}

\begin{proof}
(1) By \eqref{YAMe27Fm1}, the element $f_g$ has the right inverse $\mu_g^{-1}f_{g^{-1}m_G}f_{m_G}^{-1}$ in $L$, which is homogeneous of degree $g^{-1}$. In a domain a one-sided inverse is two-sided, since $ab=1$ implies $a(ba-1)=0$. So $f_g$ is a unit with $f_g^{-1}\in L_{g^{-1}}$, and for $b\in L_g$ we have $b=f_g(f_g^{-1}b)$ with $f_g^{-1}b\in L_e$, and similarly on the other side. Then $L_{gh}=f_gf_hL_e\subseteq L_gL_h\subseteq L_{gh}$.

(2) The element $F$ has degree $m_G^{\,o}=e$, so $F\in B_e$, and $FB_e=\nu^o(B_e)F=B_eF$, since $\nu$ preserves $B_e$. An element of $L_e$ has the form $bf_{m_G}^{-k}$ with $b\in B$, and comparing degrees gives $b\in B_{m_G^k}$. Choose $l$ with $lo\ge k$. Then $bf_{m_G}^{-k}=\big(bf_{m_G}^{\,lo-k}\big)F^{-l}$ with $bf_{m_G}^{\,lo-k}\in B_{m_G^{lo}}=B_e$. Hence $L_e=B_e[F^{-1}]$.

(3) For $g,h,h',k\in G$, the element $f_{h'}f_k^{-1}$ is homogeneous of degree $h'k^{-1}$, so \eqref{YAMe82Sm1} gives
\[
u_{g,h}u_{h',k}=f_gf_h^{-1}\,\big(f_{h'}f_k^{-1}\big)_{hk^{-1}}\#p_k=\begin{cases}f_gf_k^{-1}\#p_k=u_{g,k}&h=h',\\0&h\neq h'.\end{cases}
\]
Moreover $\sum_gu_{g,g}=\sum_g1\#p_g=1$. A ring with a complete set of $N\times N$ matrix units is isomorphic to the matrix ring of size $N$ over the corner of $u_{e,e}$. Here $u_{e,e}=e_e$ since $f_e=1$, and $e_e(L\#\kk^G)e_e\cong L_e$ by \eqref{YAMe84Co3} for $L$.
\end{proof}

The identity \eqref{YAMe27Fm1} is due to Kirkman and Zhang \cite[Lemma 2.11(2)]{KZ}. The Morita equivalence between $L\#\kk^G$ and $L_e$ also follows from the fact that $L$ is strongly graded \cite[Theorem 2.12]{CoMo}, and the matrix units of Proposition \ref{YAMlo4Pe9}(3) make it an explicit isomorphism. In Example \ref{YAMex1Cm4}, $L=\kk[x_i^{\pm1}\ (r_i>1),\ x_j\ (r_j=1)]$, the algebra $\Lambda$ is the skew group algebra of the diagonal reflection group $\prod_i\mu_{r_i}$, and Proposition \ref{YAMlo4Pe9} recovers the Morita equivalence between $L\#\kk G^\vee$ and $L^{G^\vee}=L_e$ on the complement of the reflecting hyperplanes $x_i=0$ with $r_i>1$, on which the group acts freely \cite[Theorem 1.5]{Ste}.

\begin{proof}[\textnormal{\bfseries Proof of Theorem \ref{YAMtd5Qv1}}]
Part (1) is Corollary \ref{YAMc13Cq6} and Proposition \ref{YAMp11Mh2}. Part (2) is Propositions \ref{YAMc14Pm3} and \ref{YAMp9Cq3}, and part (3) is Proposition \ref{YAMlo4Pe9}.
\end{proof}

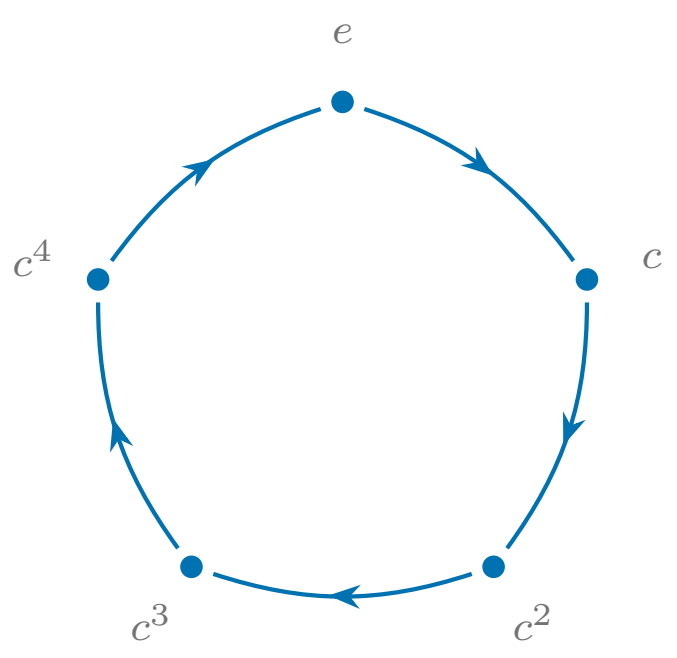
\begin{figure}[H]
\centering
\begin{tikzpicture}
\foreach \k in {0,...,4}{
  \node[vx,fill=mtone] (c\k) at ({90-72*\k}:1.7) {};
}
\node[vname] at ({90}:2.15) {$e$};
\node[vname] at ({90-72}:2.15) {$c$};
\node[vname] at ({90-144}:2.15) {$c^2$};
\node[vname] at ({90-216}:2.15) {$c^3$};
\node[vname] at ({90-288}:2.15) {$c^4$};
\foreach \k in {0,...,4}{
  \pgfmathtruncatemacro{\l}{mod(\k+1,5)}
  \draw[midarrow,thick,mtone,shorten <=2pt,shorten >=2pt] (c\k) to[bend left=18] (c\l);
}
\end{tikzpicture}
\caption{The McKay quiver of $\kk[x]$ graded by $\mathbb Z_5$.}
\label{YAMfg1Cy2}
\end{figure}

\begin{example}\label{YAMex19Qd9}
For $\kk[x]$ graded by $\mathbb Z_r=\langle c\rangle$ with $x$ of degree $c$, the quiver $Q_G$ is the cyclic quiver with $r$ vertices and no relations, so $\kk[x]\#\kk^{\mathbb Z_r}$ is its path algebra (Figure \ref{YAMfg1Cy2}), and \eqref{YAMe92Hm1} reads $H_\Lambda(s)^{-1}=I-s\Pi_c$, where $\Pi_c$ is the adjacency matrix of the cyclic quiver. The matrix $P_G(s)$ is the circulant matrix with first column $(1,s,\dots,s^{r-1})$, and \eqref{YAMe93Pg2} reduces to the classical formula $\det P_G(s)=(1-s^r)^{r-1}$. The algebra $\bar\Lambda$ is the path algebra of the cyclic quiver modulo all paths of length $r$, a self-injective Nakayama algebra, and the socle of its projective at $c^j$ sits at $c^{j+r-1}$, in agreement with $m_G=c^{r-1}$. Moreover $L=\kk[x^{\pm1}]$ and $L_e=\kk[x^{\pm r}]$ for $r\ge2$.
\end{example}

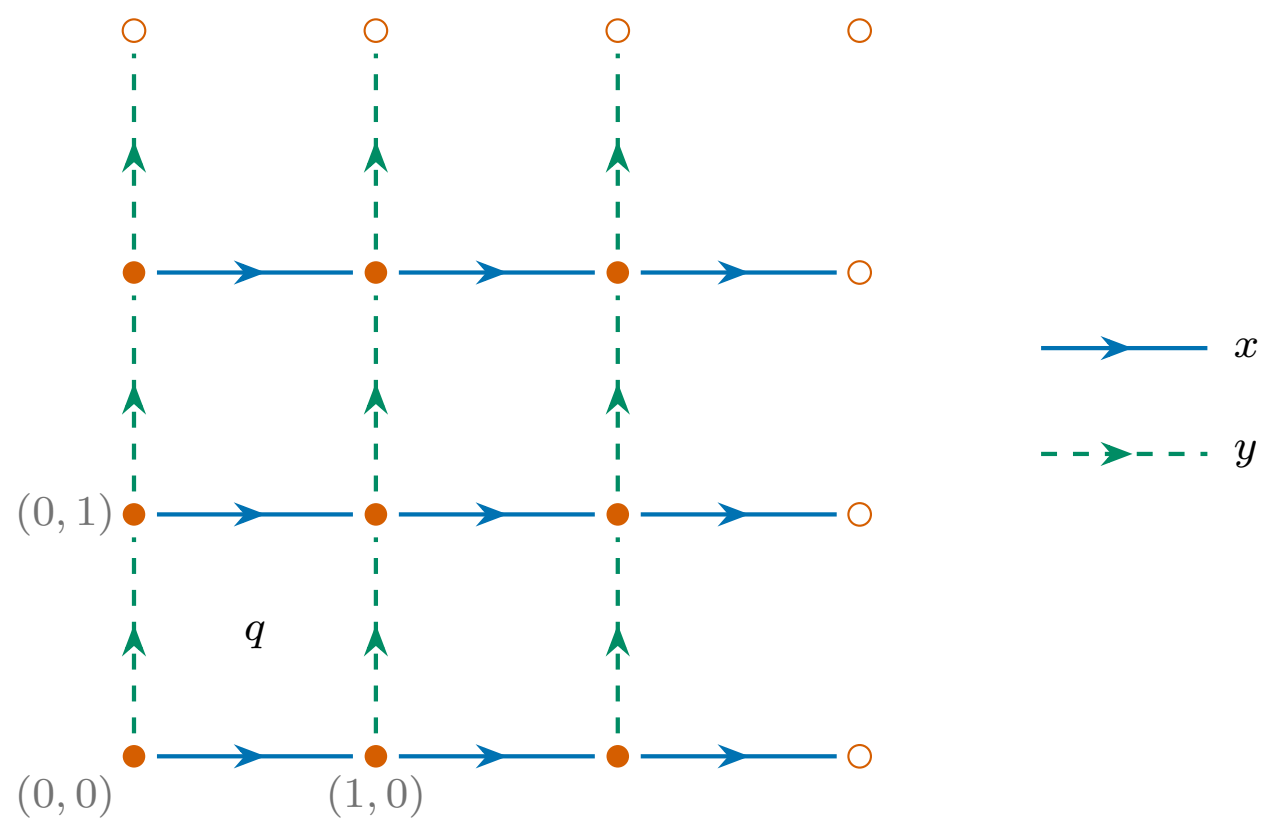
\begin{figure}[H]
\centering
\begin{tikzpicture}
\foreach \a in {0,1,2,3}{\foreach \b in {0,1,2,3}{
  \pgfmathtruncatemacro{\bd}{(\a==3)||(\b==3)}
  \ifnum\bd=1 \node[vx,draw=mttwo,fill=white] (n\a\b) at (1.6*\a,1.6*\b) {}; \else \node[vx,fill=mttwo] (n\a\b) at (1.6*\a,1.6*\b) {}; \fi
}}
\foreach \a in {0,1,2}{\foreach \b in {0,1,2}{
  \pgfmathtruncatemacro{\c}{\a+1}
  \draw[midarrow,thick,mtone,shorten <=2pt,shorten >=2pt] (n\a\b) -- (n\c\b);
}}
\foreach \a in {0,1,2}{\foreach \b in {0,1,2}{
  \pgfmathtruncatemacro{\c}{\b+1}
  \draw[midarrow,thick,mtthree,dashed,shorten <=2pt,shorten >=2pt] (n\a\b) -- (n\a\c);
}}
\node[font=\footnotesize] at (0.8,0.8) {$q$};
\node[vname,below left] at (n00) {$(0,0)$};
\node[vname,below] at (n10) {$(1,0)$};
\node[vname,left] at (n01) {$(0,1)$};
\begin{scope}[shift={(6.0,2.0)}]
\draw[midarrow,thick,mtone] (0,0.7) -- (1.1,0.7); \node[font=\footnotesize,right] at (1.15,0.7) {$x$};
\draw[midarrow,thick,mtthree,dashed] (0,0) -- (1.1,0); \node[font=\footnotesize,right] at (1.15,0) {$y$};
\end{scope}
\end{tikzpicture}
\caption{The McKay quiver of $\kk_q[x,y]$ graded by $\mathbb Z_3\times\mathbb Z_3$. The white vertices are identified with the vertices on the opposite side.}
\label{YAMfg3Tq4}
\end{figure}

\begin{example}\label{YAMex19Qd8}
For the quantum plane $\kk_q[x,y]$ with $yx=qxy$, graded by $\mathbb Z_{r_1}\times\mathbb Z_{r_2}$ with $\deg x=(1,0)$ and $\deg y=(0,1)$, the quiver $Q_G$ is the toroidal grid with vertices $(a,b)$ and arrows
\[
x^{(a,b)}\colon(a,b)\to(a+1,b),\quad y^{(a,b)}\colon(a,b)\to(a,b+1),
\]
and $I_G$ is generated by the $q$-commutative squares
\[
y^{(a+1,b)}x^{(a,b)}=q\,x^{(a,b+1)}y^{(a,b)}\quad\big((a,b)\in\mathbb Z_{r_1}\times\mathbb Z_{r_2}\big),
\]
see Figure \ref{YAMfg3Tq4}. Crawford describes the same toroidal pictures for coactions on two-dimensional AS regular algebras \cite[Section 4]{Cra1}, using the presentation of \cite[Theorem 2.6]{Cra1}, with arrows $g\to a^{-1}g$ and paths composed from left to right. For $g=(a,b)$ we have $f_g=x^ay^b$ and $g^{-1}m_G=(r_1-1-a,r_2-1-b)$. Since $y^bx^c=q^{bc}x^cy^b$,
\[
f_gf_{g^{-1}m_G}=x^ay^b\,x^{r_1-1-a}y^{r_2-1-b}=q^{b(r_1-1-a)}\,x^{r_1-1}y^{r_2-1},
\]
so $\mu_g=q^{b(r_1-1-a)}$ in \eqref{YAMe27Fm1}. If $r_1,r_2\ge2$, then $L$ is the quantum torus $\kk_q[x^{\pm1},y^{\pm1}]$, with identity component $L_e=\kk_{q^{r_1r_2}}[x^{\pm r_1},y^{\pm r_2}]$.
\end{example}

\begin{example}\label{YAMex19Qd6}
We describe the McKay quiver of the $\mathbb Z_2$-grading of Example \ref{YAMx7Lp1}, where $\delta_3\neq0$. Here $B=\kk[x][y;\sigma_2][z;\sigma_3,\delta_3]$ with relations \eqref{YAMe12Br6}, $\mathbb Z_2=\{e,\tau\}$, $c_1=c_3=e$ and $c_2=\tau$. The quiver $Q_G$ has the two vertices $e$ and $\tau$, a loop $x^{(h)}$ and a loop $z^{(h)}$ at each vertex $h$, and arrows $y^{(e)}\colon e\to\tau$ and $y^{(\tau)}\colon\tau\to e$. By Corollary \ref{YAMc13Cq6}, the ideal $I_G$ is generated by the six relations
\begin{gather*}
y^{(h)}x^{(h)}=-x^{(\tau h)}y^{(h)},\quad z^{(\tau h)}y^{(h)}=-y^{(h)}z^{(h)}-2x^{(\tau h)}y^{(h)},\\
z^{(h)}x^{(h)}=x^{(h)}z^{(h)}+x^{(h)}x^{(h)}-y^{(\tau h)}y^{(h)},
\end{gather*}
for $h\in\{e,\tau\}$. The last family contains the lift of $\delta_3(x)=x^2-y^2$, whose two terms are paths from $h$ to $h$. By Propositions \ref{YAMp11Mh2} and \ref{YAMc14Pm3},
\[
H_\Lambda(s)=\frac{1}{(1-s)^2}(I-s\Pi_\tau)^{-1}=\frac{1}{(1-s)^2(1-s^2)}\begin{pmatrix}1&s\\s&1\end{pmatrix},\quad P_G(s)=I+s\Pi_\tau,
\]
which agrees with $H_{B_e}(s)=(1-s)^{-2}(1-s^2)^{-1}$. The loops $x^{(h)}$ and $z^{(h)}$ correspond to generators of trivial degree, so they are killed in $\bar\Lambda$, and $\bar\Lambda=C(B)\#\kk^{\mathbb Z_2}$ is the path algebra of the cyclic quiver with two vertices modulo the paths of length two. The Jacobian $f_{m_G}=y$ has degree $\tau$ of order two, $F=y^2=w$, and Proposition \ref{YAMlo4Pe9} gives
\[
L_e=B_e[w^{-1}]=\kk[x,w^{\pm1}][z;\delta_e],\quad L\#\kk^{\mathbb Z_2}\cong M_2(L_e),
\]
with $\delta_e$ as in \eqref{YAMe12Br7}.
\end{example}

\begin{figure}[H]
\centering
\begin{tikzpicture}[
  xa/.style={midarrow,thick,mtone},
  xb/.style={midarrow,thick,mttwo},
  tt/.style={midarrow,thick,mtthree}]
\node[vx,fill=mtstab] (b00) at (0,0) {};
\node[vx,fill=mtstab] (b10) at (3.4,0) {};
\node[vx,fill=mtstab] (b01) at (1.6,1.35) {};
\node[vx,fill=mtstab] (b11) at (5.0,1.35) {};
\node[vx,fill=mtslice] (t00) at (0,3.6) {};
\node[vx,fill=mtslice] (t10) at (1.6,4.95) {};
\node[vx,fill=mtslice] (t01) at (3.4,3.6) {};
\node[vx,fill=mtslice] (t11) at (5.0,4.95) {};
\node[vname,below left] at (b00) {$e$};
\node[vname,below right] at (b10) {$\varepsilon_1$};
\node[vname,left=2pt] at (b01) {$\varepsilon_2$};
\node[vname,right] at (b11) {$\varepsilon_1\varepsilon_2$};
\node[vname,left] at (t00) {$\gamma$};
\node[vname,above left] at (t10) {$\varepsilon_1\gamma$};
\node[vname,right=2pt] at (t01) {$\varepsilon_2\gamma$};
\node[vname,above right] at (t11) {$\varepsilon_1\varepsilon_2\gamma$};
\foreach \u/\v/\s in {b00/b10/xa,b01/b11/xa,b00/b01/xb,b10/b11/xb,t00/t10/xa,t01/t11/xa,t00/t01/xb,t10/t11/xb,b00/t00/tt,b10/t01/tt,b01/t10/tt,b11/t11/tt}{
  \draw[\s,shorten <=2pt,shorten >=2pt] (\u) to[bend left=12] (\v);
  \draw[\s,shorten <=2pt,shorten >=2pt] (\v) to[bend left=12] (\u);
}
\begin{scope}[shift={(7.8,1.4)}]
\draw[xa] (0,1.4) -- (1.1,1.4); \node[font=\footnotesize,right] at (1.15,1.4) {$x_1$};
\draw[xb] (0,0.7) -- (1.1,0.7); \node[font=\footnotesize,right] at (1.15,0.7) {$x_2$};
\draw[tt] (0,0) -- (1.1,0); \node[font=\footnotesize,right] at (1.15,0) {$t$};
\end{scope}
\end{tikzpicture}
\caption{The McKay quiver of the dihedral group of order eight in Example \ref{YAMex19Qd7}.}
\label{YAMfg2Dh3}
\end{figure}
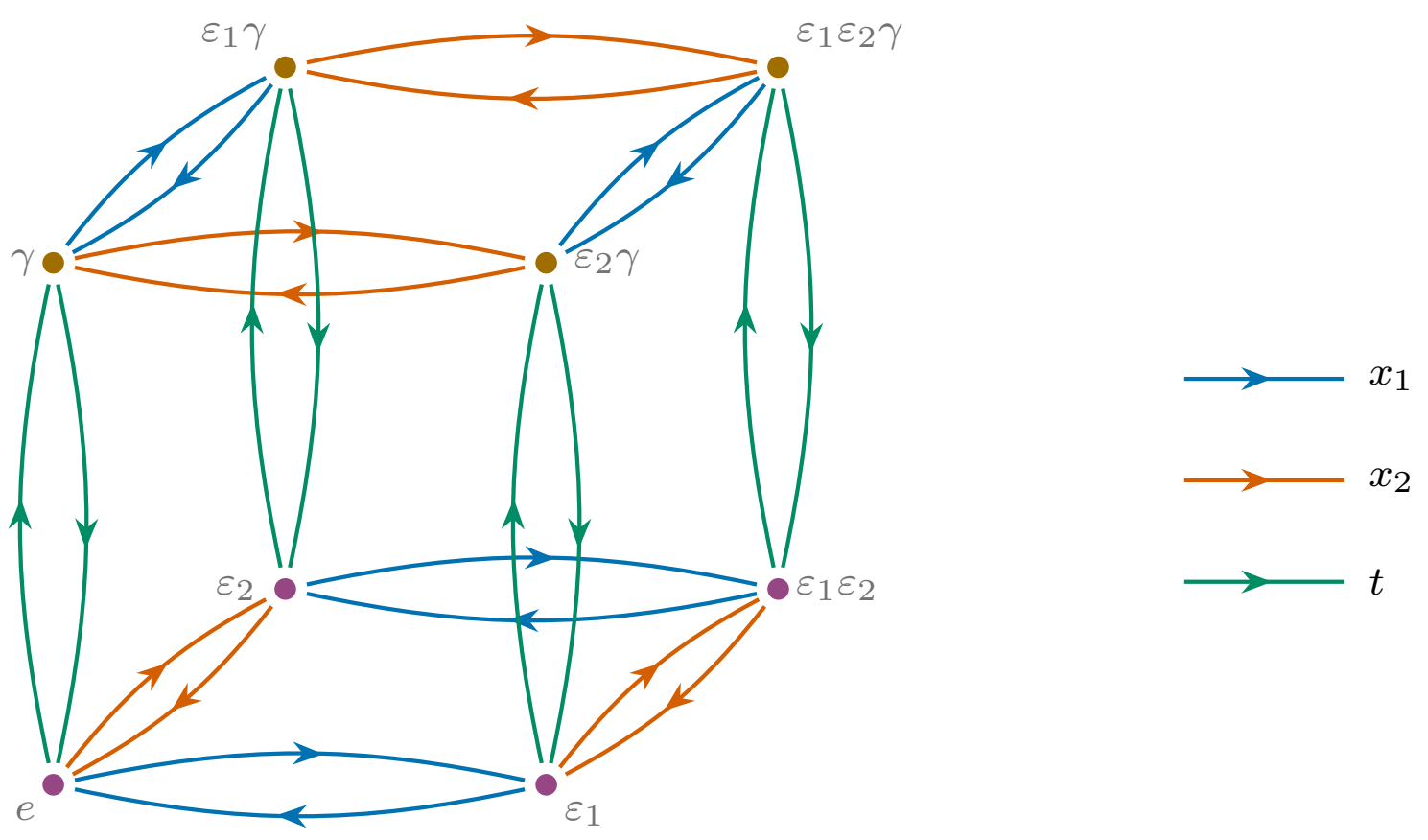

\begin{example}\label{YAMex19Qd7}
For the dihedral group of Example \ref{YAMex13Dt6}, the quiver $Q_G$ has eight vertices and three arrows $x_1^{(h)}\colon h\to\varepsilon_1h$, $x_2^{(h)}\colon h\to\varepsilon_2h$ and $t^{(h)}\colon h\to\gamma h$ at every vertex. By Corollary \ref{YAMc13Cq6} and \eqref{YAMe66Wr9}, the relations are the twenty-four commutative squares
\begin{gather*}
x_2^{(\varepsilon_1h)}x_1^{(h)}=x_1^{(\varepsilon_2h)}x_2^{(h)},\\
t^{(\varepsilon_1h)}x_1^{(h)}=x_2^{(\gamma h)}t^{(h)},\quad t^{(\varepsilon_2h)}x_2^{(h)}=x_1^{(\gamma h)}t^{(h)}\quad(h\in G).
\end{gather*}
Here $\gamma\varepsilon_1=\varepsilon_2\gamma$ and $\gamma\varepsilon_2=\varepsilon_1\gamma$, so both sides of each square have the same endpoints. Since $\varepsilon_1$, $\varepsilon_2$ and $\gamma$ are involutions, the arrows come in opposite pairs, and $Q_G$ is a cube (Figure \ref{YAMfg2Dh3}). The $t$-arrows join $h$ and $\gamma h$, so they join $e$ to $\gamma$, $\varepsilon_1$ to $\varepsilon_2\gamma$, $\varepsilon_2$ to $\varepsilon_1\gamma$ and $\varepsilon_1\varepsilon_2$ to $\varepsilon_1\varepsilon_2\gamma$. On the bottom face $\langle\varepsilon_1\rangle\times\langle\varepsilon_2\rangle$, the $x_1$-arrows join $e$ to $\varepsilon_1$, while on the top face they join $\gamma$ to $\varepsilon_1\gamma$, which lies above $\varepsilon_2$. So the arrows of type $x_1$ and $x_2$ exchange their directions between the two faces. This twist is the conjugation by $\gamma$, and it is the reason why $G$ is not abelian.

Order $G$ as $e,\varepsilon_1,\varepsilon_2,\varepsilon_1\varepsilon_2,\gamma,\varepsilon_1\gamma,\varepsilon_2\gamma,\varepsilon_1\varepsilon_2\gamma$. The table of Example \ref{YAMex13Dt6} gives
\[
P_G(s)=\begin{pmatrix}
1&s&s&s^2&s&s^2&s^2&s^3\\
s&1&s^2&s&s^2&s^3&s&s^2\\
s&s^2&1&s&s^2&s&s^3&s^2\\
s^2&s&s&1&s^3&s^2&s^2&s\\
s&s^2&s^2&s^3&1&s&s&s^2\\
s^2&s^3&s&s^2&s&1&s^2&s\\
s^2&s&s^3&s^2&s&s^2&1&s\\
s^3&s^2&s^2&s&s^2&s&s&1
\end{pmatrix},
\]
which factors as
\[
P_G(s)=(I+s\Pi_{\varepsilon_1})(I+s\Pi_{\varepsilon_2})(I+s\Pi_\gamma),
\]
with $\det P_G(s)=(1-s^2)^{12}$ and $H_\Lambda(s)=(1-s^2)^{-3}P_G(s)$. The algebra $\bar\Lambda=C(B)\#\kk^G$ is $\kk Q_G$ modulo the twenty-four squares and the twenty-four paths
\[
x_1^{(\varepsilon_1h)}x_1^{(h)},\quad x_2^{(\varepsilon_2h)}x_2^{(h)},\quad t^{(\gamma h)}t^{(h)}\quad(h\in G).
\]
There are $72$ paths of length two, and $\dim\bar\Lambda_n=8,24,24,8$ for $n=0,1,2,3$, so $\dim\bar\Lambda=64=|G|^2$. Every indecomposable projective has Loewy layers of dimensions $1,3,3,1$, and by Proposition \ref{YAMp9Cq3} the socle of the projective at $h$ is the simple module at $\varepsilon_1\varepsilon_2\gamma h$, the vertex opposite to $h$ in the cube. The mass element $m_G=\varepsilon_1\varepsilon_2\gamma$ has order two, $F=(x_1x_2t)^2=x_1^2x_2^2t^2$, and by Proposition \ref{YAMlo4Pe9}
\[
\begin{aligned}
L_e&=\kk[x_1^2,x_2^2,t^2]\big[(x_1^2x_2^2t^2)^{-1}\big]=\kk[x_1^{\pm2},x_2^{\pm2},t^{\pm2}],\\
L\#\kk^G&\cong M_8(L_e).
\end{aligned}
\]
So the localization $L\#\kk^G$ of $\Lambda$ is Morita equivalent to $L_e$, the coordinate ring of the complement of the three coordinate planes in $\operatorname{Spec}B_e$.
\end{example}

\section{Noncommutative quotients}\label{YAMs9Gq1}
In this section we interpret the extension $B_e\subseteq B$, the localization on which it becomes Hopf--Galois, and the point scheme of the dihedral example in the language of noncommutative algebraic geometry.

Artin and Zhang associate a noncommutative projective scheme to a noetherian connected graded algebra $B$ through the quotient category $\operatorname{qgr}B$ \cite{AZ}. We use the phrase \emph{noncommutative affine space} only as a geometric analogy for an AS regular algebra $B$, and we write $X/G$ informally for the object represented by the identity component $B_e$. A $G$-grading of $B$ is an action of the Hopf algebra $\kk^G$, whose invariant subalgebra is $B_e$. For $G$ abelian it is an action of the character group $G^\vee$ by graded automorphisms (Section \ref{YAMs1Pr0}). In this algebraic sense Theorem \ref{YAMtc3Mn8} is a Chevalley--Shephard--Todd theorem for iterated Ore extensions. It characterizes when $B_e$ is AS regular, and it identifies $B_e$ with the iterated Ore extension of Proposition \ref{YAMp7Ie3} on the generators $x_i^{r_i}$.

Let $G$ be a dual reflection group for $B$. By \eqref{YAMq9Fb2}, the extension $B_e\subseteq B$ is free of rank $|G|$ on both sides, with bases $\{f_g\}$. The augmentation of $B_e$ gives
\[
B\otimes_{B_e}\kk=B/B(B_e)_{\ge1}=C(B),
\]
an algebra of dimension $|G|$, since $B(B_e)_{\ge1}$ is a two-sided ideal by \eqref{YAMt4Kv0}. It is the algebraic analogue of the fibre over the vertex of a commutative reflection quotient. For iterated Ore extensions, Theorem \ref{YAMtd9Ci4} shows that $C(B)$ is a complete intersection of growth type whose Ext-algebra has the same bigraded Hilbert series as the covariant ring of the diagonal reflection group $\prod_i\mu_{r_i}$, also when $G$ is not abelian.

We identify a localization on which the extension becomes Hopf--Galois. A $G$-graded algebra $R$ is a \emph{$\kk G$-Galois extension} of $R_e$ if the canonical map
\[
\beta\colon R\otimes_{R_e}R\longrightarrow R\otimes\kk G,\quad a\otimes b\longmapsto\sum_{g\in G}ab_g\otimes g,
\]
is bijective. This holds if and only if $R$ is strongly graded, that is, $R_gR_h=R_{gh}$ for all $g,h\in G$ \cite[Theorem 8.1.7]{Mon}. It is a noncommutative analogue of the principal bundle condition. When $G$ is abelian and $R$ is commutative, it is the classical torsor condition for the action of the character group $G^\vee$. For nonabelian $G$ it is understood in the Hopf--Galois sense.

\begin{proposition}\label{YAMp12Gl3}
Let $B$ be a noetherian AS regular domain generated in degree one, $G$ a dual reflection group for $B$, and $L=B[f_{m_G}^{-1}]$.
\begin{enumerate}
\item The algebra $L$ is a $\kk G$-Galois extension of $L_e$.
\item Suppose that $B$ is the iterated Ore extension \eqref{YAMq8Hi1} with every $x_i$ homogeneous. A ring homomorphism $\varphi\colon B\to R$ inverts $f_{m_G}$ if and only if it inverts every $x_i$ with $c_i\neq e$. In particular $L$ is the localization of $B$ at the generators of nontrivial $G$-degree, and $x_i^{-1}\in L_{c_i^{-1}}$ for every such $i$.
\end{enumerate}
\end{proposition}

\begin{proof}
(1) By Proposition \ref{YAMlo4Pe9}(1), $L$ is strongly graded, so $L$ is a $\kk G$-Galois extension of $L_e$ by \cite[Theorem 8.1.7]{Mon}.

(2) If $c_i\neq e$, then $r_i\ge2$, and the normal form of $c_i$ in Theorem \ref{YAMtc3Mn8}(3) is $c_i$ itself, so $x_i=f_{c_i}$ by Theorem \ref{YAMtc3Mn8}. By Theorem \ref{YAMtc3Mn8}, $f_{m_G}=x_1^{r_1-1}\cdots x_n^{r_n-1}$, where only the factors with $c_i\neq e$ occur. If $\varphi$ inverts these $x_i$, it inverts their product $f_{m_G}$. Conversely, if $\varphi$ inverts $f_{m_G}$, then $\varphi$ extends to the Ore localization $L$, and $\varphi(x_i)=\varphi(f_{c_i})$ is invertible since $f_{c_i}$ is a unit of $L$ by Proposition \ref{YAMlo4Pe9}(1). The last statement follows since $x_i^{-1}=f_{c_i}^{-1}\in L_{c_i^{-1}}$.
\end{proof}

The generators $x_i$ with $c_i\neq e$ need not be normal individually. In Example \ref{YAMex13Dt6}, for instance, $tx_1=x_2t$ does not lie in $x_1B$, so we do not interpret their zero loci as noncommutative hypersurfaces. For the polynomial ring of Example \ref{YAMex1Cm4}, the grading is the diagonal action of $\prod_i\mu_{r_i}$, the algebra $L$ is the coordinate ring of the complement of the reflecting hyperplanes $x_i=0$ with $r_i>1$, and the action is free on this complement \cite[Theorem 1.5]{Ste}. This agrees with the reflection arrangement $f_{m_G^{-1}}$ of Kirkman and Zhang \cite[Definition 2.1]{KZ}, which in Examples \ref{YAMex1Cm4} and \ref{YAMex13Dt6} is the product of the generators of nontrivial $G$-degree.

\begin{example}\label{YAMex20Wl4}
For the family of Example \ref{YAMex7Wr3} with $r\ge2$, every generator has nontrivial $G$-degree. By Proposition \ref{YAMp12Gl3}(2), $L$ is the skew Laurent polynomial ring
\[
L=\kk[x_1^{\pm1},\dots,x_n^{\pm1}][t^{\pm1};\sigma],\quad L_e=\kk[x_1^{\pm r},\dots,x_n^{\pm r},t^{\pm np}].
\]
We have $L_e=B_e[F^{-1}]$ by Proposition \ref{YAMlo4Pe9}(2). The element $D=\varepsilon_1\cdots\varepsilon_n$ is central of order $r$, and $m_G=D^{r-1}\gamma^{np-1}=D^{-1}\gamma^{-1}$ has order $o=\operatorname{lcm}(r,np)$. Since $x_1\cdots x_n$ commutes with $t$,
\[
F=f_{m_G}^{\,o}=(x_1\cdots x_n)^{o(r-1)}\,t^{o(np-1)},
\]
a monomial in $x_1^r,\dots,x_n^r,t^{np}$ in which every variable occurs. Inverting $F$ in the polynomial ring $B_e=\kk[x_1^r,\dots,x_n^r,t^{np}]$ therefore inverts each of these generators. Thus $L_e$ is the coordinate ring of a commutative algebraic torus, while $L$ is a noncommutative skew Laurent polynomial ring. The extension $L_e\subseteq L$ is Hopf--Galois for the coaction of $\kk G$ with $G=\mathbb Z_r^n\rtimes\mathbb Z_{np}$, and $L\#\kk^G\cong M_{npr^n}(L_e)$.
\end{example}

Let $B=T(V)/(\mathcal R)$ be a quadratic algebra with $V=B_1$. A \emph{point module} is a graded right $B$-module $M=\bigoplus_{j\ge0}M_j$, generated by $M_0$, with $\dim M_j=1$ for all $j$ \cite{ATV}. Choosing generators $m_j$ of $M_j$, we have $m_jx=p_j(x)m_{j+1}$ for $x\in V$, where $p_j$ is a nonzero linear form on $V$. Since the relations of $B$ are quadratic, the sequence $(p_j)$ of points of $\mathbb P(V^*)$ defines a point module if and only if
\begin{equation}\label{YAMe97Pt2}
\sum_{a,b}\lambda_{ab}\,p_j(x_a)\,p_{j+1}(x_b)=0\quad\text{for every relation }\sum_{a,b}\lambda_{ab}x_ax_b\in\mathcal R\text{ and every }j\ge0.
\end{equation}

\begin{example}\label{YAMex21Pt5}
Let $B=\kk[x_1,x_2][t;\sigma]$ with $\sigma(x_1)=x_2$ and $\sigma(x_2)=x_1$, the algebra of the dihedral group in Example \ref{YAMex13Dt6}, with relations $x_2x_1=x_1x_2$, $tx_1=x_2t$ and $tx_2=x_1t$. Write $p=(a_1:a_2:a_3)$ and $q=(b_1:b_2:b_3)$ for points of $\mathbb P(V^*)=\mathbb P^2$, with $a_1=p(x_1)$, $a_2=p(x_2)$ and $a_3=p(t)$. By \eqref{YAMe97Pt2}, consecutive points $p=p_j$ and $q=p_{j+1}$ of a point module satisfy $M(p)q=0$, where
\[
M(p)=\begin{pmatrix}a_2&-a_1&0\\a_3&0&-a_2\\0&a_3&-a_1\end{pmatrix},\quad\det M(p)=a_3\,(a_2^2-a_1^2).
\]
So $p_0$ lies in the cubic
\[
E=\{a_3=0\}\cup\{a_1=a_2\}\cup\{a_1=-a_2\}\subseteq\mathbb P^2,
\]
a triangle of lines. The matrix $M(p)$ has rank exactly two at every point of $E$. If $a_3\neq0$, the minor formed by the second and third rows and the first and second columns equals $a_3^2$. If $a_3=0$, the first row is nonzero and is independent of at least one of the other two rows. So the kernel of $M(p)$ is one-dimensional on $E$, and $q$ is uniquely determined by $p$. Since $\det M(p)$ is a product of three distinct linear forms, $E$ is reduced. Solving $M(p)q=0$ on each line gives an automorphism $\sigma_E$ of $E$ with $q=\sigma_E(p)$, namely
\[
\sigma_E(a_1:a_2:0)=(a_1:a_2:0),\quad\sigma_E(a:a:c)=(a:a:c),\quad\sigma_E(a:-a:c)=(a:-a:-c).
\]
On the line $a_1=-a_2$, for instance, the equations read $-ab_1-ab_2=0$, $cb_1+ab_3=0$ and $cb_2-ab_3=0$, whose solution is $q=(a:-a:-c)$ for $a\neq0$, and $q=(0:0:1)$ for $a=0$. The two formulas agree at the vertices $(1:-1:0)$ and $(0:0:1)$, which are fixed. Hence the point modules of $B$ are parametrized by $E$, and the point module of $p\in E$ is given by $p_j=\sigma_E^j(p)$.

The three components of $E$ are the zero loci of the normal linear elements $t$, $x_1-x_2$ and $x_1+x_2$, since $t(x_1\pm x_2)=\pm(x_1\pm x_2)t$. On the component $a_1=-a_2$, the restriction of $\sigma_E$ coincides with the projective action of the automorphism $t\mapsto-t$, $x_i\mapsto x_i$, attached to the character of $G$ that is trivial on $\varepsilon_1,\varepsilon_2$ and takes the value $-1$ on $\gamma$. This automorphism does not induce $\sigma_E$ on all of $E$. In the basis $u=x_1+x_2$, $v=x_1-x_2$, $t$, the relations read $uv=vu$, $tu=ut$ and $tv=-vt$, so $B$ is a skew polynomial ring and $E$ is its classical point scheme \cite{ATV}. The $D_8$-grading is not diagonal in this basis, since $\gamma$ exchanges $x_1$ and $x_2$, and its identity component $B_e=\kk[x_1^2,x_2^2,t^2]$ is a commutative polynomial ring. The computation agrees with the description of point schemes of quadratic Ore extensions of twists of polynomial rings by Lozano and Vancliff \cite{LoVa}.
\end{example}

\section{Homological determinants}\label{YAMs7Hd0}

In this section we compute homological determinants of graded Ore extensions and of identity components of dual reflection groups, using trace series.

Let $R$ be a connected graded algebra and $\theta$ a graded automorphism of $R$. The \emph{trace series} of $\theta$ is
\[
\Tr_R(\theta,s)=\sum_{n\ge0}\tr\big(\theta|_{R_n}\big)\,s^n.
\]
For $\theta=\mathrm{id}$ this is the Hilbert series $H_R(s)$. We use the homological determinant $\hdet_R(\theta)\in\kk^\times$ of J{\o}rgensen and Zhang \cite{JZ} through the following property. Let $R$ be a noetherian AS regular algebra of dimension $d$ with Gorenstein parameter $\ell$, and $\theta$ a graded automorphism of $R$. Then $\theta$ is rational in the sense of \cite[1.3]{JZ} by \cite[4.2]{JZ}, so $\Tr_R(\theta,s)$ is a rational function of $s$, and its Laurent expansion in $s^{-1}$ has the form \cite[Lemma 2.6]{JZ}
\begin{equation}\label{YAMr7Ye4}
\Tr_R(\theta,s)=(-1)^d\hdet_R(\theta)^{-1}s^{-\ell}+\sum_{j>\ell}a_js^{-j}.
\end{equation}
Since a rational function has a unique Laurent expansion in $s^{-1}$, \eqref{YAMr7Ye4} determines $\ell$ and $\hdet_R(\theta)$.

\begin{proposition}\label{YAMhd2Ow5}
Let $S=R[t;\sigma,\delta]$ be a graded Ore extension with $\deg t=m$, with $R$ noetherian and AS regular, and let $\theta$ be a graded automorphism of $S$ with $\theta(R)=R$. Then $\theta(t)=\lambda t+a$ with $\lambda\in\kk^\times$ and $a\in R_m$, and
\begin{equation}\label{YAMe48Hd4}
\Tr_S(\theta,s)=\frac{\Tr_R(\theta|_R,s)}{1-\lambda s^m},\quad\hdet_S(\theta)=\lambda\,\hdet_R(\theta|_R).
\end{equation}
\end{proposition}

\begin{proof}
By Corollary \ref{YAMc1Rg5}, $S$ is a noetherian AS regular algebra of dimension $d+1$ with Gorenstein parameter $\ell_R+m$, where $d=\dim R$. Since $S_m=R_m\oplus\kk t$, we can write $\theta(t)=\lambda t+a$ with $a\in R_m$, and $\lambda\neq0$ since otherwise $\theta(S)\subseteq R$. Consider the filtration of $S$ by the graded subspaces
\[
F_jS=\bigoplus_{i\le j}Rt^i.
\]
By \eqref{YAMe10Tn1}, $(F_iS)(F_jS)\subseteq F_{i+j}S$. Hence $(\lambda t+a)^i\in\lambda^it^i+F_{i-1}S$, and $\theta(bt^i)=\theta(b)(\lambda t+a)^i\in\lambda^i\theta(b)t^i+F_{i-1}S$ for $b\in R$. So $\theta$ preserves the filtration and acts on $F_iS/F_{i-1}S\cong R(-mi)$ as $\lambda^i\theta|_R$. Since $S_n\subseteq F_jS$ for $mj\ge n$, we get
\[
\Tr_S(\theta,s)=\sum_{i\ge0}\lambda^is^{mi}\,\Tr_R(\theta|_R,s)=\frac{\Tr_R(\theta|_R,s)}{1-\lambda s^m}
\]
as power series, hence as rational functions. The Laurent expansion in $s^{-1}$ is multiplicative. Using \eqref{YAMr7Ye4} for $R$ and
\[
\frac1{1-\lambda s^m}=-\lambda^{-1}s^{-m}\big(1+\lambda^{-1}s^{-m}+\lambda^{-2}s^{-2m}+\cdots\big),
\]
the leading term of the right side is $(-1)^{d+1}\big(\lambda\hdet_R(\theta|_R)\big)^{-1}s^{-(\ell_R+m)}$. Comparing with \eqref{YAMr7Ye4} for $S$ gives \eqref{YAMe48Hd4}.
\end{proof}

When $\delta=0$ and $a=0$, the element $t$ is normal with $\theta(t)=\lambda t$, and Proposition \ref{YAMhd2Ow5} is a case of the normal-element formula of Jing and Zhang \cite[Proposition 2.4]{JiZ}, see also \cite[Lemma 1.6]{KKZdu}. For general $\delta$ and $a$, neither the $\sigma$-derivation nor the translation part $a$ contributes to the homological determinant. Nakayama automorphisms of graded Ore extensions were computed by Liu, Wang and Wu \cite[Theorem 3.3]{LWW}, by Shen and Guo for Koszul AS regular algebras \cite{SG}, and by Liu and Ma for Ore extensions of polynomial algebras \cite{LM}. Combined with the theorem of J{\o}rgensen and Zhang that invariants of automorphisms with trivial homological determinant are AS Gorenstein \cite[Theorem 3.3]{JZ}, Proposition \ref{YAMhd2Ow5} gives the following criterion.

\begin{corollary}\label{YAMc3Ig7}
Let $S=R[t;\sigma,\delta]$ be as in Proposition \ref{YAMhd2Ow5}, and let $\Theta$ be a finite group of graded automorphisms of $S$ with $\theta(R)=R$ for every $\theta\in\Theta$. For $\theta\in\Theta$ write $\theta(t)=\lambda_\theta t+a_\theta$. If $\hdet_R(\theta|_R)=\lambda_\theta^{-1}$ for every $\theta\in\Theta$, then $S^\Theta$ is AS Gorenstein.
\end{corollary}

\begin{proof}
By \eqref{YAMe48Hd4}, $\hdet_S(\theta)=1$ for every $\theta\in\Theta$, and \cite[Theorem 3.3]{JZ} applies.
\end{proof}

For example, let $S=\kk_q[x,y]=\kk[x][y;\sigma]$ and let $\zeta\in\kk^\times$ be a root of unity of order $N$. The automorphism $\theta(x)=\zeta x$, $\theta(y)=\zeta^{-1}y$ generates a group of order $N$, with $\theta(\kk[x])=\kk[x]$, $\lambda_\theta=\zeta^{-1}$ and $\hdet_{\kk[x]}(\theta|_{\kk[x]})=\zeta=\lambda_\theta^{-1}$. So $\hdet_S(\theta)=1$, and $S^{\langle\theta\rangle}$ is AS Gorenstein. For $q=1$ this is the invariant ring $\kk[x^N,xy,y^N]$ of the cyclic subgroup of $SL_2$ of order $N$.

\begin{example}\label{YAMex15Jt8}
Proposition \ref{YAMhd2Ow5} applies when $\delta\neq0$ and $a\neq0$. Let $S=\kk[x][y;\delta]$ be the Jordan plane, and let $\theta$ be given by
\[
\theta(x)=\lambda x,\quad\theta(y)=\lambda y+\mu x\quad(\lambda\in\kk^\times,\ \mu\in\kk).
\]
It preserves the relation, since $\theta(y)\theta(x)-\theta(x)\theta(y)-\theta(x)^2=\lambda^2(yx-xy-x^2)$. Here $a=\mu x$, and $\theta$ acts on $\kk[x]_n=\kk x^n$ by $\lambda^n$. So
\[
\Tr_S(\theta,s)=\frac{1}{(1-\lambda s)^2}=\lambda^{-2}s^{-2}+2\lambda^{-3}s^{-3}+\cdots,
\]
and $\hdet_S(\theta)=\lambda^2=\lambda\,\hdet_{\kk[x]}(\theta|_{\kk[x]})$, as predicted by \eqref{YAMe48Hd4}. The term $\mu x$ does not contribute.
\end{example}

We relate the homological determinants of an algebra and of the identity component of a dual reflection group.

\begin{proposition}\label{YAMhd7Iu1}
Let $A$ be a noetherian AS regular domain generated in degree one, $G$ a dual reflection group for $A$, and $\theta$ a graded automorphism of $A$ with $\theta(A_g)=A_g$ for every $g\in G$. Write $\theta(f_g)=\lambda_gf_g$ with $\lambda_g\in\kk^\times$. Then
\begin{equation}\label{YAMe51Pt7}
\Tr_A(\theta,s)=p_\theta(s)\,\Tr_{A_e}(\theta|_{A_e},s),\quad p_\theta(s)=\sum_{g\in G}\lambda_gs^{\deg f_g},
\end{equation}
and
\begin{equation}\label{YAMe50Hr6}
\hdet_{A_e}(\theta|_{A_e})=\lambda_{m_G}\,\hdet_A(\theta),\quad\ell_{A_e}=\ell_A+N_G.
\end{equation}
\end{proposition}

\begin{proof}
The scalars $\lambda_g$ exist by Section \ref{YAMs1Pr0}. We have $\theta(A_e)=A_e$ and $\theta(f_gb)=\lambda_gf_g\theta(b)$ for $b\in A_e$. Since $A=\bigoplus_gf_gA_e$ and left multiplication by $f_g$ is injective, this gives \eqref{YAMe51Pt7}. The algebra $A_e$ is AS regular by definition and noetherian by Section \ref{YAMs1Pr0}, and $\dim A_e=\dim A$ by \cite[Lemma 3.4(1)]{KKZ}, since $A$ is free of finite rank over $A_e$ on both sides. So \eqref{YAMr7Ye4} applies to $\theta$ and to $\theta|_{A_e}$, with the same sign $(-1)^d$. By \eqref{YAMt3Kk6}, the leading term of $p_\theta$ is $\lambda_{m_G}s^{N_G}$. Comparing the leading terms of the Laurent expansions in $s^{-1}$ of both sides of \eqref{YAMe51Pt7} gives
\[
(-1)^d\hdet_A(\theta)^{-1}s^{-\ell_A}=(-1)^d\lambda_{m_G}\hdet_{A_e}(\theta|_{A_e})^{-1}s^{N_G-\ell_{A_e}},
\]
which is \eqref{YAMe50Hr6}.
\end{proof}

The identity $\ell_{A_e}=\ell_A+N_G$ is \cite[Lemma 3.4(3)]{KKZ}, since the Frobenius algebra $C(A)$ has Gorenstein parameter $-N_G$.

\begin{example}\label{YAMex6Hq8}
In Example \ref{YAMex8Qv2}, let $\theta(x)=\alpha x$ and $\theta(y)=\beta y$. Then $\lambda_{m_G}=\lambda_\tau=\alpha$, and
\[
\Tr_A(\theta,s)=\frac{1}{(1-\alpha s)(1-\beta s)},\quad\Tr_{A_e}(\theta|_{A_e},s)=\frac{1}{(1-\alpha^2s^2)(1-\beta s)},\quad p_\theta(s)=1+\alpha s.
\]
The identity \eqref{YAMe51Pt7} holds since $(1+\alpha s)(1-\alpha s)=1-\alpha^2s^2$. The leading terms of the expansions in $s^{-1}$ are $(\alpha\beta)^{-1}s^{-2}$ and $(\alpha^2\beta)^{-1}s^{-3}$, so $\hdet_A(\theta)=\alpha\beta$ and $\hdet_{A_e}(\theta|_{A_e})=\alpha^2\beta$, as predicted by \eqref{YAMe50Hr6}. Also $\ell_{A_e}=3=\ell_A+\deg f_\tau$.
\end{example}

For iterated Ore extensions, Propositions \ref{YAMhd2Ow5} and \ref{YAMhd7Iu1} give closed formulas.

\begin{corollary}\label{YAMc11Hi5}
Let $B$ be the iterated Ore extension \eqref{YAMq8Hi1}, and let $\theta$ be a graded automorphism of $B$ with $\theta(B^{(i)})=B^{(i)}$ for every $i$. Then $\theta(x_i)=\lambda_ix_i+a_i$ with $\lambda_i\in\kk^\times$ and $a_i\in(B^{(i-1)})_1$, and
\begin{equation}\label{YAMe79Hi7}
\hdet_B(\theta)=\lambda_1\lambda_2\cdots\lambda_n.
\end{equation}
If moreover $G$ is a dual reflection group for $B$, with every $x_i$ homogeneous, and $\theta(B_g)=B_g$ for every $g\in G$, then
\begin{equation}\label{YAMe80Hi8}
\theta(f_{m_G})=\lambda_1^{r_1-1}\cdots\lambda_n^{r_n-1}f_{m_G},\quad\hdet_{B_e}(\theta|_{B_e})=\lambda_1^{r_1}\lambda_2^{r_2}\cdots\lambda_n^{r_n}.
\end{equation}
\end{corollary}

\begin{proof}
Since $(B^{(i)})_1=(B^{(i-1)})_1\oplus\kk x_i$ and $\theta$ preserves $B^{(i)}$ and $B^{(i-1)}$, we have $\theta(x_i)=\lambda_ix_i+a_i$ with $a_i\in(B^{(i-1)})_1$, and $\lambda_i\neq0$ since $\theta(B^{(i)})=B^{(i)}$. By Proposition \ref{YAMhd2Ow5} applied to $B^{(i)}=B^{(i-1)}[x_i;\sigma_i,\delta_i]$, we get $\hdet_{B^{(i)}}(\theta)=\lambda_i\hdet_{B^{(i-1)}}(\theta)$, and \eqref{YAMe79Hi7} follows by induction, since $\hdet_\kk(\theta)=1$.

For the second statement, put $g_i=x_1^{r_1-1}\cdots x_i^{r_i-1}$, so that $f_{m_G}=g_n$ by Theorem \ref{YAMtc3Mn8}, and $\theta(g_n)=\lambda_{m_G}g_n$. Expanding $(\lambda_nx_n+a_n)^{r_n-1}=\lambda_n^{r_n-1}x_n^{r_n-1}+(\text{lower powers of }x_n)$ and comparing the coefficients of $x_n^{r_n-1}$ in the free left $B^{(n-1)}$-module $B$, we get $\theta(g_{n-1})\lambda_n^{r_n-1}=\lambda_{m_G}g_{n-1}$. So $g_{n-1}$ is an eigenvector of $\theta$. Since $g_{n-1}=f_{m_{K_{n-1}}}$ for the dual reflection group $K_{n-1}$ of $B^{(n-1)}$, the same argument and induction give $\lambda_{m_G}=\prod_i\lambda_i^{r_i-1}$. Now Proposition \ref{YAMhd7Iu1} gives $\hdet_{B_e}(\theta|_{B_e})=\lambda_{m_G}\hdet_B(\theta)$, and \eqref{YAMe79Hi7} gives \eqref{YAMe80Hi8}.
\end{proof}

The formula \eqref{YAMe80Hi8} can also be obtained from Proposition \ref{YAMp7Ie3}. The automorphism $\theta$ preserves every $B^{(i)}_e$, and $\theta(x_i^{r_i})=(\lambda_ix_i+a_i)^{r_i}$ has leading coefficient $\lambda_i^{r_i}$ in the variable $y_i=x_i^{r_i}$ of the iterated Ore extension \eqref{YAMe77Ie5}, so Proposition \ref{YAMhd2Ow5}, applied to each step of \eqref{YAMe77Ie5}, gives the same value.

\begin{example}\label{YAMex7Hw9}
In Example \ref{YAMex7Wr3}, let $\mu,\lambda\in\kk^\times$. The assignment $\theta(x_i)=\mu x_i$ and $\theta(t)=\lambda t$ defines a graded automorphism of $B$, since $\theta(tx_i-x_{i+1}t)=\lambda\mu(tx_i-x_{i+1}t)$. It preserves every $B^{(i)}$ and every component $B_g$, since it multiplies each monomial by a scalar. By Corollary \ref{YAMc11Hi5},
\[
\hdet_B(\theta)=\lambda\mu^n,\quad\theta(f_{m_G})=\mu^{n(r-1)}\lambda^{np-1}f_{m_G},\quad\hdet_{B_e}(\theta|_{B_e})=\mu^{nr}\lambda^{np}.
\]
This agrees with the direct computation in the polynomial ring $B_e=\kk[x_1^r,\dots,x_n^r,t^{np}]$, on which $\theta$ acts by $\mu^r$ on each $x_i^r$ and by $\lambda^{np}$ on $t^{np}$. If $\theta$ has finite order and $\lambda\mu^n=1$, then $B^{\langle\theta\rangle}$ is AS Gorenstein by Corollary \ref{YAMc3Ig7}. For $n=r=2$ and $p=1$, the table of Example \ref{YAMex13Dt6} gives $p_\theta(s)=(1+\mu s)^2(1+\lambda s)$, and
\[
\Tr_B(\theta,s)=\frac{1}{(1-\mu s)^2(1-\lambda s)},\quad\Tr_{B_e}(\theta|_{B_e},s)=\frac{1}{(1-\mu^2s^2)^2(1-\lambda^2s^2)}.
\]
The leading terms of the expansions in $s^{-1}$ are $-(\mu^2\lambda)^{-1}s^{-3}$ and $-(\mu^4\lambda^2)^{-1}s^{-6}$, so $\hdet_B(\theta)=\lambda\mu^2$ and $\hdet_{B_e}(\theta|_{B_e})=\lambda^2\mu^4=\lambda_{m_G}\hdet_B(\theta)$ with $\lambda_{m_G}=\mu^2\lambda$, and $\ell_{B_e}=6$.
\end{example}

\bibliographystyle{amsplain}
\bibliography{dgoe}

\end{document}